\documentclass{article}
\usepackage{subfig}
\usepackage[section] {placeins}
\usepackage{amsthm}
\usepackage{amssymb,amsmath}
\usepackage{anysize}
\usepackage{psfrag}
\newtheorem{thm}{Theorem}[section]
\newtheorem{def.}{Definition}[section]

\numberwithin{table}{section}
\begin{document}

\title{Removing Colors 2k, 2k-1, and k}
\author{Pedro Lopes\\
        Department of Mathematics\\
        Instituto Superior T\'ecnico\\
        Universidade  de Lisboa\\
        Av. Rovisco Pais\\
        1049-001 Lisbon\\
        Portugal\\
        \texttt{pelopes@math.ist.utl.pt}\\
}
\date{August 23, 2013}
\maketitle

\begin{abstract}
We prove that if a link admits non-trivial $(2k+1)$-colorings, with prime $2k+1>7$, it also admits non-trivial $(2k+1)$-colorings not involving colors $2k, 2k-1$, nor $k$.
\end{abstract}

\bigbreak

\bigbreak


\bigbreak

\section{Introduction}

\noindent
In \cite{Oshiro} Oshiro proved that colors  $6, 5$ and $3$ can always be removed from a non-trivial $7$-coloring. We mimic Oshiro's approach to prove the following.

\begin{thm}\label{thm:main}
Let $p=2k+1>7$ be a prime. Let $L$ be a non-split link which admits non-trivial $(2k+1)$-colorings. Then colors $2k, 2k-1$, and $k$  are not needed in order to assemble a non-trivial $(2k+1)$-coloring of $L$.
\end{thm}

We now give an outline of the proof.
We consider a diagram endowed with a non-trivial (2k+1)-coloring. For each color, $c$, to be removed the strategy is the following. First we check if the color is used in the coloring. If it is not we are done. Else we consider the following three possibilities,  in turn. The first is to ascertain if the color $c$ is used in a monochromatic crossing. If it is not we go to the next step. Else we remove the color $c$ from each monochromatic crossing in the diagram until there is no monochromatic crossing involving color $c$. The next step is to check if $c$ is used in an over-arc. If it is not we go to the next step. Else we remove $c$ from each over-arc until  there is no over-arc with color $c$. Finally we check if $c$ is used in an under-arc. If it is not we are done. Else we remove it from each under-arc until there is no under-arc using color $c$. This ends the removal of color $c$ from the coloring. The actual removal of the color at each step is done by performing Reidemeister moves on the diagram and consistently re-assigning colors in the neighborhoods affected by the transformations.

The proof is split into three Sections. In Section \ref{sect:2k} we show how to remove $2k$; in Section \ref{sect:remove2k-1} we show how to remove color $2k-1$, provided color $2k$ is no longer being used; and in Section \ref{sect:removek} we show how to remove color $k$ provided colors $2k$ and $2k-1$ are no longer being used.

In the sequel, in any Figure the right-most diagram is the candidate to a diagram without the color under scrutiny, see Figure \ref{fig:eps1} below, for example. For that the colors in the new arcs are analyzed. These colors are forced to equal the forbidden colors so that we identify the new situations we have to tackle. In the first cases we do the explicit calculations but later on we rely on displaying in tables the forbidden identities (in the first row) and their consequences (in the second row). These consequences can be either an inconsistency (which is marked by an $X$) and we can dismiss the situation, or  an expression which helps us identify the new situation to be tackled. These new situations, in turn, are dealt with below in the text. We remark also that the inconsistencies are of the sort ``arc assuming a forbidden color'' or ``$2k+1=$ a prime less than $11$ or a composite integer'', among others.

{\bf Expressions and equalities will be understood modulo $2k+1$}.

\section{Removal of Color $2k$}\label{sect:2k}

\noindent

\begin{figure}[!ht]
    \psfrag{a}{\huge $a$}
    \psfrag{2k}{\huge $2k$}
    \psfrag{2a+1}{\huge $2a+1$}
    \centerline{\scalebox{.350}{\includegraphics{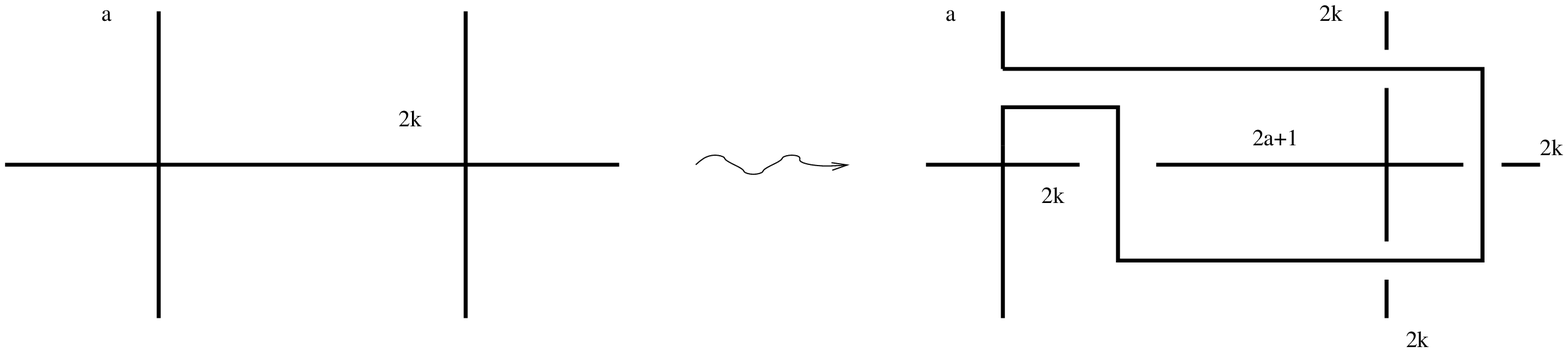}}}
    \caption{Removing color $2k$ from a monochromatic crossing}\label{fig:eps1}
\end{figure}

\subsection{Removal of Color $2k$ from a Monochromatic Crossing (Figure \ref{fig:eps1}).}\label{subsect:2kmono}

\noindent

View Figure \ref{fig:eps1}. We remark that the details of some crossings are irrelevant for the reasonings here. Since we are trying to remove color $2k$ we would not want $2a+1$ in the diagram on the right-hand side to equal $2k$. Let us look into the consequences of such identity. If $2a+1=2k$ then $2a+1=2k+2k+1$ which amounts to $a=2k$. But if there is no pair of crossings in the diagram as depicted in Figure \ref{fig:eps1} such that $a\neq 2k$, then the whole coloring is monochromatic which conflicts with the standing assumption of a non-trivial coloring. (This corresponds to the inconsistencies we discussed above. Below in the tables these instances will be marked with an $X$.) Hence, by repeating this operation we can assume we remove all the monochromatic crossings with color $2k$.

\begin{figure}[!ht]
    \psfrag{a}{\huge $a$}
    \psfrag{2k}{\huge $2k$}
    \psfrag{2a+1}{\huge $2a+1$}
    \psfrag{-2-a}{\huge $-2-a$}
    \psfrag{3a+2}{\huge $3a+2$}
    \centerline{\scalebox{.350}{\includegraphics{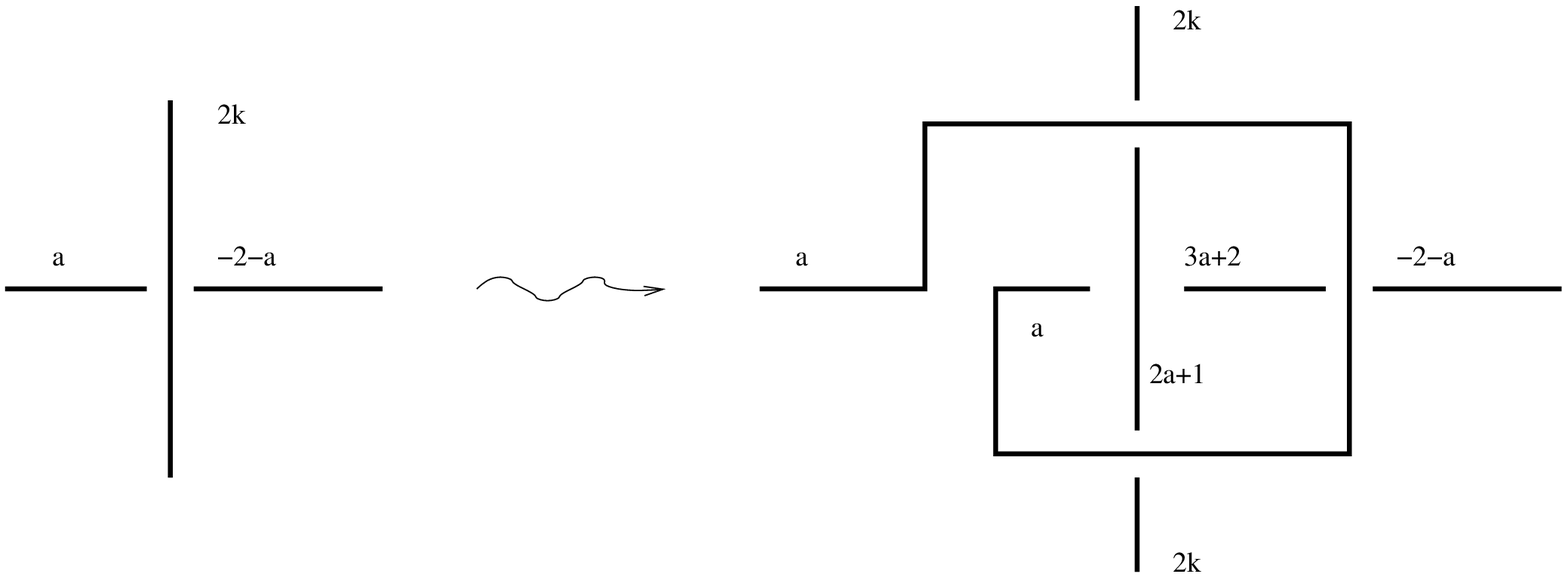}}}
    \caption{Removing color $2k$ from an over-arc}\label{fig:eps2}
\end{figure}

\subsection{Removal of Color $2k$ from an Over-Arc (Figure \ref{fig:eps2}).}\label{subsect:2kover}

\noindent

View Figure \ref{fig:eps2}. We have already seen that $a= 2k$ is equivalent to $2a+1= 2k$ and we now note that $a= 2k$ is also equivalent to $3a+2= 2k$. In this way we have achieved the goal of removing color $2k$ from over-arcs.

\begin{figure}[!ht]
    \psfrag{a}{\huge $a$}
    \psfrag{2k}{\huge $2k$}
    \psfrag{2a+1}{\huge $2a+1$}
    \psfrag{2a-b}{\huge $2a-b$}
    \psfrag{2a-2b-1}{\huge $2a-2b-1$}
    \psfrag{b}{\huge $b$}
    \psfrag{2b+1}{\huge $2b+1$}
    \centerline{\scalebox{.350}{\includegraphics{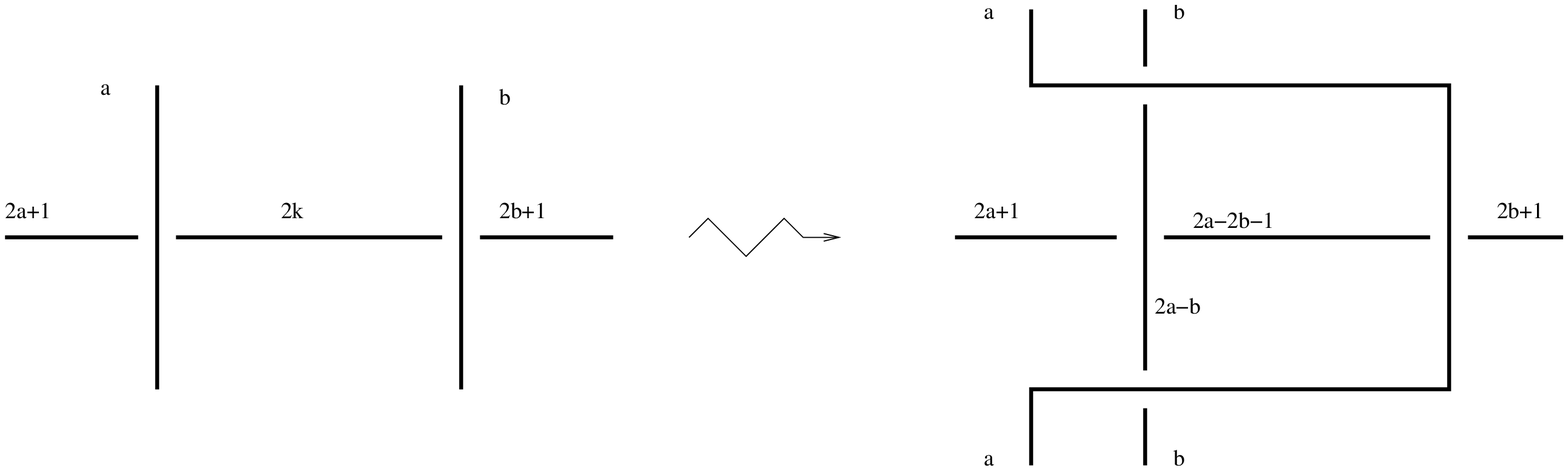}}}
    \caption{Removing color $2k$ from an under-arc: the initial setting}\label{fig:eps3}
\end{figure}

\subsection{Removal of Color $2k$ from an Under-Arc (Figure \ref{fig:eps3}).}\label{subsect:2kunder}

\noindent

View Figure \ref{fig:eps3}.

\begin{table}[h!]
\begin{center}
    \begin{tabular}{| c | c |}\hline
    2a-b=-1 &   2a-2b-1=-1     \\ \hline
            b=2a+1 &   b=a    \\ \hline
    \end{tabular}
\caption{Equalities which should not occur in Figure \ref{fig:eps3} (1st row) and their consequences (2nd row). These situations are dealt with in \ref{subsubsect:2kunderb2a+1} and \ref{subsubsect:2kunderba} below.}\label{Ta:fig:red3bis}
\end{center}
\end{table}

\begin{figure}[!ht]
    \psfrag{a}{\huge $a$}
    \psfrag{2k}{\huge $2k$}
    \psfrag{2a+1}{\huge $2a+1$}
    \psfrag{3a+2}{\huge $3a+2$}
    \psfrag{4a+3}{\huge $4a+3$}
    \centerline{\scalebox{.350}{\includegraphics{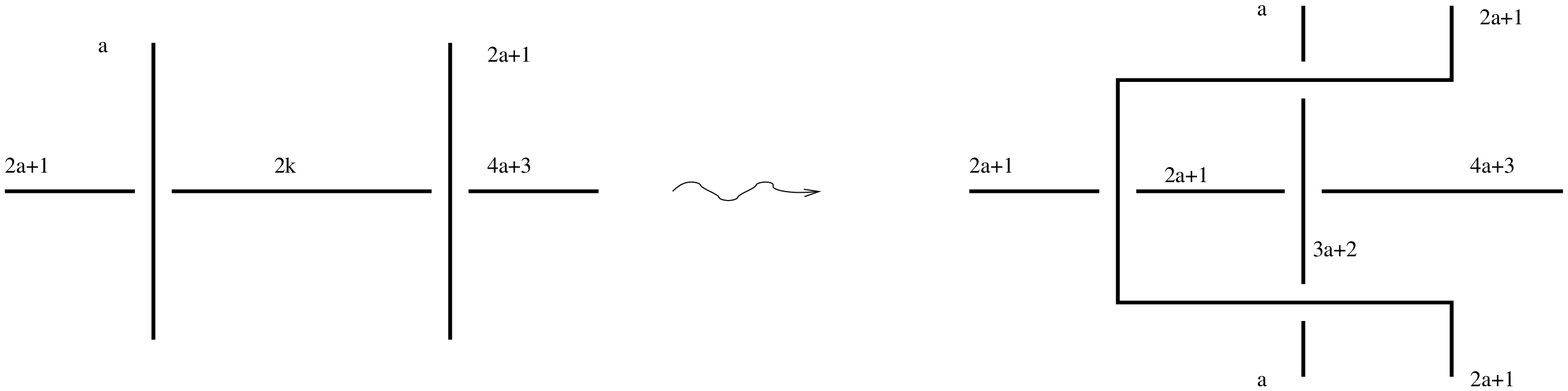}}}
    \caption{Removing color $2k$ from an under-arc: the $b=2a+1$ instance.}\label{fig:red3b=2a+1}
\end{figure}

\subsubsection{\bf The $b=2a+1$ instance (Figure \ref{fig:red3b=2a+1}).}\label{subsubsect:2kunderb2a+1}

 View Figure \ref{fig:red3b=2a+1}.

\begin{table}[h!]
\begin{center}
    \begin{tabular}{| c |}\hline
    3a+2=-1      \\ \hline
    a=-1 \, X   \\ \hline
    \end{tabular}
\caption{Equalities which should not occur in Figure \ref{fig:red3b=2a+1} (1st row) and their consequences (2nd row). $X$'s stand for  conflicts with   assumptions thus not requiring further considerations.}\label{Ta:fig:red3b=2a+1}
\end{center}
\end{table}

\begin{figure}[!ht]
    \psfrag{a}{\huge $a$}
    \psfrag{2k}{\huge $2k$}
    \psfrag{2a+1}{\huge $2a+1$}
    \psfrag{3a+2}{\huge $3a+2$}
    \psfrag{4a+3}{\huge $4a+3$}
    \centerline{\scalebox{.350}{\includegraphics{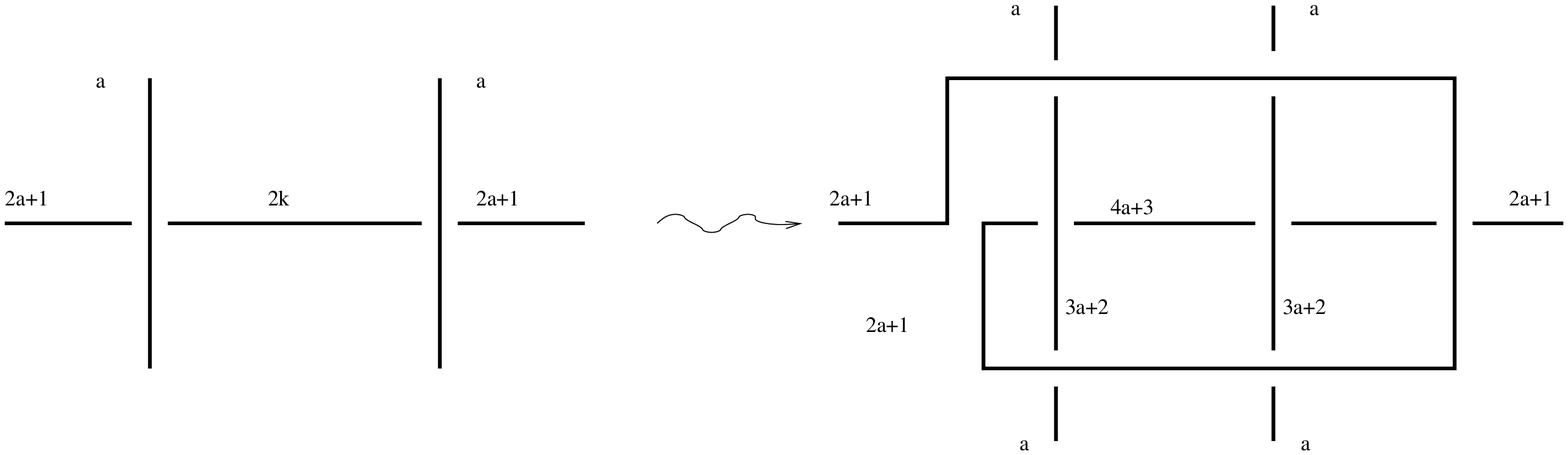}}}
    \caption{Removing color $2k$ from an under-arc: the $b=a$ instance}\label{fig:red3b=a}
\end{figure}

\subsubsection{\bf The $b=a$ instance (Figure \ref{fig:red3b=a}).}\label{subsubsect:2kunderba}

 View Figure \ref{fig:red3b=a}.

\begin{table}[h!]
\begin{center}
    \begin{tabular}{| c | c |}\hline
    3a+2=-1 & 4a+3=-1     \\ \hline
            a=-1 \, X  &   a=-1  \, X  \\ \hline
    \end{tabular}
\caption{Equalities which should not occur in Figure \ref{fig:red3b=a} (1st row) and their consequences (2nd row). $X$'s stand for  conflicts with   assumptions thus not requiring further considerations.}\label{Ta:fig:red3b=a}
\end{center}
\end{table}

This concludes the removal of color $2k$ from under-arcs and consequently concludes the removal of color $2k$ from the coloring.

\section{Removal of Color $2k-1$, provided  color $2k$ has already been removed}\label{sect:remove2k-1}

\noindent

\begin{figure}[!ht]
    \psfrag{a}{\huge $a$}
    \psfrag{2k-1}{\huge $2k-1$}
    \psfrag{2a+2}{\huge $2a+2$}
    \centerline{\scalebox{.350}{\includegraphics{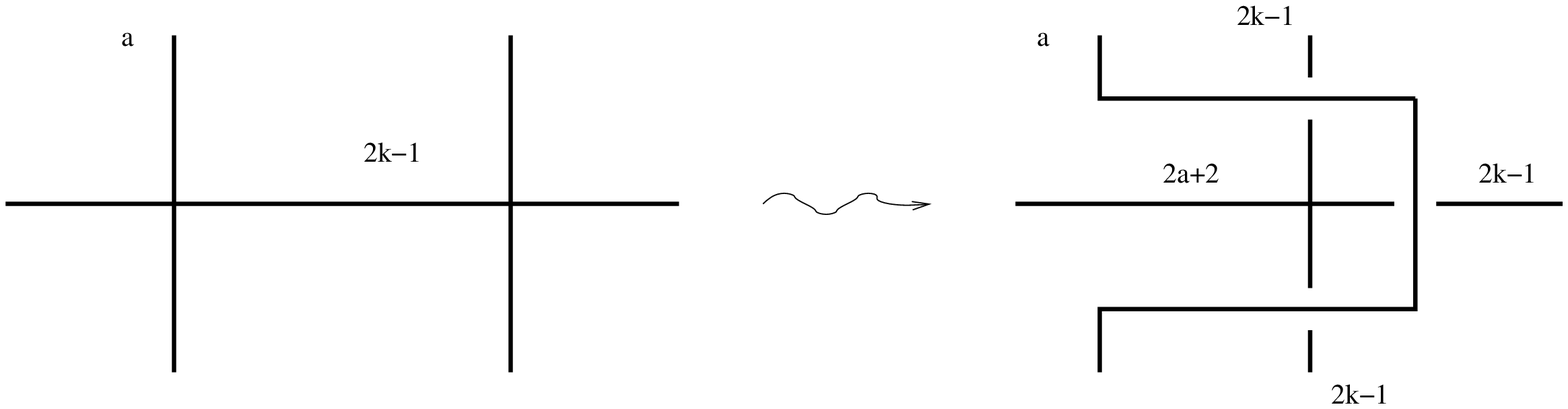}}}
    \caption{Removing color $2k-1$ from a monochromatic crossing}\label{fig:red4}
\end{figure}

\begin{figure}[!ht]
    \psfrag{a}{\huge $a$}
    \psfrag{2k-1}{\huge $2k-1$}
    \psfrag{k-1}{\huge $k-1$}
    \psfrag{k-2}{\huge $k-2$}
    \psfrag{-3}{\huge $-3$}
    \psfrag{...}{\huge $\vdots$}
    \centerline{\scalebox{.350}{\includegraphics{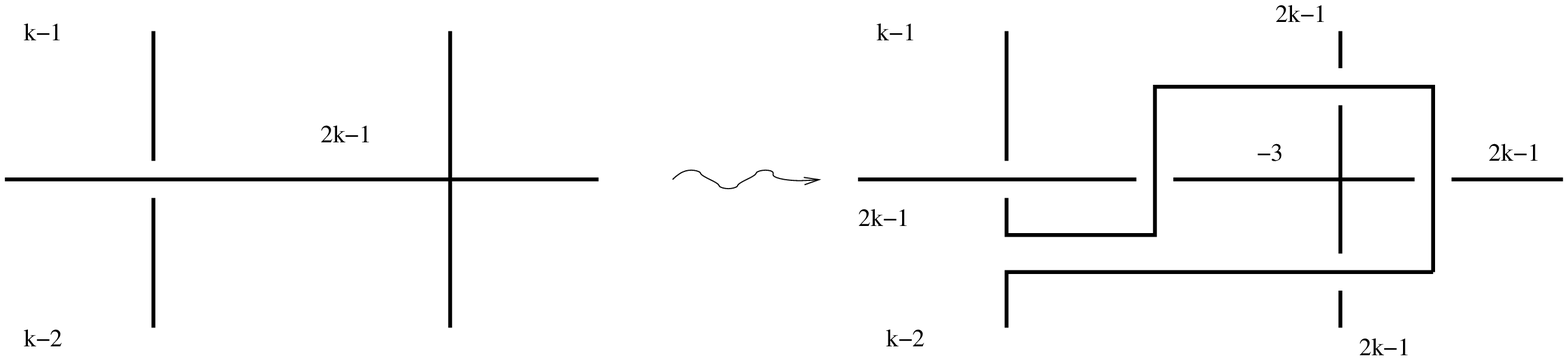}}}
    \caption{Removing color $2k-1$ from a monochromatic crossing (remaining instance)}\label{fig:red4bis}
\end{figure}

\subsection{Removal of Color $2k-1$ from a Monochromatic Crossing (Figures \ref{fig:red4} and \ref{fig:red4bis}).}\label{subsect:2k-1mono}

\noindent

View Figure \ref{fig:red4}. If $2a+2 =-2$ then $a=-2$ which conflicts with the assumptions. If $2a+2 =-1 +2k+1$ then $a=k-1$. In Figure \ref{fig:red4} if the arc receiving color $a=k-1$ is an over-arc then the under-arc to its left receives color $2(k-1)-(2k-1)=-1$ which conflicts with the standing assumptions. The remaining instance is contemplated in Figure \ref{fig:red4bis}.

\begin{figure}[!ht]
    \psfrag{a}{\huge $a$}
    \psfrag{2k-1}{\huge $2k-1$}
    \psfrag{2a+2}{\huge $2a+2$}
    \psfrag{-4-a}{\huge $-4-a$}
    \psfrag{3a+4}{\huge $3a+4$}
    \centerline{\scalebox{.350}{\includegraphics{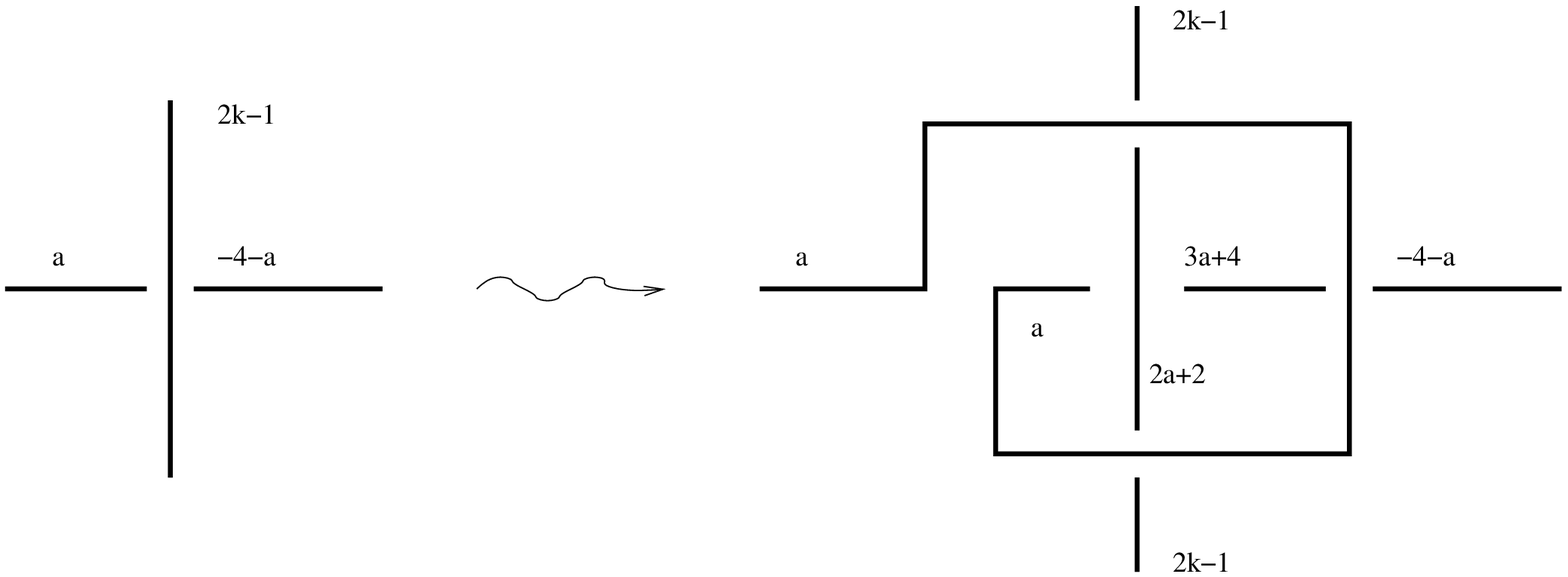}}}
    \caption{Removing color $2k-1$ from an over-arc}\label{fig:red5}
\end{figure}

\begin{figure}[!ht]
    \psfrag{a}{\huge $a$}
    \psfrag{2k-1}{\huge $2k-1$}
    \psfrag{-4-a}{\huge $-4-a$}
    \psfrag{-6-2a}{\huge $-6-2a$}
    \psfrag{-8-3a}{\huge $-8-3a$}
    \centerline{\scalebox{.350}{\includegraphics{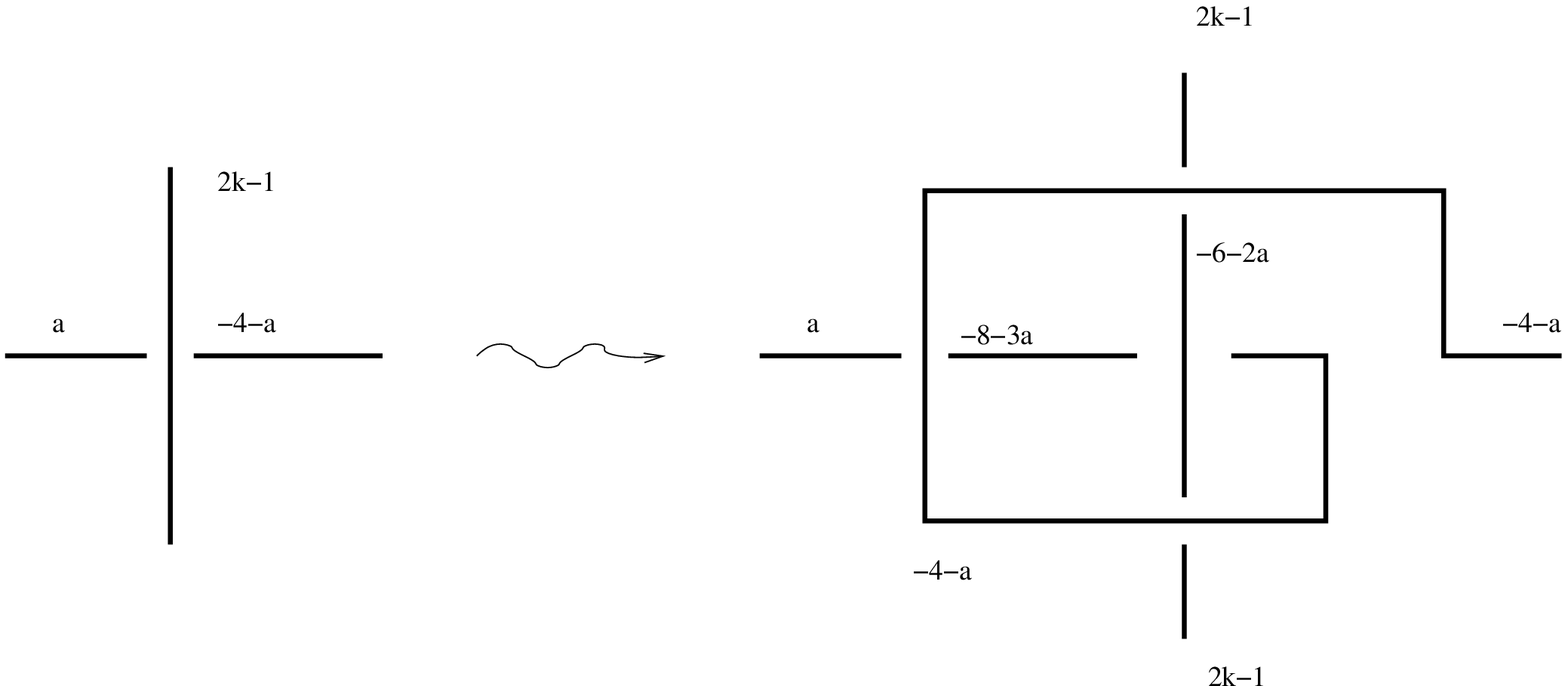}}}
    \caption{Removing color $2k-1$ from an over-arc}\label{fig:red5bis}
\end{figure}

\subsection{Removal of Color $2k-1$ from an Over-Arc (Figures \ref{fig:red5} and \ref{fig:red5bis}).}\label{subsect:2k-1over}

\noindent

View Figure \ref{fig:red5}.

\begin{table}[h!]
\begin{center}
    \begin{tabular}{| c | c | c | c |}\hline
2a+2=-1 & 2a+2=-2 & 3a+4=-1 & 3a+4=-2    \\ \hline
a=k-1 & a=-2\, X  &   3a=-5  &  a=-2 \, X  \\ \hline
    \end{tabular}
\caption{Equalities which should not occur in Figure \ref{fig:red5} (1st row) and their consequences (2nd row). $X$'s stand for  conflicts with assumptions thus not requiring further considerations. The other situations are dealt with below in \ref{subsubsect:over-arc2k-1ak-1} and \ref{subsubsect:over-arc2k-13a-5}.}\label{Ta:fig:red5}
\end{center}
\end{table}

\subsubsection{The $a=k-1$ instance (Figure \ref{fig:red5bis}).}\label{subsubsect:over-arc2k-1ak-1} View Figure \ref{fig:red5bis}. In this case
\[
-6-2a=-6-2k+2=-3
\]
\[
-8-3a=-8-3k+3=k-3
\]

\begin{table}[h!]
\begin{center}
    \begin{tabular}{| c | c |}\hline
    k-3=-1 & k-3=-2     \\ \hline
            2k+1=5 \, X  &   2k+1=3  \, X  \\ \hline
    \end{tabular}
\caption{Equalities which should not occur in Figure \ref{fig:red5bis} in the $a=k-1$ instance (1st row) and their consequences (2nd row). $X$'s stand for  conflicts with assumptions thus not requiring further considerations.}\label{Ta:fig:red5bis2}
\end{center}
\end{table}

\subsubsection{The $3a=-5$ instance (Figure \ref{fig:red5bis})}\label{subsubsect:over-arc2k-13a-5} View Figure \ref{fig:red5bis} again. In this case

\[
-8-3a=-8+5=-3
\]
\[
-6-2a=-6-3a+a=-6+5+a=a-1
\]

\begin{table}[h!]
\begin{center}
    \begin{tabular}{| c | c |}\hline
    a-1=-1 $\Rightarrow$ a=0 $\Rightarrow$ -5=3a=0  & a-1=-2     \\ \hline
            2k+1=5 \, X  &   a=-1  \, X  \\ \hline
    \end{tabular}
\caption{Equalities which should not occur in Figure \ref{fig:red5bis} in the $-3a=5$ instance (1st row) and their consequences (2nd row). $X$'s stand for  conflicts with assumptions thus not requiring further considerations.}\label{Ta:fig:red5bis}
\end{center}
\end{table}

\begin{figure}[!ht]
    \psfrag{a}{\huge $a$}
    \psfrag{2k-1}{\huge $2k-1$}
    \psfrag{2a+2}{\huge $2a+2$}
    \psfrag{2b+2}{\huge $2b+2$}
    \psfrag{2a-2b-2}{\huge $2a-2b-2$}
    \psfrag{b}{\huge $b$}
    \psfrag{2a-b}{\huge $2a-b$}
    \centerline{\scalebox{.350}{\includegraphics{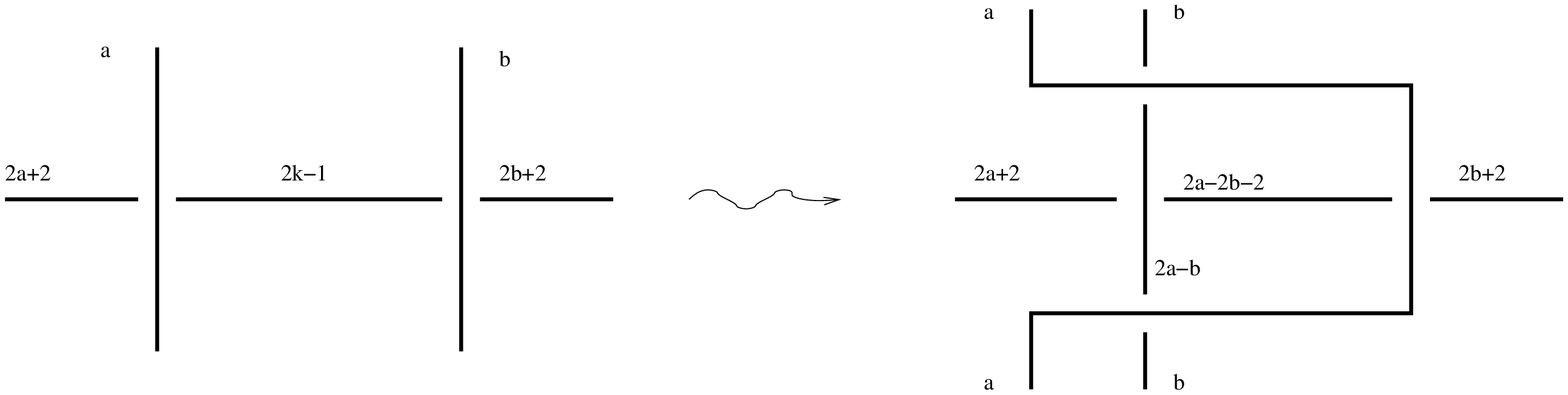}}}
    \caption{Removing color $2k-1$ from an under-arc}\label{fig:red6}
\end{figure}

\subsection{Removal of Color $2k-1$ from an Under-Arc (Figure \ref{fig:red6}).}\label{subsect:2k-1under}

\noindent

View Figure \ref{fig:red6}.

\begin{table}[h!]
\begin{center}
    \begin{tabular}{| c | c | c | c |}\hline
2a-b=2k &  2a-b=2k-1 &  2a-2b-2=2k &  2a-2b-2=2k-1   \\ \hline
b=2a+1 & b=2a+2  &  b=a+k  &   b=a  \\ \hline
    \end{tabular}
\caption{Equalities which should not occur in Figure \ref{fig:red6} (1st row) and their consequences (2nd row). $X$'s stand for  conflicts with assumptions thus not requiring further considerations. These situations are dealt with below in \ref{subsubsect:2k-1underb2a+1}, \ref{subsubsect:2k-1underb2a+2}, \ref{subsubsect:2k-1underba+k}, and \ref{subsubsect:2k-1underba}.}\label{Ta:fig:red6}
\end{center}
\end{table}

\begin{figure}[!ht]
    \psfrag{a}{\huge $a$}
    \psfrag{2k-1}{\huge $2k-1$}
    \psfrag{2a+1}{\huge $2a+1$}
    \psfrag{2a}{\huge $2a$}
    \psfrag{2a+2}{\huge $2a+2$}
    \psfrag{3a+2}{\huge $3a+2$}
    \psfrag{4a+4}{\huge $4a+4$}
    \centerline{\scalebox{.350}{\includegraphics{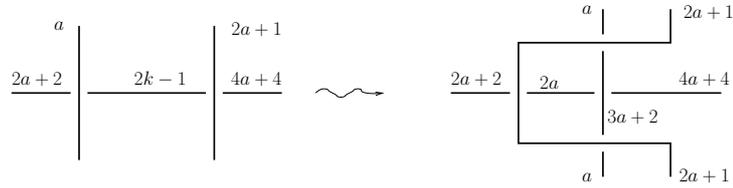}}}
    \caption{Removing color $2k-1$ from an under-arc: the $b=2a+1$ instance}\label{fig:red6b2a+1}
\end{figure}

\subsubsection{\bf The $b=2a+1$ instance (Figure \ref{fig:red6b2a+1}).}\label{subsubsect:2k-1underb2a+1}

View Figure \ref{fig:red6b2a+1}. A few words about the sub-instance $3a+2=2k-1$, see Table \ref{Ta:fig:red6b2a+1}. In order to solve this equation for $a$ it is useful to know how the modulus $2k+1$ can be expressed mod $3$. If $2k+1=3l$ for some integer $l$ then $2k+1$ would be composite which conflicts with the assumptions. So the possibilities are $2k+1=3l+1$ or $2k+1=2l+2$. This is how the two consequences in the bottom right slot of Table \ref{Ta:fig:red6b2a+1} come about. Similar situations occur below with a parameter being multiplied by a constant. They are dealt with analogously.

\bigbreak

\bigbreak

\begin{table}[h!]
\begin{center}
    \begin{tabular}{| c | c | c | c |}\hline
2a=2k & 2a=2k-1 & 3a+2=2k & 3a+2=2k-1   \\ \hline
a=k & a=-1\, X  &  a=-1 \, X &   a$=_{3l+1}$l-1 or a$=_{3l+2}$2l  \\ \hline
    \end{tabular}
\caption{Equalities which should not occur in Figure \ref{fig:red6b2a+1} (1st row) and their consequences (2nd row). $X$'s stand for  conflicts with assumptions thus not requiring further considerations. The other situations are dealt with below in \ref{para:2k-1underb2a+1ak}, \ref{para:2k-1underb2a+1al-1}, and \ref{para:2k-1underb2a+1a2l}}\label{Ta:fig:red6b2a+1}
\end{center}
\end{table}
\bigbreak

\bigbreak

\bigbreak

\begin{figure}[!ht]
    \psfrag{0}{\huge $0$}
    \psfrag{1}{\huge $1$}
    \psfrag{-1}{\huge $-1$}
    \psfrag{2}{\huge $2$}
    \psfrag{k}{\huge $k$}
    \psfrag{k+1}{\huge $k+1$}
    \psfrag{-2}{\huge $-2$}
    \centerline{\scalebox{.350}{\includegraphics{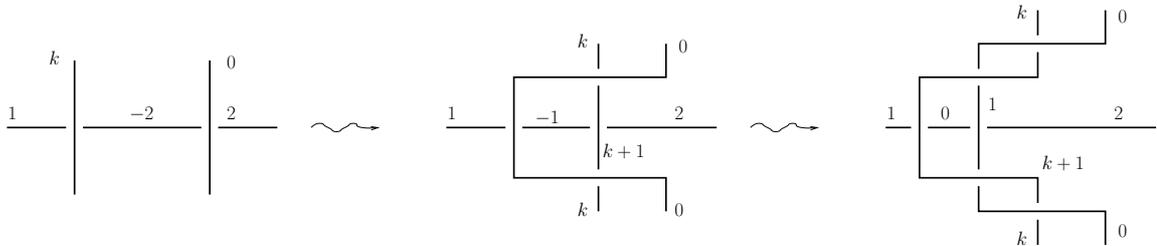}}}
    \caption{Removing color $2k-1$ from an under-arc: the $b=2a+1$ with $a=k$ sub-instance }\label{fig:red7}
\end{figure}

\paragraph{\bf The $a=k$ instance (Figure \ref{fig:red7})}\label{para:2k-1underb2a+1ak}

Working over modulus $2k+1$ with $a=k$, then $2a+2=2k+2=1$, $2a+1=2k+1=0$ and $4a+4=4k+4=2(2k+1)+2=2$. Now Figure \ref{fig:red7} should be self-explanatory.

\bigbreak

\bigbreak

\begin{figure}[!ht]
    \psfrag{0}{\huge $0$}
    \psfrag{-2}{\huge $-2$}
    \psfrag{2l}{\huge $2l$}
    \psfrag{2l-1}{\huge $2l-1$}
    \psfrag{l+1}{\huge $l+1$}
    \psfrag{2l-2}{\huge $2l-2$}
    \psfrag{l-1}{\huge $l-1$}
    \psfrag{l+1}{\huge $l+1$}
    \psfrag{3l-3}{\huge $3l-3$}
    \psfrag{3l+1}{\huge $\mathbf{3l+1}$}
    \centerline{\scalebox{.350}{\includegraphics{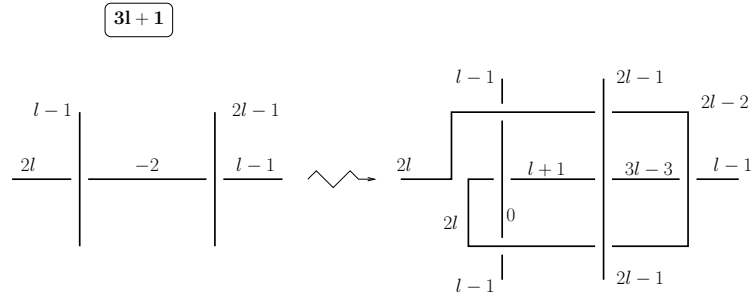}}}
    \caption{Removing color $2k-1$ from an under-arc: the $b=2a+1$ with $a=l-1$ sub-instance }\label{fig:red8}
\end{figure}

\paragraph{\bf The $a=l-1$ mod $3l+1$ instance (Figure \ref{fig:red8})}\label{para:2k-1underb2a+1al-1}
Working over modulus $3l+1$ with $a=l-1$, then $2a+1=2(l-1)+1=2l-1$, $2a+2=2l$ and $4a+4=4l=3l+1+l-1=l-1$. Remarking that $l\geq 4$, figure \ref{fig:red8} should now be self-explanatory.

\begin{table}[h!]
\begin{center}
    \begin{tabular}{| c | c | c | c |}\hline
l+1=3l & l+1=3l-1& 3l-3=3l & 3l-3=3l-1   \\ \hline
l+1($\leq$ 2l)$<$3l \, X & l=1\, X  &  3=0 \, X &   2=0 \, X  \\ \hline
    \end{tabular}
\caption{Equalities which should not occur in Figure \ref{fig:red8} (1st row) and their consequences (2nd row). $X$'s stand for  conflicts with assumptions thus not requiring further considerations.}\label{Ta:fig:red8}
\end{center}
\end{table}

\bigbreak

\begin{figure}[!ht]
    \psfrag{0}{\huge $0$}
    \psfrag{-2}{\huge $-2$}
    \psfrag{2l}{\huge $2l$}
    \psfrag{2l+2}{\huge $2l+2$}
    \psfrag{l-1}{\huge $l-1$}
    \psfrag{l-2}{\huge $l-2$}
    \psfrag{-4}{\huge $-4$}
    \psfrag{l}{\huge $l$}
    \psfrag{3l+2}{\huge $\mathbf{3l+2}$}
    \centerline{\scalebox{.350}{\includegraphics{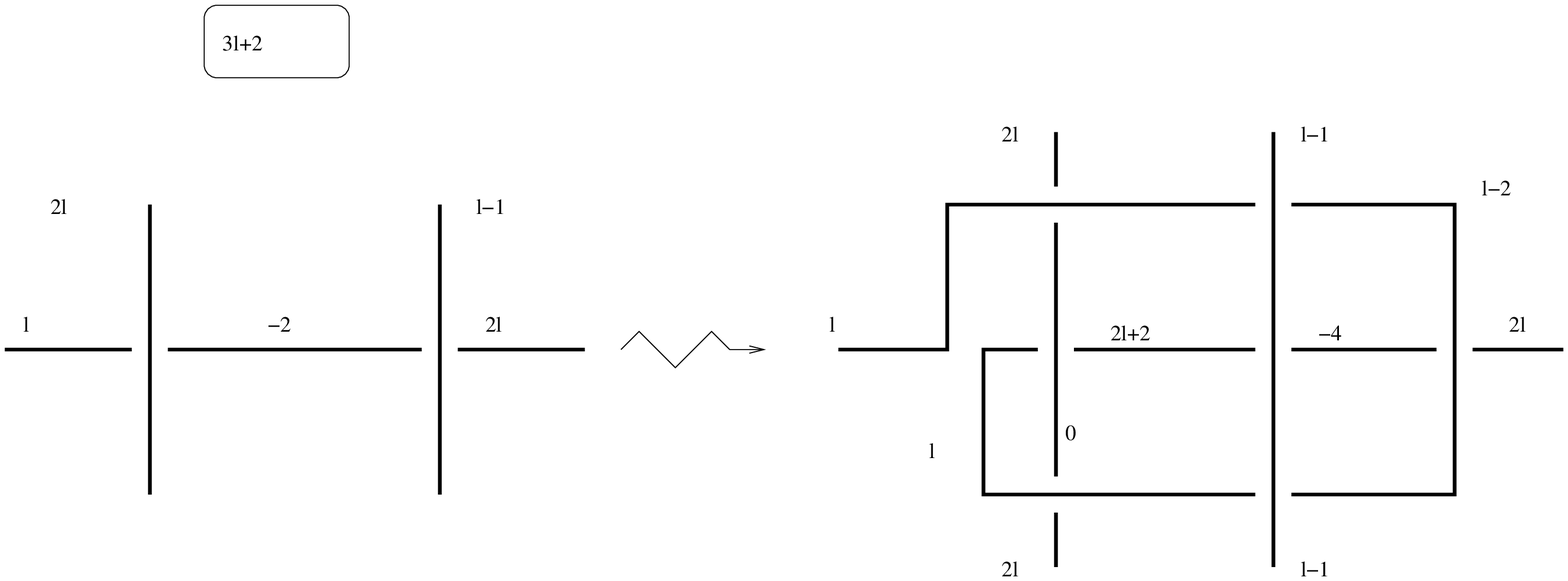}}}
    \caption{Removing color $2k-1$ from an under-arc: the $b=2a+1$ with $a=2l$ sub-instance }\label{fig:red9}
\end{figure}

\paragraph{\bf The $a=2l$ mod $3l+2$ instance (Figure \ref{fig:red9})}\label{para:2k-1underb2a+1a2l}

Working over modulus $3l+2$ with $a=2l$, then $2a+1=2(2l)+1=4l+1=3l+2+l-1=l-1$, $2a+2=l$,  and $4a+4=2l$. $l\geq 3$ and Figure \ref{fig:red9} should now be self-explanatory.

\bigbreak

\bigbreak

\begin{table}[h!]
\begin{center}
    \begin{tabular}{| c | c |}\hline
2l+2=3l+1 & 2l+2=3l    \\ \hline
l=1 \, X & l=2\, X    \\ \hline
    \end{tabular}
\caption{Equalities which should not occur in Figure \ref{fig:red9} (1st row) and their consequences (2nd row). $X$'s stand for  conflicts with assumptions thus not requiring further considerations.}\label{Ta:fig:red9}
\end{center}
\end{table}

\bigbreak

\begin{figure}[!ht]
    \psfrag{a}{\huge $a$}
    \psfrag{2a+2}{\huge $2a+2$}
    \psfrag{3a+4}{\huge $3a+4$}
    \psfrag{4a+6}{\huge $4a+6$}
    \psfrag{-2}{\huge $-2$}
    \centerline{\scalebox{.350}{\includegraphics{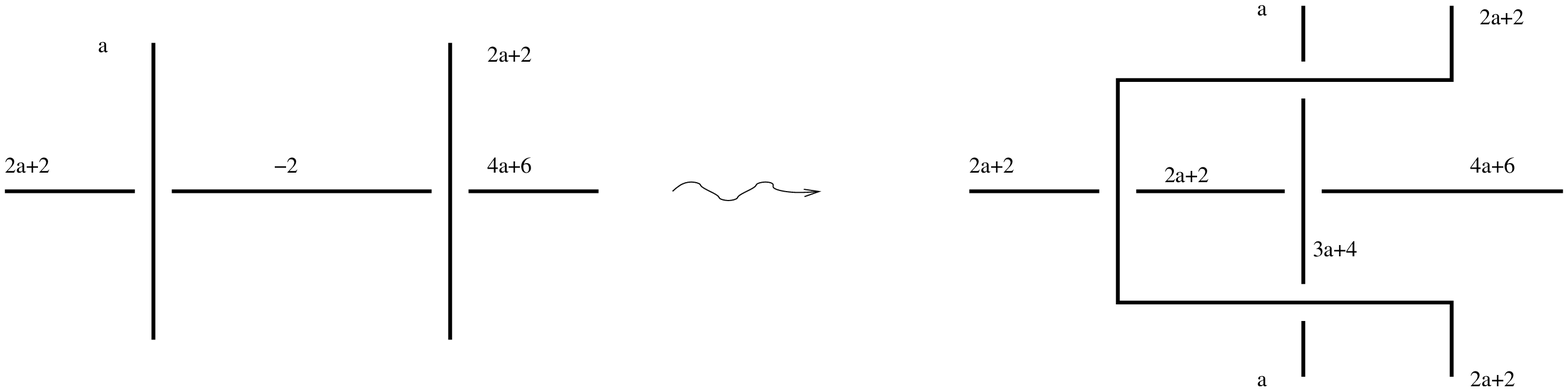}}}
    \caption{Removing color $2k-1$ from an under-arc: the $b=2a+2$ instance }\label{fig:red10}
\end{figure}

\subsubsection{\bf The $b=2a+2$ instance (Figure \ref{fig:red10}).}\label{subsubsect:2k-1underb2a+2}

\bigbreak

 View Figure \ref{fig:red10}.

\begin{table}[h!]
\begin{center}
    \begin{tabular}{| c | c | }\hline
3a+4 = 2k & 3a+4=2k-1   \\ \hline
a$=_{3l+1}$2l-1 or a$=_{3l+2}$l-1 &  a=-2 \, X  \\ \hline
    \end{tabular}
\caption{Equalities which should not occur in Figure \ref{fig:red10} (1st row) and their consequences (2nd row). $X$'s stand for  conflicts with assumptions thus not requiring further considerations. The other situations are dealt with below in \ref{para:2k-1underb2a+2a2l-1}, and \ref{para:2k-1underb2a+2al-1}}\label{Ta:fig:red10}
\end{center}
\end{table}
\bigbreak

\begin{figure}[!ht]
    \psfrag{0}{\huge $0$}
    \psfrag{-2}{\huge $-2$}
    \psfrag{2l}{\huge $2l$}
    \psfrag{2l-1}{\huge $2l-1$}
    \psfrag{2l+1}{\huge $2l+1$}
    \psfrag{3l+1}{\huge $\mathbf{3l+1}$}
    \psfrag{l+1}{\huge $l+1$}
    \psfrag{2l-2}{\huge $2l-2$}
    \psfrag{2l+2}{\huge $2l+2$}
    \psfrag{l-1}{\huge $l-1$}
    \psfrag{l+1}{\huge $l+1$}
    \psfrag{3l-3}{\huge $3l-3$}
    \psfrag{l-2}{\huge $l-2$}
    \psfrag{-4}{\huge $-4$}
    \psfrag{l}{\huge $l$}
    \centerline{\scalebox{.350}{\includegraphics{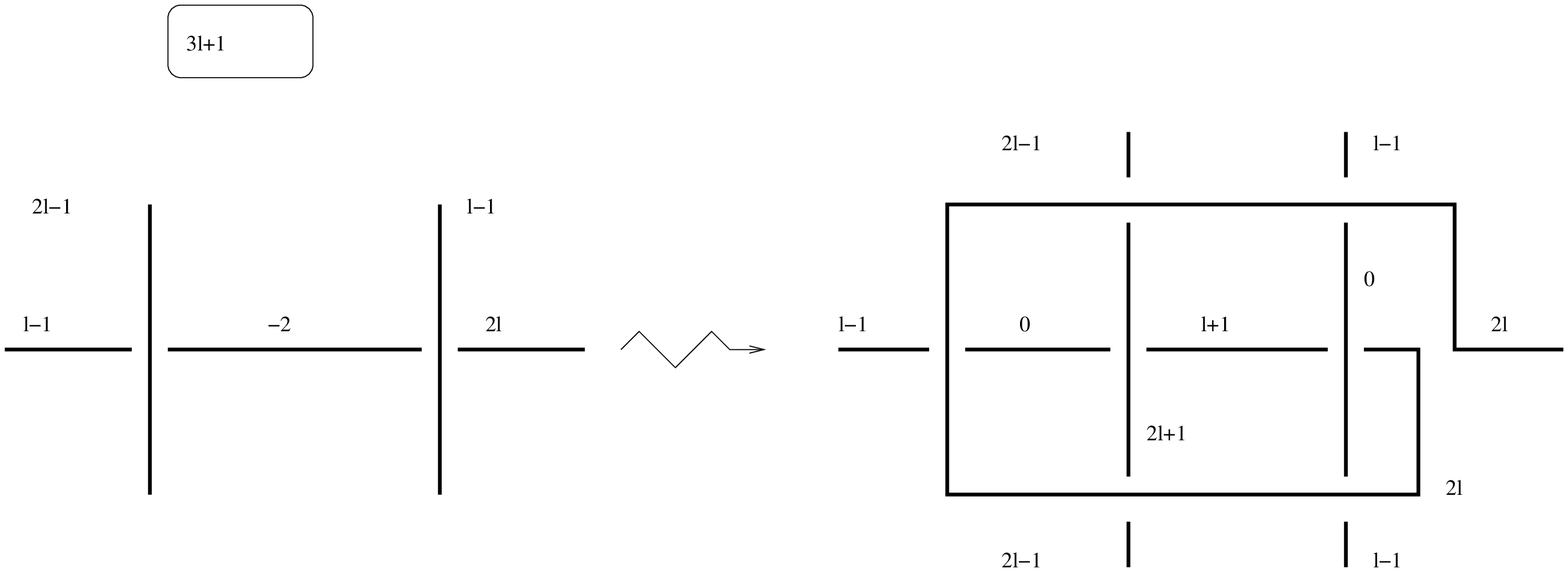}}}
    \caption{Removing color $2k-1$ from an under-arc: the $b=2a+2$ with $a=2l-1$ sub-instance}\label{fig:red11}
\end{figure}

\paragraph{\bf The $b=2a+2$ with $a=2l-1$ mod $3l+1$ instance (Figure \ref{fig:red11})}\label{para:2k-1underb2a+2a2l-1}

Working over modulus $3l+1$ with $a=2l-1$, then $2a+2=2(2l-1)+2=4l=3l+1+l-1=l-1$ and $4a+6=2l$

\begin{table}[h!]
\begin{center}
    \begin{tabular}{| c | c | c | c |}\hline
l+1=3l & l+1=3l-1& 2l+1=3l & 2l+1=3l-1   \\ \hline
l+1($\leq$ 2l)$<$3l \, X & l=1\, X  &  l=1 \, X &   l=2 \, X  \\ \hline
    \end{tabular}
\caption{Equalities which should not occur in Figure \ref{fig:red11} (1st row) and their consequences (2nd row). $X$'s stand for  conflicts with assumptions thus not requiring further considerations.}\label{Ta:fig:red11}
\end{center}
\end{table}

\begin{figure}[!ht]
    \psfrag{0}{\huge $0$}
    \psfrag{-2}{\huge $-2$}
    \psfrag{2l}{\huge $2l$}
    \psfrag{2l-1}{\huge $2l-1$}
    \psfrag{2l+1}{\huge $2l+1$}
    \psfrag{l+1}{\huge $l+1$}
    \psfrag{2l-2}{\huge $2l-2$}
    \psfrag{2l+2}{\huge $2l+2$}
    \psfrag{3l+2}{\huge $\mathbf{3l+2}$}
    \psfrag{l-1}{\huge $l-1$}
    \psfrag{l+1}{\huge $l+1$}
    \psfrag{3l-3}{\huge $3l-3$}
    \psfrag{l-2}{\huge $l-2$}
    \psfrag{-4}{\huge $-4$}
    \psfrag{l}{\huge $l$}
    \centerline{\scalebox{.350}{\includegraphics{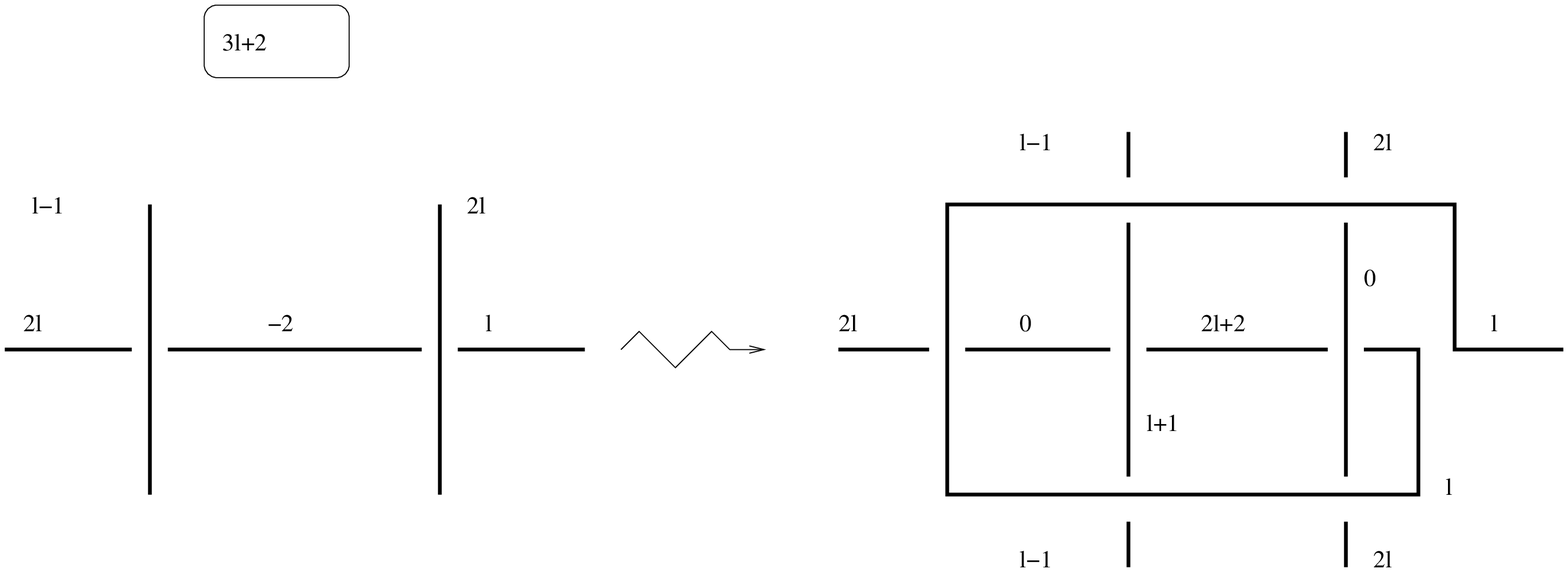}}}
    \caption{Removing color $2k-1$ from an under-arc: the $b=2a+2$ with $a=l-1$ sub-instance}\label{fig:red12}
\end{figure}

\paragraph{\bf The $b=2a+2$  with $a=l-1$ mod $3l+2$ instance (Figure \ref{fig:red12})}\label{para:2k-1underb2a+2al-1}

\bigbreak

Working over modulus $3l+2$ with $a=l-1$, then $2a+2=2(l-1)+2=2l$ and $4a+6=4(l-1)+6=4l+2=3l+2+l=l$

\bigbreak
\begin{table}[h!]
\begin{center}
    \begin{tabular}{| c | c | c | c |}\hline
l+1=3l+1 & l+1=3l& 2l+2=3l+1 & 2l+2=3l   \\ \hline
l=0 or 2=0 \, X  &l+1($\leq$ 2l)$<$3l \, X &   l=1 \, X &   l=2 \, X  \\ \hline
    \end{tabular}
\caption{Equalities which should not occur in Figure \ref{fig:red12} (1st row) and their consequences (2nd row). $X$'s stand for  conflicts with assumptions thus not requiring further considerations.}\label{Ta:fig:red12}
\end{center}
\end{table}

\begin{figure}[!ht]
    \psfrag{a}{\huge $a$}
    \psfrag{2a+1}{\huge $2a+1$}
    \psfrag{2a+2}{\huge $2a+2$}
    \psfrag{-2}{\huge $-2$}
    \psfrag{a+k}{\huge $a+k$}
    \psfrag{a-1}{\huge $a-1$}
    \psfrag{-3}{\huge $-3$}
    \centerline{\scalebox{.350}{\includegraphics{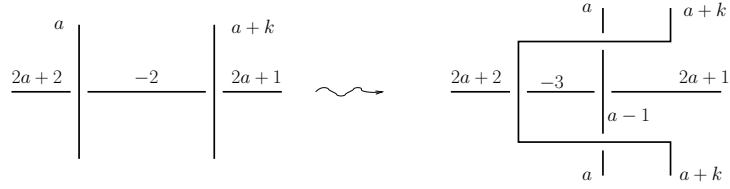}}}
    \caption{Removing color $2k-1$ from an under-arc: the $b=a+k$ instance. The $a=0$ situation has been contemplated in
Figure \ref{fig:red7}; $a-1=-2$ conflicts with the assumption $a\neq-1$.}\label{fig:red13}
\end{figure}

\subsubsection{\bf The $b=a+k$ instance (Figure \ref{fig:red13}).}\label{subsubsect:2k-1underba+k}

\bigbreak

 View Figure \ref{fig:red13}.

\bigbreak

\begin{figure}[!ht]
    \psfrag{a}{\huge $a$}
    \psfrag{2a+2}{\huge $2a+2$}
    \psfrag{3a+4}{\huge $3a+4$}
    \psfrag{4a+6}{\huge $4a+6$}
    \psfrag{-2}{\huge $-2$}
    \centerline{\scalebox{.350}{\includegraphics{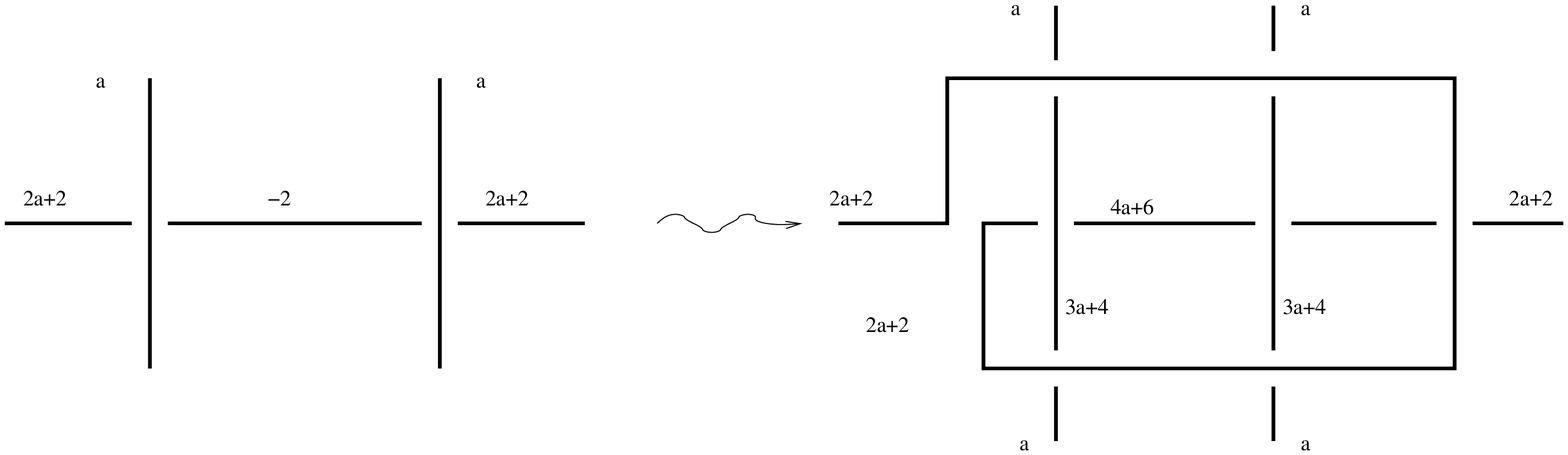}}}
    \caption{Removing color $2k-1$ from an under-arc: the $b=a$ instance}\label{fig:red14}
\end{figure}

\subsubsection{\bf The $b=a$ instance (Figure \ref{fig:red14}).}\label{subsubsect:2k-1underba}

\bigbreak

 View Figure \ref{fig:red14}.

\begin{table}[h!]
\begin{center}
    \begin{tabular}{| c | c |  c | c | }\hline
4a+6=2k &  4a+6=2k-1 & 3a+4=2k & 3a+4=2k-1 \\ \hline
a$=_{4l+1}$3l-1 or a$=_{4l+3}$l-1 & a=-2 \, X  & a$=_{3l+1}$2l-1 or a$=_{3l+2}$l-1 & a=-2 \, X \\ \hline
    \end{tabular}
\caption{Equalities which should not occur in Figure \ref{fig:red14} (1st row) and their consequences (2nd row). $X$'s stand for  conflicts with assumptions thus not requiring further considerations. The other situations are dealt with below in \ref{para:2k-1underbaa3l-1} and in \ref{para:2k-1underbaal-14l+3}.}\label{Ta:fig:red14}
\end{center}
\end{table}

\bigbreak

\begin{figure}[!ht]
    \psfrag{4l+1}{\huge $\mathbf{4l+1}$}
    \psfrag{-2}{\huge $-2$}
    \psfrag{-3}{\huge $-3$}
    \psfrag{2l}{\huge $2l$}
    \psfrag{2l-1}{\huge $2l-1$}
    \psfrag{2l+1}{\huge $2l+1$}
    \psfrag{l+1}{\huge $l+1$}
    \psfrag{3l-1}{\huge $3l-1$}
    \psfrag{2l-3}{\huge $2l-3$}
    \psfrag{l-3}{\huge $l-3$}
    \psfrag{2l-2}{\huge $2l-2$}
    \psfrag{2l+2}{\huge $2l+2$}
    \psfrag{l-1}{\huge $l-1$}
    \psfrag{l+1}{\huge $l+1$}
    \psfrag{3l-3}{\huge $3l-3$}
    \psfrag{l-2}{\huge $l-2$}
    \centerline{\scalebox{.35}{\includegraphics{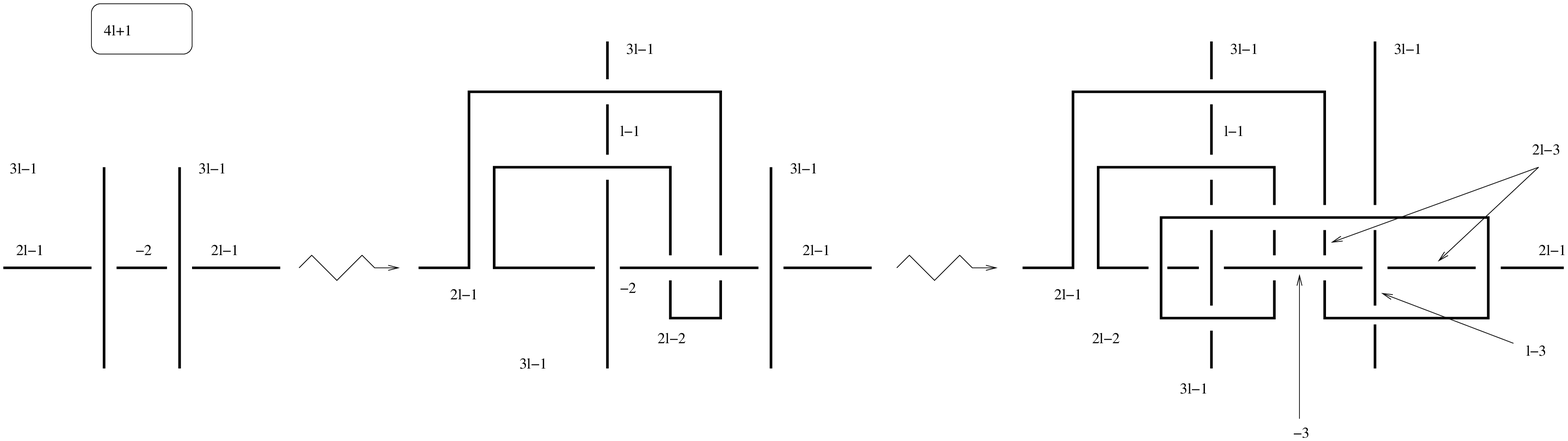}}}
    \caption{Removing color $2k-1$ from an under-arc: the $b=a$ with $a=3l-1$ sub-instance}\label{fig:red15}
\end{figure}

\paragraph{\bf The $b=a$ with $a=3l-1$ mod $4l+1$ instance (Figure \ref{fig:red15}).}\label{para:2k-1underbaa3l-1}

\bigbreak
With $a=3l-1$ mod $4l+1$, $2a+2=6l=4l+1+2l-1=2l-1$. Noting that $l\geq 4$, Figure \ref{fig:red15} is self-explanatory.

\begin{table}[h!]
\begin{center}
  \scalebox{.91}{  \begin{tabular}{| c | c |  c | c | c | c |  c | c |}\hline
l-3=-1 & l-3=-2  & l-1=-1 & l-1=4l-1 & 2l-3=-1 & 2l-3=4l-1  & 2l-2=4l & 2l-2=-2 \\ \hline
l=2 \, X & l=1\, X  & l=0\, X & l=0 or 3=0 \, X & l=1  \, X & 2l-3($<$3l-2)$<$4l-1 \, X & 2l-1($<$3l)$<$4l\, X  & l=0 or 2=0\, X \\ \hline
    \end{tabular}}
\caption{Equalities which should not occur in Figure \ref{fig:red15} (1st row) and their consequences (2nd row). $X$'s stand for  conflicts with assumptions thus not requiring further considerations.}\label{Ta:fig:red15}
\end{center}
\end{table}

\bigbreak

\begin{figure}[!ht]
    \psfrag{-2}{\huge $-2$}
    \psfrag{-3}{\huge $-3$}
    \psfrag{2l}{\huge $2l$}
    \psfrag{2l-1}{\huge $2l-1$}
    \psfrag{2l+1}{\huge $2l+1$}
    \psfrag{l+1}{\huge $l+1$}
    \psfrag{3l-1}{\huge $3l-1$}
    \psfrag{3l+1}{\huge $3l+1$}
    \psfrag{2l-3}{\huge $2l-3$}
    \psfrag{l-3}{\huge $l-3$}
    \psfrag{2l-2}{\huge $2l-2$}
    \psfrag{2l+2}{\huge $2l+2$}
    \psfrag{l-1}{\huge $l-1$}
    \psfrag{l+1}{\huge $l+1$}
    \psfrag{3l-3}{\huge $3l-3$}
    \psfrag{l-2}{\huge $l-2$}
    \psfrag{4l+3}{\huge $\mathbf{4l+3}$}
    \centerline{\scalebox{.35}{\includegraphics{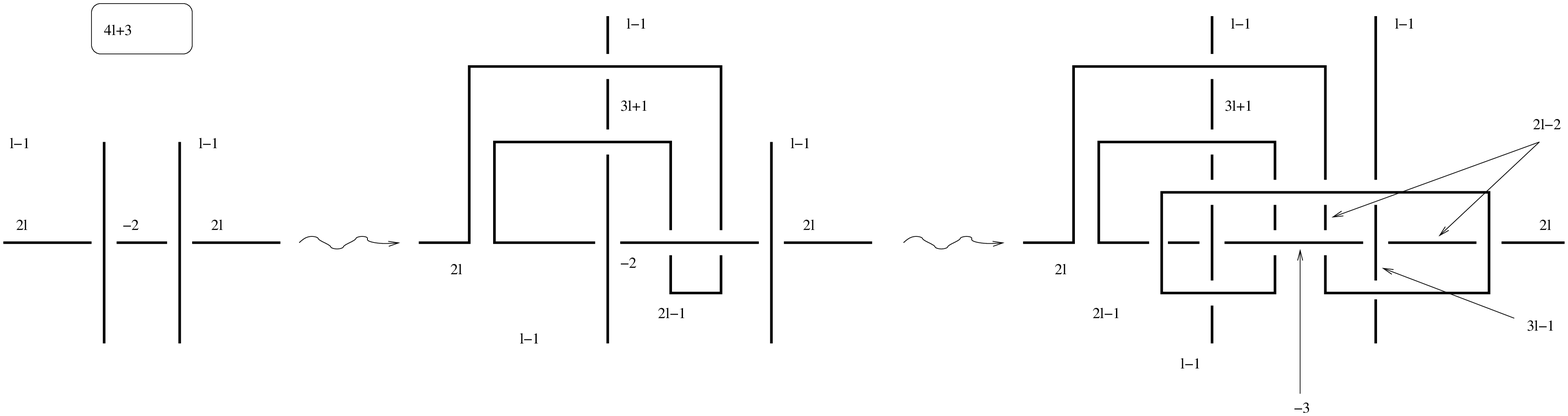}}}
    \caption{Removing color $2k-1$ from an under-arc: the $b=a$ with $a=l-1$ sub-instance}\label{fig:red16}
\end{figure}

\paragraph{\bf The $b=a$ with $a=l-1$ mod $4l+3$ instance (Figure \ref{fig:red16}).}\label{para:2k-1underbaal-14l+3}

With $a=l-1$ mod $4l+3$, $2a+2=2l$.

\begin{table}[h!]
\begin{center}
\scalebox{.87}{\begin{tabular}{| c | c |  c | c | c | c |  c | c |}\hline
2l-2=-1 & 2l-2=-2  & 2l-1=-1 & 2l-1=-2 & 3l-1=-1 & 3l-1=-2  & 3l+1=-1 & 3l+1=-2 \\ \hline
2l-2($<$3l)$<$4l+2 \, X & l=0\, X  & l=0\, X & 2l-1($<$3l)$<$4l+1\, X & l=0 \, X & 3l-1$<$4l+1 \, X & 3l+1$<$4l+2\, X  & l=0 \, X \\ \hline
    \end{tabular}}
\caption{Equalities which should not occur in Figure \ref{fig:red16} (1st row) and their consequences (2nd row). $X$'s stand for  conflicts with assumptions thus not requiring further considerations.}\label{Ta:fig:red16}
\end{center}
\end{table}

\begin{figure}[!ht]
    \psfrag{0}{\huge $0$}
    \psfrag{-2}{\huge $-2$}
    \psfrag{-1}{\huge $-1$}
    \psfrag{2l}{\huge $2l$}
    \psfrag{2l-1}{\huge $2l-1$}
    \psfrag{2l+1}{\huge $2l+1$}
    \psfrag{l+1}{\huge $l+1$}
    \psfrag{3l-1}{\huge $3l-1$}
    \psfrag{3l+1}{\huge $3l+1$}
    \psfrag{2l-3}{\huge $2l-3$}
    \psfrag{l-3}{\huge $l-3$}
    \psfrag{2l-2}{\huge $2l-2$}
    \psfrag{2l+2}{\huge $2l+2$}
    \psfrag{l-1}{\huge $l-1$}
    \psfrag{l+1}{\huge $l+1$}
    \psfrag{3l-3}{\huge $3l-3$}
    \psfrag{l-2}{\huge $l-2$}
    \psfrag{l}{\huge $l$}
    \psfrag{3l+1}{\huge $\mathbf{3l+1}$}
    \centerline{\scalebox{.35}{\includegraphics{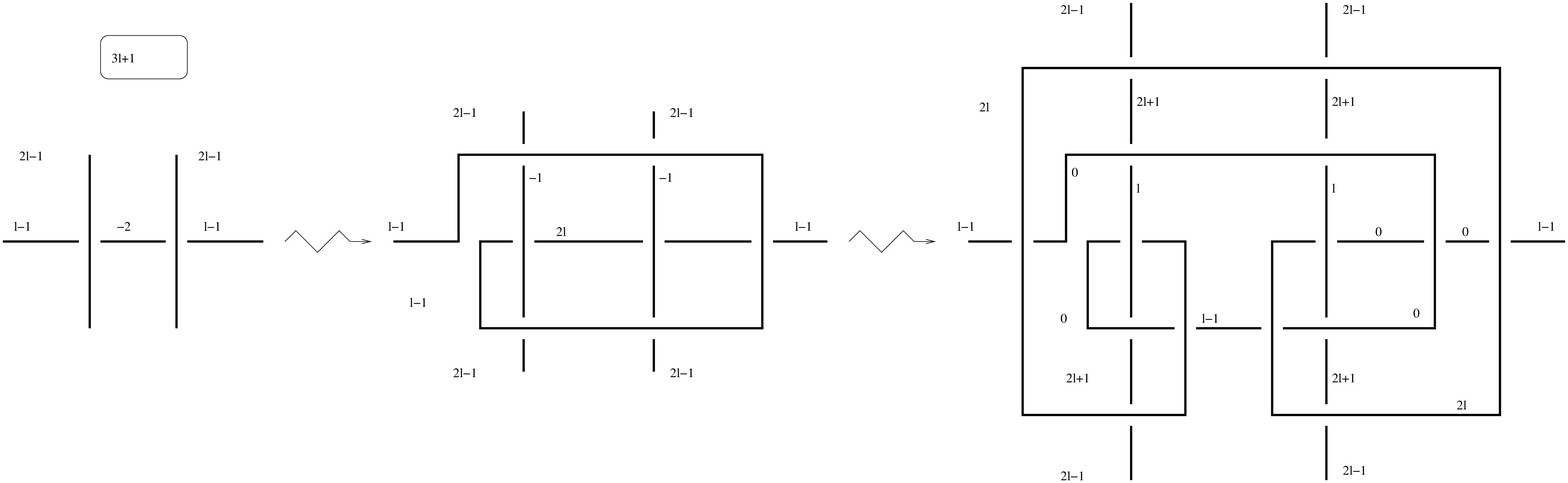}}}
    \caption{Removing color $2k-1$ from an under-arc: the $b=a$ with $a=2l-1$ sub-instance}\label{fig:red17}
\end{figure}

\paragraph{\bf The $b=a$ with $a=2l-1$ mod $3l+1$ instance (Figure \ref{fig:red17}).}\label{para:2k-1underbaa2l-1}

With $a=2l-1$ mod $3l+1$, $2a+2=4l=l-1$.

\begin{table}[h!]
\begin{center}
    \begin{tabular}{| c | c |  c | c | }\hline
l=-1 & l=-2  & 2l+1=-1 & 2l+1=-2  \\ \hline
l=0 or 2=0 \, X & l($<$2l)$<$3l-1 \, X  & l=1\, X & l=2\, X  \\ \hline
    \end{tabular}
\caption{Equalities which should not occur in Figure \ref{fig:red17} (1st row) and their consequences (2nd row). $X$'s stand for  conflicts with assumptions thus not requiring further considerations.}\label{Ta:fig:red17}
\end{center}
\end{table}

\begin{figure}[!ht]
    \psfrag{0}{\huge $0$}
    \psfrag{-2}{\huge $-2$}
    \psfrag{-1}{\huge $-1$}
    \psfrag{2l}{\huge $2l$}
    \psfrag{2l-1}{\huge $2l-1$}
    \psfrag{2l+1}{\huge $2l+1$}
    \psfrag{l+1}{\huge $l+1$}
    \psfrag{3l-1}{\huge $3l-1$}
    \psfrag{3l+1}{\huge $3l+1$}
    \psfrag{2l-3}{\huge $2l-3$}
    \psfrag{l-3}{\huge $l-3$}
    \psfrag{2l-2}{\huge $2l-2$}
    \psfrag{2l+2}{\huge $2l+2$}
    \psfrag{l-1}{\huge $l-1$}
    \psfrag{l+1}{\huge $l+1$}
    \psfrag{3l-3}{\huge $3l-3$}
    \psfrag{l-2}{\huge $l-2$}
    \psfrag{l}{\huge $l$}
    \psfrag{3l+2}{\huge $\mathbf{3l+2}$}
    \centerline{\scalebox{.35}{\includegraphics{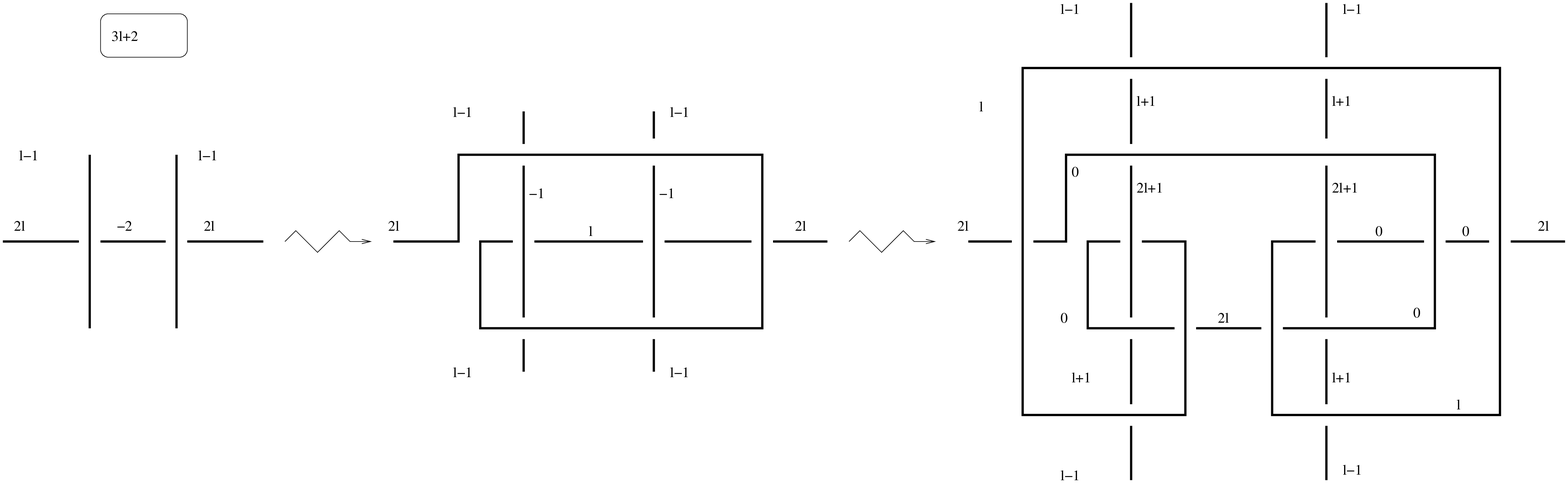}}}
    \caption{Removing color $2k-1$ from an under-arc: the $b=a$ with $a=l-1$ sub-instance}\label{fig:red18}
\end{figure}

\paragraph{\bf The $b=a$ with $a=l-1$ mod $3l+2$ instance (Figure \ref{fig:red18}).}\label{para:2k-1underbaal-13l+2}

With $a=l-1$ mod $3l+2$, $2a+2=2l$.

\begin{table}[h!]
\begin{center}
    \begin{tabular}{| c | c | c | c |  c | c | }\hline
l=-1  &  l=-2  &  l+1=-1 & l+1=-2  & 2l+1=-1 & 2l+1=-2  \\ \hline
l$<$3l+1\, X & l$<$3l\, X & l=0 or 2=0 \, X & l+1($\leq$2l)$<$3l\, X  & l=0\, X & l=1\, X  \\ \hline
    \end{tabular}
\caption{Equalities which should not occur in Figure \ref{fig:red18} (1st row) and their consequences (2nd row). $X$'s stand for  conflicts with assumptions thus not requiring further considerations.}\label{Ta:fig:red18}
\end{center}
\end{table}

This completes the proof that colors $2k-1$ and $2k$ can be removed from a non-trivial $2k+1$-coloring.

\section{Removal of color $k$ provided colors $2k$ and $2k-1$ have been removed}\label{sect:removek}

\noindent

\subsection{Removal of Color $k$ from a Monochromatic Crossing}\label{subsect:kmono}

\noindent

\begin{figure}[!ht]
    \psfrag{a}{\huge $a$}
    \psfrag{k}{\huge $k$}
    \psfrag{2k-a}{\huge $2k-a$}
    \psfrag{2a-k}{\huge $2a-k$}
    \centerline{\scalebox{.35}{\includegraphics{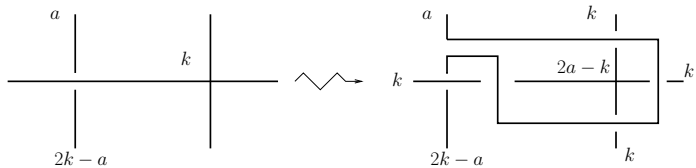}}}
    \caption{Removing color $k$ from a monochromatic crossing (1st instance).}\label{fig:49h}
\end{figure}

View Figure \ref{fig:49h}.

\begin{table}[h!]
\begin{center}
    \begin{tabular}{| c | c |  c |  }\hline
2a-k = -1 & 2a-k=-2 & 2a-k=k  \\ \hline
2a=k-1 & 2a=k-2 & a=k\, X  \\ \hline
    \end{tabular}
\caption{Equalities which should not occur in Figure \ref{fig:49h} (1st row) and their consequences (2nd row). $X$'s stand for  conflicts with assumptions thus not requiring further considerations.}\label{Ta:fig:49h}
\end{center}
\end{table}

\begin{figure}[!ht]
    \psfrag{a}{\huge $a$}
    \psfrag{2k-a}{\huge $2k-a$}
    \psfrag{-2a-2-k}{\huge $-2a-2-k$}
    \psfrag{k}{\huge $k$}
    \centerline{\scalebox{.35}{\includegraphics{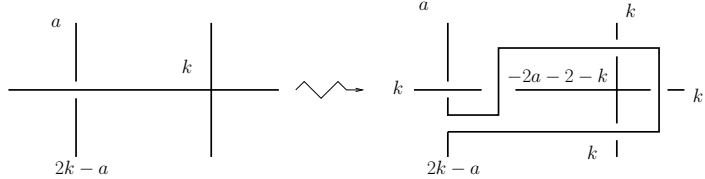}}}
    \caption{Removing color $k$ from a monochromatic crossing (2nd instance).}\label{fig:49v}
\end{figure}

\begin{table}[!ht]
\begin{center}
    \begin{tabular}{| c | c |   }\hline
2a-k = -1 & 2a-k=-2   \\ \hline
-2a-2-k=0 \, V & -2a-2-k=1 \, V  \\ \hline
    \end{tabular}
\caption{Equalities which should not occur in Figure \ref{fig:49h} (1st row) and their consequences in Figure \ref{fig:49v} (2nd row). $V$'s stand for ``admissible value'' thus not requiring further considerations. Note the different reasoning here.}\label{Ta:fig:49v}
\end{center}
\end{table}

\subsection{Removal of Color $k$ from an Over-Arc}\label{subsect:kover}

\noindent

\begin{figure}[!ht]
    \psfrag{a}{\huge $a$}
    \psfrag{3a+1}{\huge $3a+1$}
    \psfrag{2k-a}{\huge $2k-a$}
    \psfrag{2a-k}{\huge $2a-k$}
    \psfrag{k}{\huge $k$}
    \centerline{\scalebox{.35}{\includegraphics{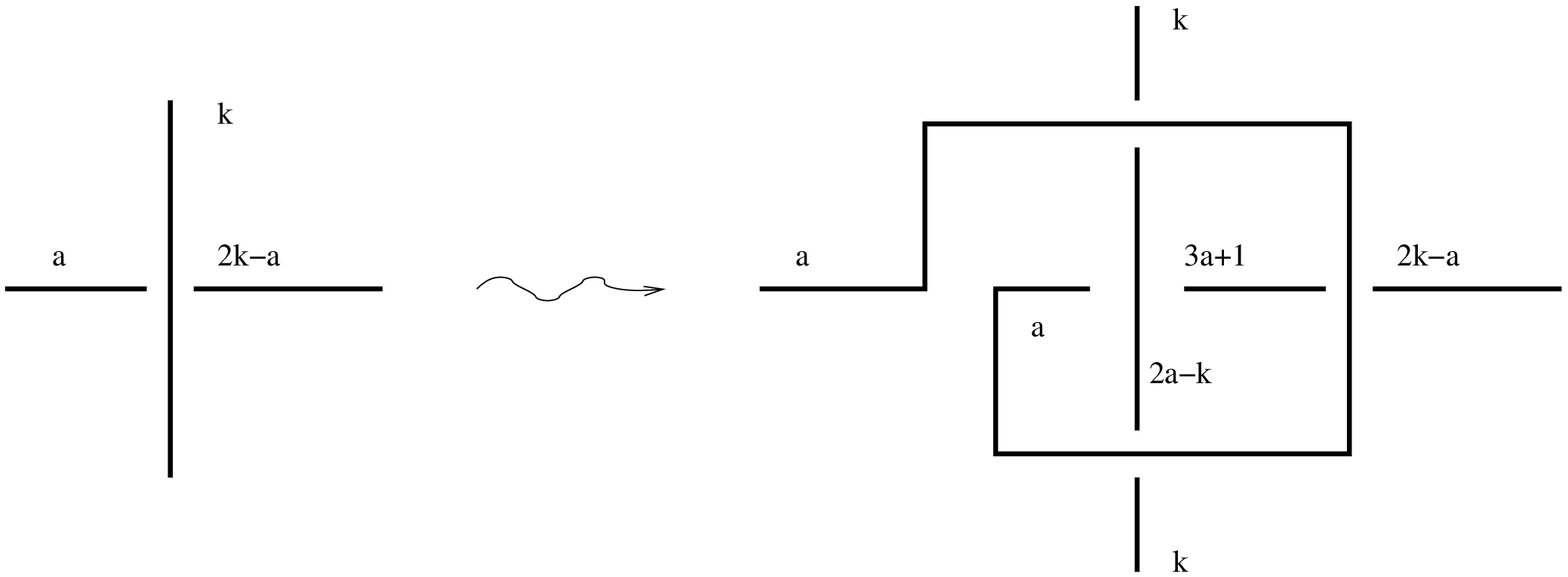}}}
    \caption{Removing color $k$ from an over-arc.}\label{fig:50-}
\end{figure}

View Figure \ref{fig:50-}.

\begin{table}[!ht]
\begin{center}
    \scalebox{.97}{\begin{tabular}{| c | c |  c | c | c | c | }\hline
2a-k=-1 & 2a-k=-2 & 2a-k=k & 3a+1=-1 &  3a+1=-2 &  3a+1=k\\ \hline
a$=_{4l+1}$3l or a$=_{4l+3}$l & a$=_{4l+1}$l-1 or a$=_{4l+3}$3l+1 & a=k \, X &  a$=_{3l+1}$2l  or a$=_{3l+2}$l  & a=-1 \, X & a=k \, X \\ \hline
    \end{tabular}}
\caption{Equalities which should not occur in Figure \ref{fig:50-} (1st row) and their consequences (2nd row). $X$'s stand for  conflicts with assumptions thus not requiring further considerations. The other situations are dealt with in \ref{para:overarck-a-k-1=-1}, \ref{para:overarck-a-k-1=-2}, \ref{para:overarc2a=k-2k=2l}, and in \ref{para:overarc2a=k-1k=2l+1}.}\label{Ta:fig:50-}
\end{center}
\end{table}

\subsubsection{The $2a=k-1$ instance.}\label{subsubsect:overarck2a=k-1}

\begin{figure}[!ht]
    \psfrag{a}{\huge $3l$}
    \psfrag{a'}{\huge $3l+1$}
    \psfrag{k}{\huge $2l$}
    \psfrag{-a-1}{\huge $l$}
    \psfrag{4l+1}{\huge $\mathbf{4l+1}$}
    \psfrag{0}{\huge $0$}
    \centerline{\scalebox{.35}{\includegraphics{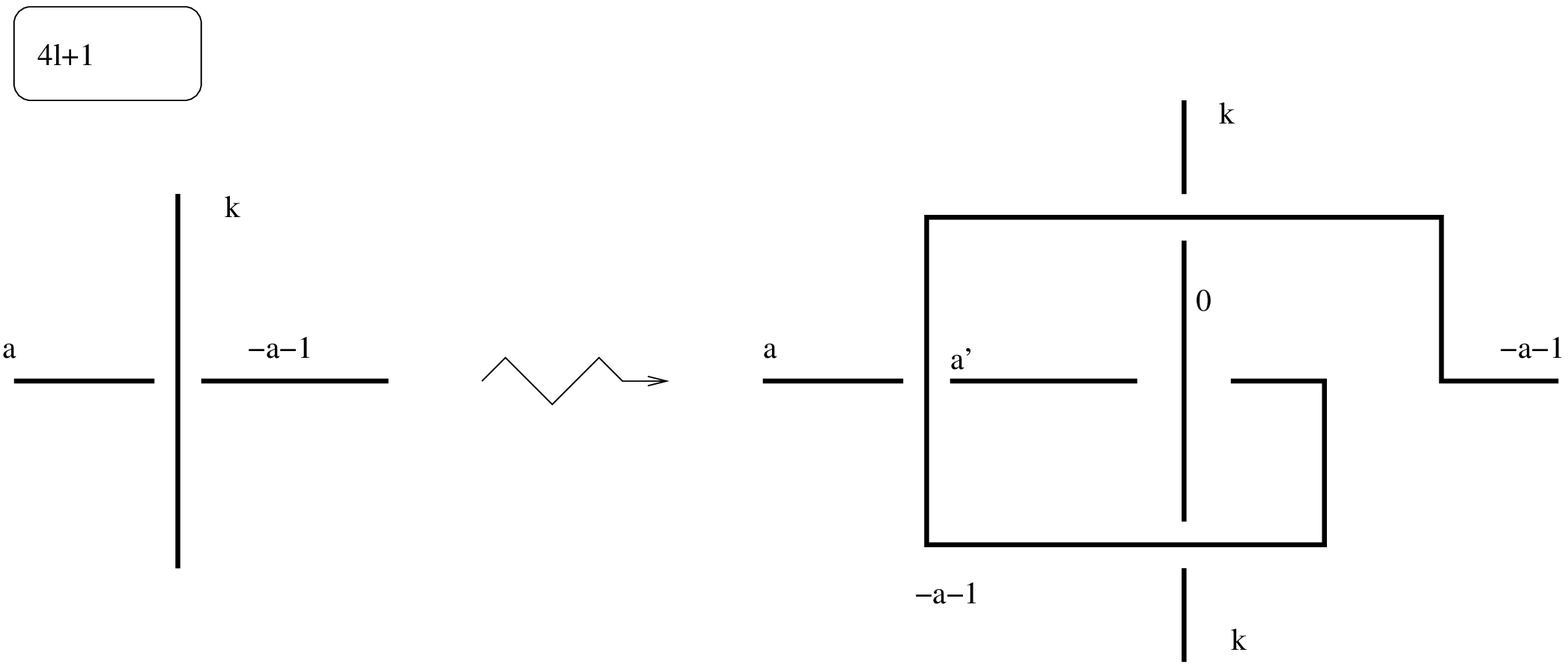}}}
    \caption{Removing color $k$ from an over-arc. The $2a=k-1$ with $k=_{4l+1}2l$ sub-instance; $(k=2l<)a=3l(<4l)$}\label{fig:50a}
\end{figure}

\paragraph{The $k=2l$ and $a=3l$, mod $4l+1$ sub-instance.}\label{para:overarck-a-k-1=-1} View Figure \ref{fig:50a}.

\begin{table}[!ht]
\begin{center}
    \begin{tabular}{| c | c |  c | c |}\hline
 3l+1=4l &  3l+1=4l-1 \, X   \\ \hline
l=1\, X & l=2 \, X  \\ \hline
    \end{tabular}
\caption{Equalities which should not occur in Figure \ref{fig:50a} (1st row) and their consequences (2nd row). $X$'s stand for  conflicts with assumptions thus not requiring further considerations.}\label{Ta:fig:red50a}
\end{center}
\end{table}

\begin{figure}[!ht]
    \psfrag{a}{\huge $l$}
    \psfrag{a'}{\huge $l+1$}
    \psfrag{k}{\huge $2l+1$}
    \psfrag{-a-1}{\huge $3l+2$}
    \psfrag{4l+1}{\huge $\mathbf{4l+3}$}
    \psfrag{0}{\huge $0$}
    \centerline{\scalebox{.35}{\includegraphics{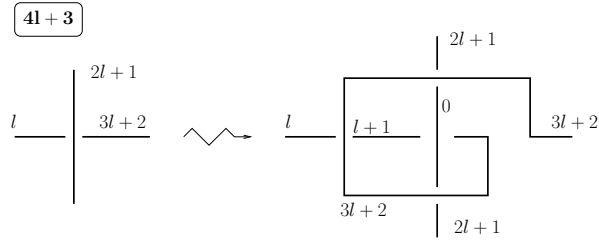}}}
    \caption{Removing color $k$ from an over-arc. The $2a=k-1$ with $k=_{4l+3}2l+1$ sub-instance; $(0<)a=l(<k=2l+1<4l+1<4l+2)$.}\label{fig:50b}
\end{figure}

\paragraph{The $k=2l+1$ and $a=l$, mod $4l+3$ sub-instance.}\label{para:overarck-a-k-1=-2}View Figure \ref{fig:50b}.

\begin{table}[!ht]
\begin{center}
    \begin{tabular}{| c | }\hline
$0<\quad l+1\quad <2l+1<3l+1<4l+1<4l+2$ \\ \hline
    \end{tabular}
\caption{ $l+1$ can only  take on admissible values, in Figure \ref{fig:50b}.}\label{Ta:fig:red50b}
\end{center}
\end{table}

\subsubsection{The $2a=k-2$ instance.}\label{subsubsect:overarck2a=k-2}

\begin{figure}[!ht]
    \psfrag{a}{\huge $l-1$}
    \psfrag{a'}{\huge $l+2$}
    \psfrag{k}{\huge $2l$}
    \psfrag{-a-1}{\huge $3l+1$}
    \psfrag{4l+1}{\huge $\mathbf{4l+1}$}
    \psfrag{0}{\huge $1$}
    \centerline{\scalebox{.35}{\includegraphics{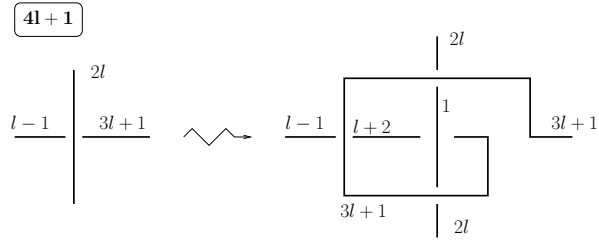}}}
    \caption{Removing color $k$ from an over-arc. The $2a=k-2$ with $k=_{4l+1}2l$ sub-instance; $a=l-1(<k=2l<4l-1<4l).$}\label{fig:50c}
\end{figure}

\paragraph{The $k=2l$ and $a=l-1$, mod $4l+1$ sub-instance.}\label{para:overarc2a=k-2k=2l} View Figure \ref{fig:50c}.

\bigbreak

\begin{table}[!ht]
\begin{center}
    \begin{tabular}{| c | c | c | }\hline
 l+2 = 2l & l+2 = 4l-1 & l+2 = 4l \\ \hline
 4l+1 = 9 \, X & 4l+1 = 5 \, X & 3l=2\, X\\ \hline
    \end{tabular}
\caption{ $l+2$ can only  take on admissible values, in Figure \ref{fig:50c}.}\label{Ta:fig:red50c}
\end{center}
\end{table}

\bigbreak

\begin{figure}[!ht]
    \psfrag{a}{\huge $3l+1$}
    \psfrag{a'}{\huge $3l+4$}
    \psfrag{k}{\huge $2l+1$}
    \psfrag{-a-1}{\huge $l+1$}
    \psfrag{4l+1}{\huge $\mathbf{4l+3}$}
    \psfrag{0}{\huge $1$}
    \centerline{\scalebox{.35}{\includegraphics{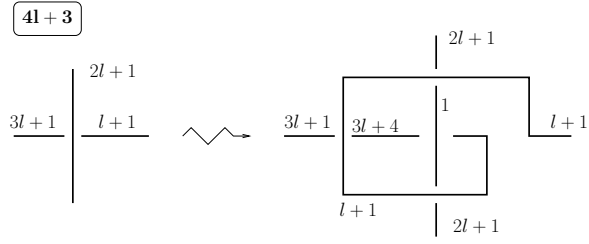}}}
    \caption{Removing color $k$ from an over-arc. The $2a=k-2$ with $k=_{4l+3}2l+1$ sub-instance; $(k=2l+1<)a=3l+1(<4l+1<4l+2)$.}\label{fig:50d}
\end{figure}

\paragraph{The $k=2l+1$ and $a=3l+1$, mod $4l+3$ sub-instance.}\label{para:overarc2a=k-1k=2l+1}. View Figure \ref{fig:50d}.

\bigbreak

\bigbreak

\begin{table}[!ht]
\begin{center}
    \begin{tabular}{| c | c |  c | c |}\hline
$2l+1<3l+4$ &    3l+4=4l+1 &  3l+4=4l+2 \, X   \\ \hline
$2l+1<3l+4$\, X & 4l+3=15 \, X & l=2 and 4l+3=11   \\ \hline
    \end{tabular}
\caption{Equalities which should not occur in Figure \ref{fig:50d} (1st row) and their consequences (2nd row). $X$'s stand for  conflicts with assumptions thus not requiring further considerations.}\label{Ta:fig:red50d}
\end{center}
\end{table}

\bigbreak

\begin{figure}[!ht]
    \psfrag{3}{\huge $3$}
    \psfrag{1}{\huge $1$}
    \psfrag{-1}{\huge $-1$}
    \psfrag{2}{\huge $2$}
    \psfrag{4}{\huge $4$}
    \psfrag{5}{\huge $5$}
    \psfrag{7}{\huge $7$}
    \psfrag{11}{\huge $\mathbf{11}$}
    \psfrag{0}{\huge $0$}
    \centerline{\scalebox{.35}{\includegraphics{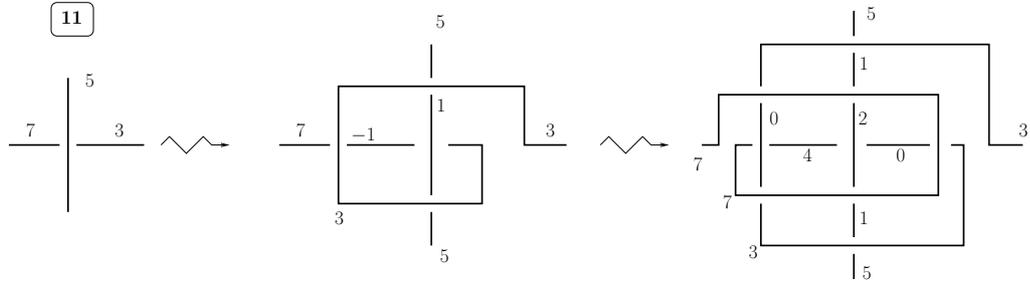}}}
    \caption{Removing color $k=5\, (\text{mod } 11)$  from an over-arc. The $2a=k-2$ with $k=_{4l+3}2l+1$ and $l=2$ sub-sub-instance.}\label{fig:50e}
\end{figure}

\begin{itemize}\label{para:overarck-a-k-1=-2sub}
\item The $l=2$ sub-instance. View Figure \ref{fig:50e}
\end{itemize}
\bigbreak

\bigbreak

\subsubsection{The $3a=-2$ instance.}\label{subsubsect:overarck3a=-2}

\begin{figure}[!ht]
    \psfrag{a}{\huge $4l$}
    \psfrag{a'}{\huge $0$}
    \psfrag{k}{\huge $3l$}
    \psfrag{-a-1}{\huge $2l$}
    \psfrag{4l+1}{\huge $\mathbf{6l+1}$}
    \psfrag{0}{\huge $l$}
    \centerline{\scalebox{.35}{\includegraphics{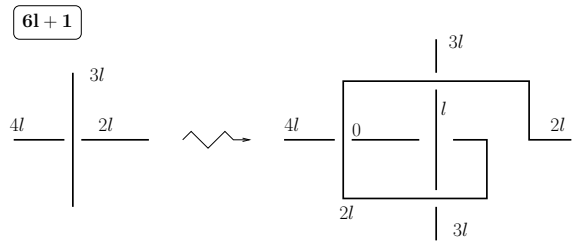}}}
    \caption{Removing color $k$ from an over-arc. The $3a=-2$ with $k=_{6l+1}3l$ sub-instance; $(k=3l<)a=4l(<6l-1<6l).$}\label{fig:50f}
\end{figure}

\paragraph{The $k=3l$ and $a=4l$ mod $6l+1$ sub-instance.}\label{para:overarc3a=-2k-a-k-1=-1} View Figure \ref{fig:50f}.

\bigbreak

\begin{table}[!ht]
\begin{center}
    \begin{tabular}{| c |  }\hline
 $l< 3l < 6l-1 < 6l$ \\ \hline
    \end{tabular}
\caption{ $l$ can only  take on admissible values, in Figure \ref{fig:50f}.}\label{Ta:fig:red50f}
\end{center}
\end{table}

\begin{figure}[!ht]
    \psfrag{a}{\huge $2l+1$}
    \psfrag{a'}{\huge $0$}
    \psfrag{k}{\huge $3l+2$}
    \psfrag{-a-1}{\huge $4l+3$}
    \psfrag{4l+1}{\huge $\mathbf{6l+5}$}
    \psfrag{0}{\huge $5l+4$}
    \centerline{\scalebox{.35}{\includegraphics{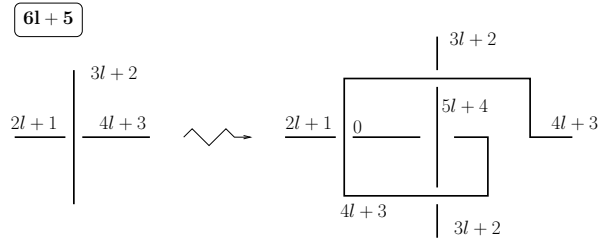}}}
    \caption{Removing color $k$ from an over-arc. The $3a=-2$ with $k=_{6l+5}3l+2$ sub-instance; $a=2l+1(<k=3l+2<6l+3<6l+4)$.}\label{fig:50g}
\end{figure}

\paragraph{The $k=3l+2$ and $a=2l+1$, mod $6l+5$ sub-instance.}\label{para:overarc3a=-2k=3l+26l+5} View Figure \ref{fig:50g}.

\bigbreak

\begin{table}[!ht]
\begin{center}
    \begin{tabular}{| c| c | c | }\hline
3l+2$<$5l+4 &  5l+4 = 6l+3 & $5l+4<6l+4$    \\ \hline
3l+2$<$5l+4 \, X  &  6l+5=11 & $5l+4<6l+4$ \, X   \\ \hline
    \end{tabular}
\caption{Equalities which should not occur in Figure \ref{fig:50g} (1st row) and their consequences (2nd row). The $6l+5=11$ ($l=1$) instance is solved in Figure \ref{fig:50e}.}\label{Ta:fig:red50g}
\end{center}
\end{table}

\subsection{Removal of Color $k$ from an Under-Arc}\label{subsect:kunder}

\noindent

\begin{figure}[!ht]
    \psfrag{a}{\huge $a$}
    \psfrag{b}{\huge $b$}
    \psfrag{k}{\huge $k$}
    \psfrag{2a-2b+k}{\huge $2a-2b+k$}
    \psfrag{2a-k}{\huge $2a-k$}
    \psfrag{2a-b}{\huge $2a-b$}
    \psfrag{2b-k}{\huge $2b-k$}
    \centerline{\scalebox{.35}{\includegraphics{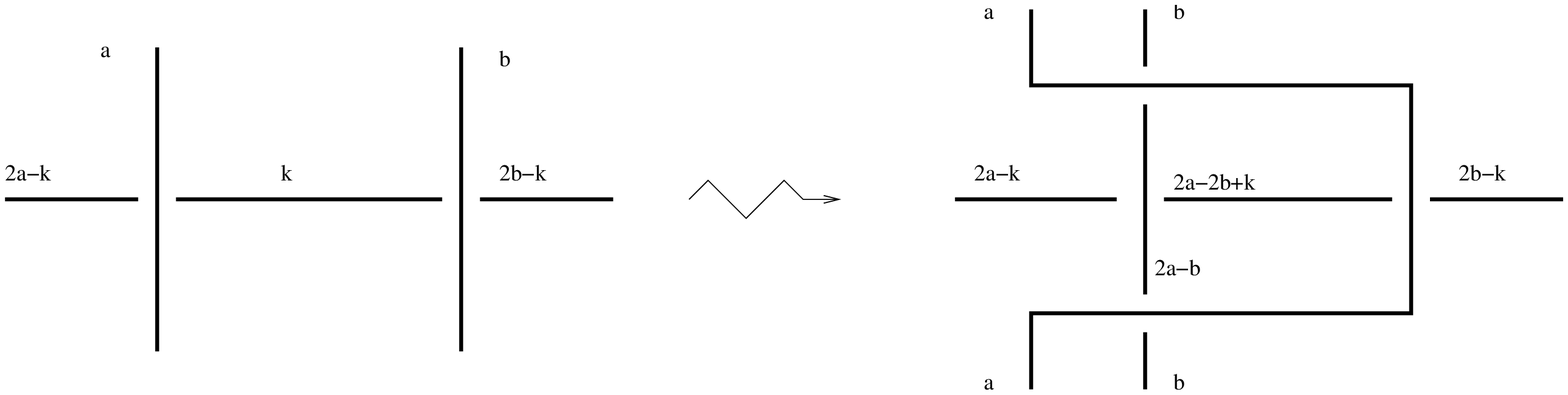}}}
    \caption{Removing color $k$ from an under-arc}\label{fig:54}
\end{figure}

View Figure \ref{fig:54}.

\begin{table}[!ht]
\begin{center}
    \scalebox{.9}{\begin{tabular}{| c | c |  c | c | c | c |}\hline
2a-b=-1 & 2a-b=-2 & 2a-b=k &  2a-2b+k=-1  & 2a-2b+k=-2 & 2a-2b+k=k \\ \hline
b=2a+1 & b=2a+2  &  b=2a-k  &   b=$_{4l+1}$a+3l+1 or b=$_{4l+3}$a+l+1  &   b=$_{4l+1}$a+l+1 or b=$_{4l+3}$a+3l+3  & b=a \\ \hline
    \end{tabular}}
\caption{Equalities which should not occur in Figure \ref{fig:54} (1st row) and their consequences (2nd row). $X$'s stand for  conflicts with assumptions thus not requiring further considerations. The other instances are dealt with in \ref{subsubsectb2a+1}, \ref{subsubsectb2a+2}, \ref{subsubsectb2a-k}, \ref{subsubsect2(b-a)=k+1},  \ref{subsubsect2(b-a)=k+2},  \ref{subsubsectb=a}   }\label{Ta:fig:red54}
\end{center}
\end{table}

\subsubsection{The $b=2a+1$ instance}\label{subsubsectb2a+1} 
\begin{figure}[!ht]
    \psfrag{a}{\huge $a$}
    \psfrag{2a-k}{\huge $2a-k$}
    \psfrag{2a+2+k}{\huge $2a+2+k$}
    \psfrag{2a+1}{\huge $2a+1$}
    \psfrag{4a+2-k}{\huge $4a+2-k$}
    \psfrag{k}{\huge $k$}
    \psfrag{3a+2}{\huge $3a+2$}
    \centerline{\scalebox{.35}{\includegraphics{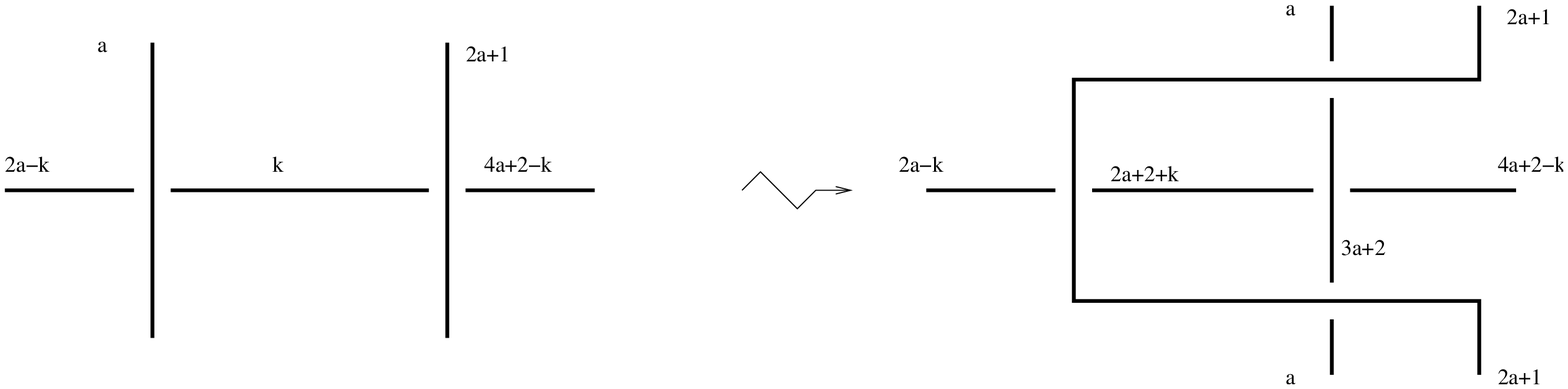}}}
    \caption{Removing color $k$ from an under-arc: the $b=2a+1$ instance}\label{fig:55}
\end{figure}

View Figure \ref{fig:55}.

\begin{table}[!ht]
\begin{center}
    \scalebox{.92}{\begin{tabular}{ | c | c |  c | c | c | c |}\hline
2a+2+k=-1 & 2a+2+k=-2 & 2a+2+k=k &  3a+2=-1  & 3a+2=-2 & 3a+2=k \\ \hline
2a-k=-2 \, X & a$=_{4l+1}$3l-1 or a$=_{4l+3}$l-1  &  a=-1 \, X  &   a=-1 \, X & a$=_{3l+1}$l-1 or a$=_{3l+2}$2l &  a$=_{6l+1}$5l or a$=_{6l+5}$l   \\ \hline
    \end{tabular}}
\caption{Equalities which should not occur in Figure \ref{fig:55} (1st row) and their consequences (2nd row). $X$'s stand for  conflicts with assumptions thus not requiring further considerations. The other instances are dealt with in \ref{parab2a+12a=k-3k=2l}, \ref{parab2a+12a=k-3k=2l+1}, \ref{parab2a+13a=k-2k=3l}, \ref{parab2a+13a=k-2k=3l+2}, \ref{parab2a+13a=-4k=3l}, and \ref{parab2a+13a=k-2k=3l+2a=4l+2}.  }\label{Ta:fig:red55}
\end{center}
\end{table}

\begin{figure}[!ht]
    \psfrag{-3}{\huge $-3$}
    \psfrag{2l-5}{\huge $2l-5$}
    \psfrag{2l}{\huge $2l$}
    \psfrag{2l-4}{\huge $2l-4$}
    \psfrag{l-4}{\huge $l-4$}
    \psfrag{2l-1}{\huge $2l-1$}
    \psfrag{2l-3}{\huge $2l-3$}
    \psfrag{3l-1}{\huge $3l-1$}
    \psfrag{2l-2}{\huge $2l-2$}
    \psfrag{4l+1}{\huge $\mathbf{4l+1}$}
    \centerline{\scalebox{.35}{\includegraphics{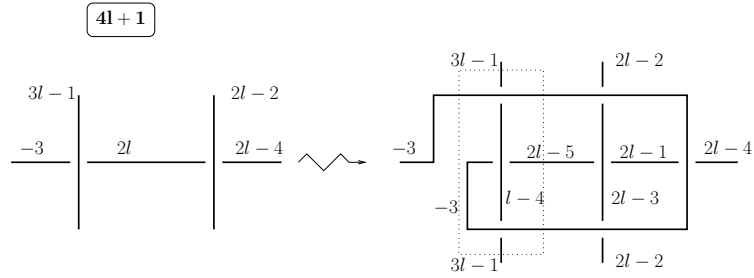}}}
    \caption{Removing color $k$ from an under-arc: the $b=2a+1$, $2a=k-3$; $k=2l$, $a=3l-1$, mod $4l+1$ sub-instance. Note that $k(=2l)\leq a (=3l-1)<4l-1<4l$}\label{fig:72aa}
\end{figure}

\paragraph{The $b=2a+1$, $2a=k-3$; $k=2l$, $a=3l-1$, mod $4l+1$ sub-instance}\label{parab2a+12a=k-3k=2l} View Figure \ref{fig:72aa}.

\begin{table}[!ht]
\begin{center}
    \scalebox{.72}{\begin{tabular}{| c | c |  c | c | c | c |c | c | c |c | c | c |}\hline
l-4 =2l & l-4=-1 & l-4=-2  & 2l-5=2l   & 2l-5=-1 &  2l-5=-2 & 2l-3=2l & 2l-3=-1 & 2l-3=-2 & 2l-1 = 2l & 2l-1=-1l=0 2l-1=-2\\ \hline
l-4$<$2l\, X & l=3  & 4l+1=9 \, X  & 2l-5$<$2l & 4l+1=9 \, X  &  2l=3 \, X & 2l-3$<$2l \, X &  4l+1=5 \, X  &  2l=1 \, X  & 2l-1 $<$ 2l\, X & \, X  2l=-1\, X\\ \hline
    \end{tabular}}
\caption{Equalities which should not occur in Figure \ref{fig:72aa} (1st row) and their consequences (2nd row). An inequality in the second row stands for an obvious statement which contradicts the equality above it. $X$'s stand for  conflicts with assumptions thus not requiring further considerations. The region boxed by the dotted lines contains a forbidden color for $l=3$ which is dealt with in Figure \ref{fig:72aaa}}\label{Ta:fig:red72aa}
\end{center}
\end{table}

\begin{figure}[!ht]
    \psfrag{-3}{\huge $-3$}
    \psfrag{8}{\huge $8$}
    \psfrag{-1}{\huge $-1$}
    \psfrag{1}{\huge $1$}
    \psfrag{7}{\huge $7$}
    \psfrag{5}{\huge $5$}
    \psfrag{3}{\huge $3$}
    \psfrag{13}{\huge $\mathbf{13}$}
    \centerline{\scalebox{.35}{\includegraphics{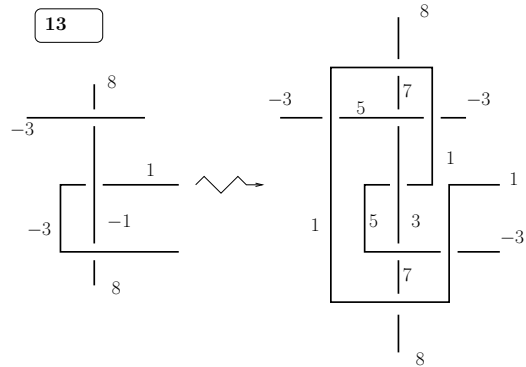}}}
    \caption{Removing color $k$ from an under-arc: the $b=2a+1$, $2a=k-3$; $k=2l$, $a=3l-1$, mod $4l+1$  sub-instance: fixing the boxed region in Figure \ref{fig:72aa} for $l=3$.}\label{fig:72aaa}
\end{figure}

\bigbreak

\begin{figure}[!ht]
    \psfrag{3l-2}{\huge $3l-2$}
    \psfrag{2l-1}{\huge $2l-1$}
    \psfrag{2l-4}{\huge $2l-4$}
    \psfrag{2l+1}{\huge $2l+1$}
    \psfrag{2l-3}{\huge $2l-3$}
    \psfrag{2l-2}{\huge $2l-2$}
    \psfrag{l-1}{\huge $l-1$}
    \psfrag{l-5}{\huge $l-5$}
    \psfrag{-3}{\huge $-3$}
    \psfrag{2l}{\huge $2l$}
    \psfrag{4l+3}{\huge $\mathbf{4l+3}$}
    \centerline{\scalebox{.35}{\includegraphics{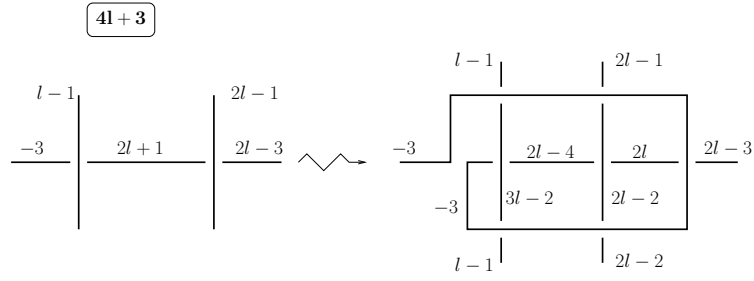}}}
    \caption{Removing color $k$ from an under-arc: the $b=2a+1$, $2a=k-3$; $k=2l+1$, $a=l-1$ mod $4l+3$ sub-instance. Note that $a(=l-1)<k(=2l+1)<4l+1<4l+2$}\label{fig:73a}
\end{figure}

\paragraph{The $b=2a+1$, $2a=k-3$; $k=2l+1$, $a=l-1$ mod $4l+3$ sub-instance.}\label{parab2a+12a=k-3k=2l+1} View Figure \ref{fig:73a}.

\begin{table}[!ht]
\begin{center}
 \begin{tabular}{| c | c  |  c  |  c  |}\hline
$2l-4<2l-2<2l<k(=2l+1) \leq 3l-2 <4l+1<4l+2$   \\ \hline
\end{tabular}
\caption{$l\geq 2$. The equality $2l+1 \leq 3l-2$ implies $4l+3=15$ which conflicts with the assumptions. There are no further considerations to be made about Figure \ref{fig:73a}.}\label{Ta:fig:red73a}
\end{center}
\end{table}
\bigbreak

\begin{figure}[!ht]
    \psfrag{5l}{\huge $5l$}
    \psfrag{4l}{\huge $4l$}
    \psfrag{5l-1}{\huge $5l-1$}
    \psfrag{3l-1}{\huge $3l-1$}
    \psfrag{3l}{\huge $3l$}
    \psfrag{4l-1}{\huge $4l-1$}
    \psfrag{l-1}{\huge $l-1$}
    \psfrag{l}{\huge $l$}
    \psfrag{6l+1}{\huge $\mathbf{6l+1}$}
    \centerline{\scalebox{.37}{\includegraphics{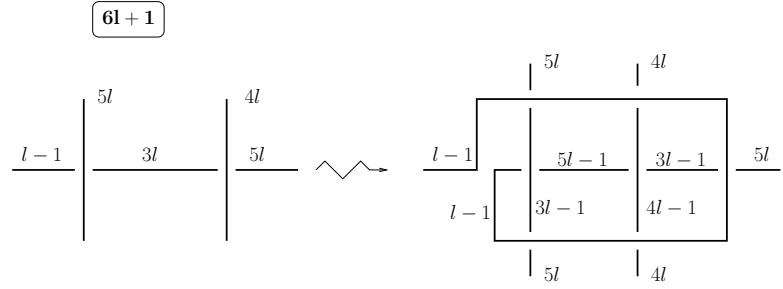}}}
    \caption{Removing color $k$ from an under-arc: the $b=2a+1$, $3a=k-2$; $k=3l$, $a=5l$ mod $6l+1$ sub-instance. $(k=3l<)a(=5l)\leq 6l-1<6l$.}\label{fig:70}
\end{figure}

\paragraph{The $b=2a+1$, $3a=k-2$; $k=3l$, $a=5l$ mod $6l+1$ sub-instance.}\label{parab2a+13a=k-2k=3l} View Figure \ref{fig:70}.

\begin{table}[!ht]
\begin{center}
 \begin{tabular}{| c |}\hline
$3l-1 < k(=3l) \leq 4l-1 < 5l-1 < 6l-1 < 6l$   \\ \hline
\end{tabular}
\caption{The equality $3l=4l-1$ leads to $6l+1=7$ which conflicts with the assumptions. There are no further considerations to be made about Figure \ref{fig:70}.}\label{Ta:fig:red70}
\end{center}
\end{table}

\bigbreak

\begin{figure}[!ht]
    \psfrag{2l+1}{\huge $2l+1$}
    \psfrag{3l+1}{\huge $3l+1$}
    \psfrag{2l+2}{\huge $2l+2$}
    \psfrag{3l+2}{\huge $3l+2$}
    \psfrag{2l}{\huge $2l$}
    \psfrag{l-1}{\huge $l-1$}
    \psfrag{l+1}{\huge $l+1$}
    \psfrag{5l+3}{\huge $5l+3$}
    \psfrag{6l+5}{\huge $\mathbf{6l+5}$}
    \psfrag{l}{\huge $l$}
    \centerline{\scalebox{.35}{\includegraphics{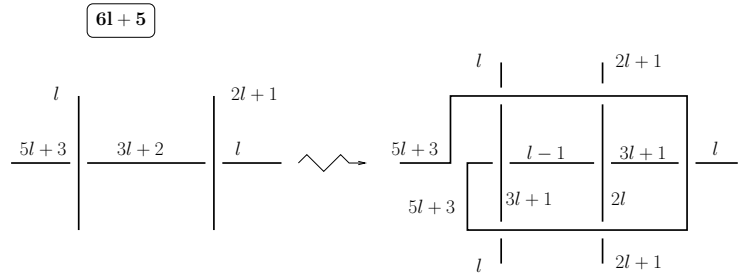}}}
    \caption{Removing color $k$ from an under-arc: the $b=2a+1$, $3a=k-2$, $k=3l$ sub-instance. $a=l(<k=3l+2<6l+3<6l+4)$.}\label{fig:71}
\end{figure}

\paragraph{The $b=2a+1$, $3a=k-2$; $k=3l+2$, $a=l$ mod $6l+5$ sub-instance}\label{parab2a+13a=k-2k=3l+2} View Figure \ref{fig:71}.

\begin{table}[!ht]
\begin{center}
    \begin{tabular}{| c |}\hline
$l-1<2l<3l+1<k(=3l+2)<6l+3<6l+4$ \\ \hline
    \end{tabular}
\caption{$l\geq 1$ and there are no further considerations to be made about Figure \ref{fig:71}.}\label{Ta:fig:red71}
\end{center}
\end{table}

\bigbreak

\begin{figure}[!ht]
    \psfrag{l-2}{\huge $l-2$}
    \psfrag{2l-1}{\huge $2l-1$}
    \psfrag{4l-1}{\huge $4l-1$}
    \psfrag{3l}{\huge $3l$}
    \psfrag{5l-2}{\huge $5l-2$}
    \psfrag{3l-3}{\huge $3l-3$}
    \psfrag{2l-4}{\huge $2l-4$}
    \psfrag{l-5}{\huge $l-5$}
    \psfrag{-4}{\huge $-4$}
    \psfrag{6l+1}{\huge $\mathbf{6l+1}$}
    \centerline{\scalebox{.35}{\includegraphics{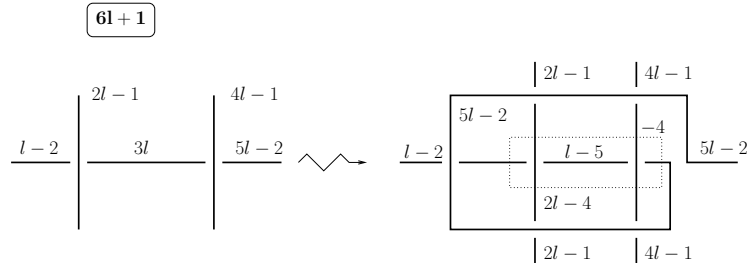}}}
    \caption{Removing color $k$ from an under-arc: the $b=2a+1$, $3a=k-2$, $k=3l$ sub-instance. $a=2l-1(<k=3l<6l-1<6l)$. The boxed region is fixed in Figure \ref{fig:70aabis} for $l=3$.}\label{fig:70aa}
\end{figure}

\paragraph{The $b=2a+1$, $3a=-4$; $k=3l$, $a=2l-1$ mod $6l+1$ sub-instance}\label{parab2a+13a=-4k=3l} View Figure \ref{fig:70aa}.

\bigbreak

\bigbreak

\bigbreak

\begin{table}[!ht]
\begin{center}
    \begin{tabular}{| c |}\hline
$l-5 < 2l-4 < 3l-3 < k(=3l) $       \\ \hline
    \end{tabular}
\caption{$l\geq2$ and the equality $l-5=-1$ leads to an inconsistency whereas $l-5=-2$ leads to the $l=3$ instance which is dealt with in Figure \ref{fig:71aabis}. Nothing else gives rise to further considerations on Figure \ref{fig:70aa}.}\label{Ta:fig:70aa}
\end{center}
\end{table}

\begin{figure}[!ht]
    \psfrag{6}{\huge $6$}
    \psfrag{2}{\huge $2$}
    \psfrag{-4}{\huge $-4$}
    \psfrag{-2}{\huge $-2$}
    \psfrag{13}{\huge $13$}
    \psfrag{8}{\huge $8$}
    \psfrag{10}{\huge $10$}
    \psfrag{19}{\huge $\mathbf{19}$}
    \centerline{\scalebox{.35}{\includegraphics{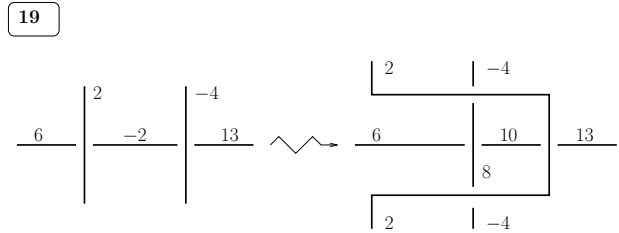}}}
    \caption{Removing color $k$ from an under-arc: the $b=2a+1$, $3a=k-2$, $k=3l$ sub-instance: fixing the boxed region in Figure \ref{fig:70aa} for  $l=3$.}\label{fig:70aabis}
\end{figure}

\paragraph{The $b=2a+1$, $3a=-4$, $k=3l$, $l=3$ sub-instance}\label{parab2a+13a=-4k=3ll=2} View Figure \ref{fig:70aabis}.

\bigbreak

\bigbreak

\bigbreak
\bigbreak

\begin{figure}[!ht]
    \psfrag{2l}{\huge $2l$}
    \psfrag{3l+1}{\huge $3l+1$}
    \psfrag{l-2}{\huge $l-2$}
    \psfrag{l-3}{\huge $l-3$}
    \psfrag{3l+2}{\huge $3l+2$}
    \psfrag{2l-1}{\huge $2l-1$}
    \psfrag{4l+2}{\huge $4l+2$}
    \psfrag{5l-1}{\huge $5l-1$}
    \psfrag{-3}{\huge $-3$}
    \psfrag{5l+2}{\huge $5l+2$}
    \psfrag{6l+5}{\huge $\mathbf{6l+5}$}
    \psfrag{l}{\huge $4l+2$}
    \centerline{\scalebox{.35}{\includegraphics{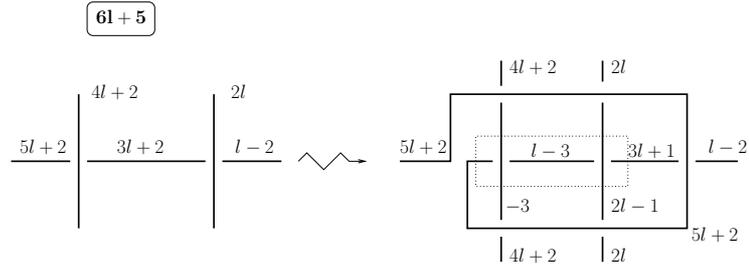}}}
    \caption{Removing color $k$ from an under-arc: the $b=2a+1$, $3a=k-2$; $k=3l+2$ $a=4l+2$ mod $6l+5$ sub-instance. $(k=3l+2<)a=4l+2(<6l+3<6l+4)$.  The boxed region is fixed in Figure \ref{fig:71aabis} for $l=2$.}\label{fig:71aa}
\end{figure}

\paragraph{The $b=2a+1$, $3a=k-2$; $k=3l+2$, $a=4l+2$ mod $6l+5$ sub-instance}\label{parab2a+13a=k-2k=3l+2a=4l+2} View Figure \ref{fig:71aa}.

\bigbreak

\bigbreak

\begin{table}[!ht]
\begin{center}
    \begin{tabular}{| c |}\hline
$l-3 < 2l-1 < 3l+1 < k(=3l+2)<6l+3<6l+4  $      \\ \hline
    \end{tabular}
\caption{$l\geq1$ and the equality $l-3=-2$ leads to an inconsistency $l-2=-1$, whereas $l-3=-1$ leads to the $l=2$ instance which is dealt with in Figure \ref{fig:71aabis}. Nothing else gives rise to further considerations on Figure \ref{fig:71aa}.}\label{Ta:fig:red71aa}
\end{center}
\end{table}

\begin{figure}[!ht]
    \psfrag{-3}{\huge $-3$}
    \psfrag{12}{\huge $12$}
    \psfrag{-1}{\huge $-1$}
    \psfrag{3}{\huge $3$}
    \psfrag{7}{\huge $7$}
    \psfrag{11}{\huge $11$}
    \psfrag{9}{\huge $9$}
    \psfrag{17}{\huge $\mathbf{17}$}
    \centerline{\scalebox{.35}{\includegraphics{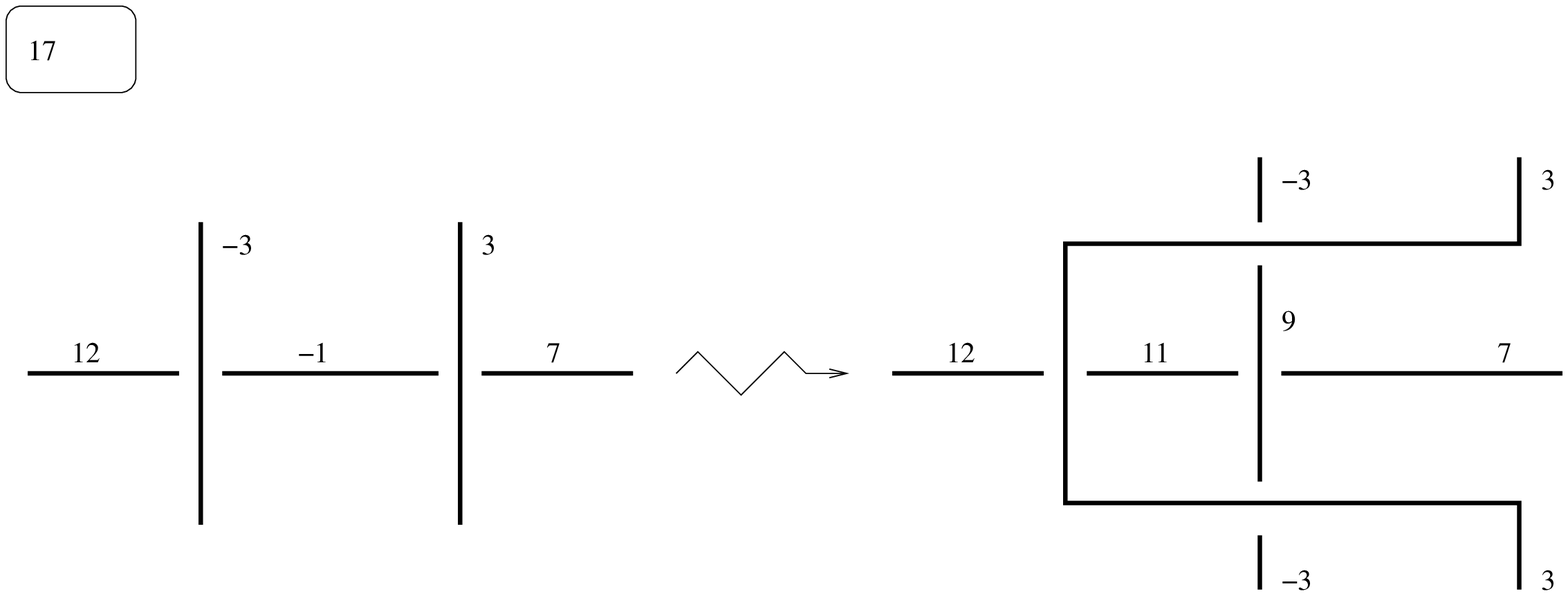}}}
    \caption{Removing color $k$ from an under-arc: the $b=2a+1$, $3a=k-2$; $k=3l+2$ $a=4l+2$ mod $6l+5$ sub-instance: fixing the boxed region in Figure \ref{fig:71aa} for  $l=2$.}\label{fig:71aabis}
\end{figure}

\bigbreak

\bigbreak

\bigbreak

\subsubsection{The $b=2a+2$ instance}\label{subsubsectb2a+2} 
\begin{figure}[!ht]
    \psfrag{a}{\huge $a$}
    \psfrag{2a-k}{\huge $2a-k$}
    \psfrag{2a+2+k}{\huge $2a+4+k$}
    \psfrag{2a+1}{\huge $2a+2$}
    \psfrag{4a+2-k}{\huge $4a+4-k$}
    \psfrag{k}{\huge $k$}
    \psfrag{3a+2}{\huge $3a+4$}
    \centerline{\scalebox{.35}{\includegraphics{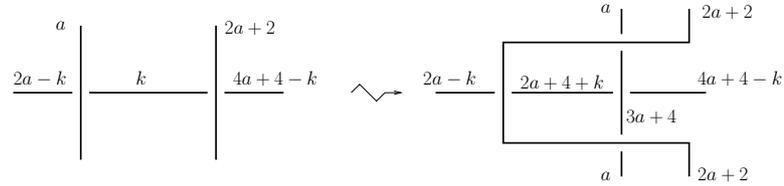}}}
    \caption{Removing color $k$ from an under-arc: the $b=2a+2$ instance}\label{fig:552a+2}
\end{figure}

View Figure \ref{fig:552a+2}.

\bigbreak

\bigbreak

\bigbreak

\bigbreak

\begin{table}[!ht]
\begin{center}
   \scalebox{.85}{\begin{tabular}{| c | c |  c | c | c | c |}\hline
2a+4+k=-1 & 2a+4+k=-2 & 2a+4+k=k &  3a+4=-1  & 3a+4=-2 & 3a+4=k(=3k+1) \\ \hline
a=$_{4l+1}$l-2 or a=$_{4l+3}$3l  & a=$_{4l+1}$3l-2 or a=$_{4l+3}$l-2 &  a=-2 \, X  &   a=$_{6l+1}$4l-1 or a=$_{6l+5}$2l& a=-2 \, X  &   2a+2=-1  \, X  \\ \hline
    \end{tabular}}
\caption{Equalities which should not occur in Figure \ref{fig:552a+2} (1st row) and their consequences (2nd row). $X$'s stand for  conflicts with assumptions thus not requiring further considerations. The other situations are dealt with in \ref{parab2a+22a=k-4k=2l}, \ref{parab2a+22a=k-4k=2l+1}, \ref{parab2a+22a=k-5k=2l}, \ref{parab2a+22a=k-5k=2l+1}, \ref{parab2a+23a=-5k=3l}, and   \ref{parab2a+23a=-5k=3l+2}.}\label{Ta:fig:red552a+2}
\end{center}
\end{table}

\bigbreak

\bigbreak

\bigbreak

\bigbreak

\bigbreak

\begin{figure}[!ht]
    \psfrag{-4}{\huge $-4$}
    \psfrag{2l-7}{\huge $2l-7$}
    \psfrag{2l}{\huge $2l$}
    \psfrag{2l-4}{\huge $2l-4$}
    \psfrag{l-2}{\huge $l-2$}
    \psfrag{2l-5}{\huge $2l-5$}
    \psfrag{2l-3}{\huge $2l-3$}
    \psfrag{3l-5}{\huge $3l-5$}
    \psfrag{2l-2}{\huge $2l-2$}
    \psfrag{4l+1}{\huge $\mathbf{4l+1}$}
    \centerline{\scalebox{.35}{\includegraphics{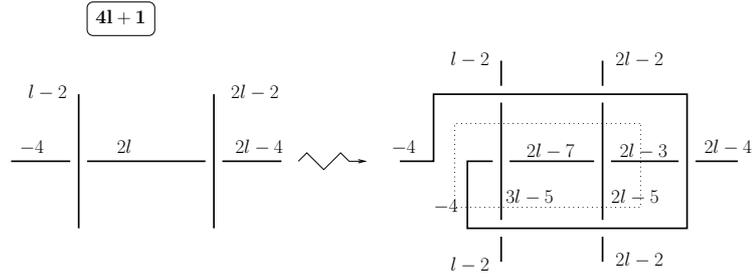}}}
    \caption{Removing color $k$ from an under-arc: the $b=2a+2$, $2a=k-4$; $k=2l$ $a=l-2$ mod $4l+1$  sub-instance. $a=l-2(<k=2l<4l-1<4l)$.}\label{fig:72aa2a+2}
\end{figure}

\paragraph{The $b=2a+2$, $2a=k-4$; $k=2l$, $a=l-2$ mod $4l+1$  sub-instance}\label{parab2a+22a=k-4k=2l} View Figure \ref{fig:72aa2a+2}.

\bigbreak

\bigbreak

\bigbreak

\bigbreak

\begin{table}[!ht]
\begin{center}
    \begin{tabular}{| c |}\hline
$l\geq 3 \qquad 2l-7 < 2l-5 < 2l-3 < k(=2l) \leq 3l-5 < 4l-1 < 4l $ \\ \hline
    \end{tabular}
\caption{Consider Figure \ref{fig:72aa2a+2}. $l\geq 3$. The equalities $2l-7=-1$,   and $2l=3l-5$ only give rise to further consideration with $l=3$. This corresponds to the elimination of color $2l-7=-1$, for $l=3$, in the boxed region in Figure \ref{fig:72aa2a+2} and is dealt with in Figure \ref{fig:72aaa2a+2}.}\label{Ta:fig:red72aa2a+2}
\end{center}
\end{table}

\bigbreak

\bigbreak

\begin{figure}[!ht]
    \psfrag{4}{\huge $4$}
    \psfrag{9}{\huge $9$}
    \psfrag{-1}{\huge $-1$}
    \psfrag{1}{\huge $1$}
    \psfrag{7}{\huge $7$}
    \psfrag{5}{\huge $5$}
    \psfrag{3}{\huge $3$}
    \psfrag{13}{\huge $\mathbf{13}$}
    \centerline{\scalebox{.35}{\includegraphics{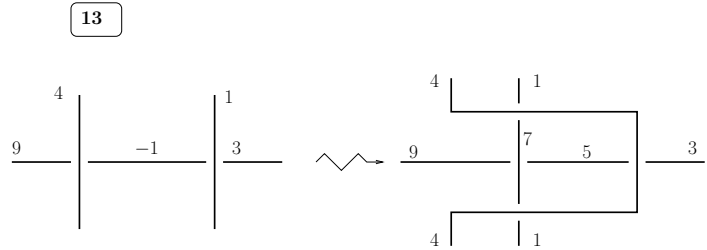}}}
    \caption{Removing color $k$ from an under-arc: the $b=2a+2$, $2a=k-4$, $k=2l$ sub-instance: fixing the boxed region in Figure \ref{fig:72aa2a+2} for $l=3$.}\label{fig:72aaa2a+2}
\end{figure}

\bigbreak

\begin{figure}[!ht]
    \psfrag{3l}{\huge $3l$}
    \psfrag{2l-1}{\huge $2l-1$}
    \psfrag{2l-4}{\huge $2l-4$}
    \psfrag{2l+1}{\huge $2l+1$}
    \psfrag{2l-3}{\huge $2l-3$}
    \psfrag{2l-2}{\huge $2l-2$}
    \psfrag{2l-6}{\huge $2l-6$}
    \psfrag{l-5}{\huge $l-5$}
    \psfrag{-4}{\huge $-4$}
    \psfrag{2l}{\huge $2l$}
    \psfrag{4l+3}{\huge $\mathbf{4l+3}$}
    \centerline{\scalebox{.35}{\includegraphics{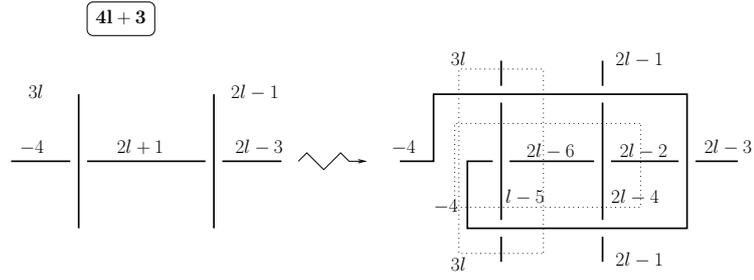}}}
    \caption{Removing color $k$ from an under-arc: the $b=2a+2$, $2a=k-4$; $k=2l+1$ $a=3l$ mod $4l+3$ sub-instance. $a=3l(<4l+1<4l+2))$. The vertically boxed region is fixed for $l=4$ in Figure \ref{fig:73a2a+219} and the horizontally boxed region is fixed for $l=2$ in Figure \ref{fig:73a2a+211} }\label{fig:73a2a+2}
\end{figure}

\paragraph{The $b=2a+2$, $2a=k-4$; $k=2l+1$, $a=3l$ mod $4l+3$ sub-instance.}\label{parab2a+22a=k-4k=2l+1} View Figure \ref{fig:73a2a+2}.

\begin{table}[!ht]
\begin{center}
   \scalebox{.95}{ \begin{tabular}{| c |}\hline
$l\geq 2 \qquad l-5\leq 2l-6 < 2l-4 < 2l-2 < k(=2l+1) < 4l+1 < 4l+2   $ \\ \hline
    \end{tabular}}
\caption{Consider Figure \ref{fig:73a2a+2}. $l\geq 2$. The equalities $l-5=2l-6$, $l-5=-1$, $l-5=-2$, $2l-6=-1$  and $2l-6=-2$ only give rise to further considerations for $l=2$ or $l=4$. This corresponds to the elimination of color $2l-6=-2$, for $l=2$, in the vertically boxed region in Figure \ref{fig:72aa2a+2} and to the horizontally boxed region in the same Figure. These situations are dealt with in Figures \ref{fig:73a2a+219} and \ref{fig:73a2a+211}.}\label{Ta:fig:red73a2a+2}
\end{center}
\end{table}

\bigbreak

\begin{figure}[!ht]
    \psfrag{-3}{\huge $-3$}
    \psfrag{-2}{\huge $-2$}
    \psfrag{2}{\huge $2$}
    \psfrag{0}{\huge $0$}
    \psfrag{3}{\huge $3$}
    \psfrag{4}{\huge $4$}
    \psfrag{-4}{\huge $-4$}
    \psfrag{11}{\huge $\mathbf{11}$}
    \centerline{\scalebox{.35}{\includegraphics{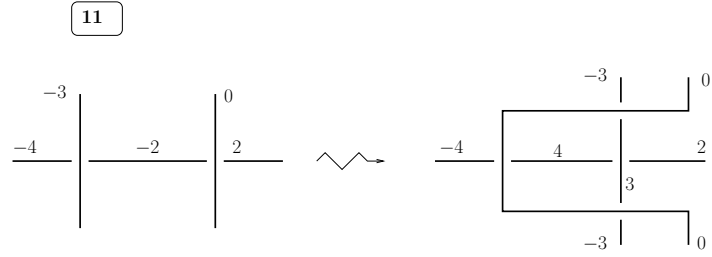}}}
    \caption{Removing color $k$ from an under-arc: the $b=2a+2$, $2a=k-4$; $k=2l+1$ $a=3l$ mod $4l+3$ sub-instance. Fixing  the horizontally boxed region for $l=2$ in Figure \ref{fig:73a2a+2} }\label{fig:73a2a+211}
\end{figure}

\bigbreak

\bigbreak

\bigbreak

\bigbreak

\begin{figure}[!ht]
    \psfrag{12}{\huge $12$}
    \psfrag{-1}{\huge $-1$}
    \psfrag{2}{\huge $2$}
    \psfrag{11}{\huge $11$}
    \psfrag{8}{\huge $8$}
    \psfrag{5}{\huge $5$}
    \psfrag{-4}{\huge $-4$}
    \psfrag{19}{\huge $\mathbf{19}$}
    \centerline{\scalebox{.35}{\includegraphics{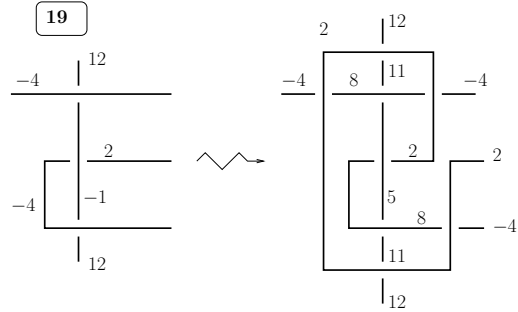}}}
    \caption{Removing color $k$ from an under-arc: the $b=2a+2$, $2a=k-4$; $k=2l+1$ $a=3l$ mod $4l+3$ sub-instance. Fixing  the vertically boxed region for $l=4$ in Figure \ref{fig:73a2a+2} }\label{fig:73a2a+219}
\end{figure}

\bigbreak

\begin{figure}[!ht]
    \psfrag{-4}{\huge $-5$}
    \psfrag{2l-7}{\huge $2l-9$}
    \psfrag{2l}{\huge $2l$}
    \psfrag{2l-4}{\huge $2l-6$}
    \psfrag{l-2}{\huge $3l-2$}
    \psfrag{2l-5}{\huge $2l-6$}
    \psfrag{2l-3}{\huge $2l-3$}
    \psfrag{3l-5}{\huge $l-7$}
    \psfrag{2l-2}{\huge $2l-3$}
    \psfrag{4l+1}{\huge $\mathbf{4l+1}$}
    \centerline{\scalebox{.35}{\includegraphics{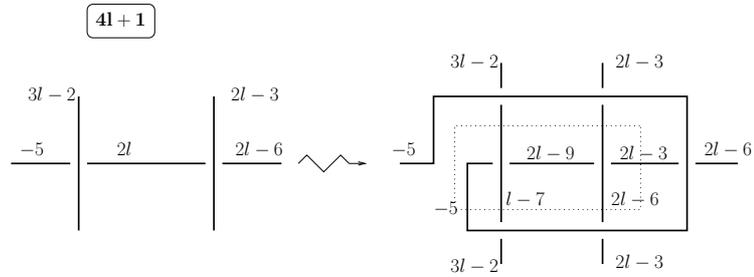}}}
    \caption{Removing color $k$ from an under-arc: the $b=2a+2$, $2a=k-5$; $k=2l$ $a=3l-2$ mod $4l+1$ sub-instance. $a=3l-2(<4l-1<4l)$.}\label{fig:72aa2a+2k-5}
\end{figure}

\paragraph{The $b=2a+2$, $2a=k-5$; $k=2l$, $a=3l-2$ mod $4l+1$ sub-instance}\label{parab2a+22a=k-5k=2l} View Figure \ref{fig:72aa2a+2k-5}.

\bigbreak

\bigbreak

\bigbreak

\begin{table}[!ht]
\begin{center}
    \scalebox{.9}{\begin{tabular}{| c  |}\hline
$l\geq 2 \qquad l-7 \leq 2l-9 < k(=2l) < 4l-1 < 4l $ \\ \hline
    \end{tabular}}
\caption{Consider Figure \ref{fig:72aa2a+2k-5}. $l\geq 2$. The equalities $l-7=2l-9$, $l-7=-1$, $l-7= -2$, $2l-9=-1$, and $2l-9=-2$ only give rise to further considerations for $l=4$. This corresponds to the elimination of color $2l-9=-1$, for $l=4$, in the boxed region in Figure \ref{fig:72aa2a+2k-5}. This situation
is dealt with in Figure \ref{fig:72aaa2a+2k-5}.}\label{Ta:fig:red72aa2a+2k-5}
\end{center}
\end{table}

\bigbreak

\begin{figure}[!ht]
    \psfrag{-3}{\huge $-3$}
    \psfrag{-2}{\huge $-1$}
    \psfrag{2}{\huge $5$}
    \psfrag{0}{\huge $2$}
    \psfrag{3}{\huge $7$}
    \psfrag{4}{\huge $9$}
    \psfrag{-4}{\huge $-5$}
    \psfrag{11}{\huge $\mathbf{17}$}
    \centerline{\scalebox{.35}{\includegraphics{eps42.eps}}}
     \caption{Removing color $k$ from an under-arc: the $b=2a+2$, $2a=k-5$; $k=2l$, $a=3l-2$ mod $4l+1$ sub-instance: fixing the boxed region in Figure \ref{fig:72aa2a+2k-5} for $l=4$.}\label{fig:72aaa2a+2k-5}
\end{figure}

\begin{figure}[!ht]
    \psfrag{3l}{\huge $l-2$}
    \psfrag{2l-1}{\huge $2l-2$}
    \psfrag{2l-4}{\huge $2l-5$}
    \psfrag{2l+1}{\huge $2l+1$}
    \psfrag{2l-3}{\huge $2l-5$}
    \psfrag{2l-2}{\huge $2l-2$}
    \psfrag{2l-6}{\huge $2l-8$}
    \psfrag{l-5}{\huge $3l-5$}
    \psfrag{-4}{\huge $-5$}
    \psfrag{2l}{\huge $2l$}
    \psfrag{4l+3}{\huge $\mathbf{4l+3}$}
    \centerline{\scalebox{.35}{\includegraphics{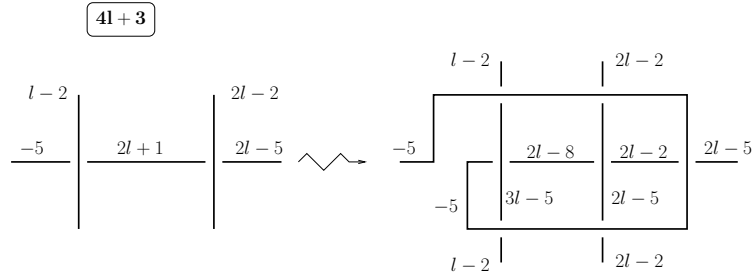}}}
    \caption{Removing color $k$ from an under-arc: the $b=2a+2$, $2a=k-4$; $k=2l+1$ $a=l-2$ mod $4l+3$ sub-instance. $a=l-2(<k=2l+1<4l+1<4l+2))$.}\label{fig:73a2a+2k-5}
\end{figure}

\paragraph{The $b=2a+2$, $2a=k-5$; $k=2l+1$, $a=l-2$ mod $4l+3$ sub-instance.}\label{parab2a+22a=k-5k=2l+1} View Figure \ref{fig:73a2a+2k-5}.

\bigbreak

\bigbreak

\begin{table}[!ht]
\begin{center}
   \scalebox{.95}{ \begin{tabular}{| c |}\hline
$ 2l-8 < k(=2l+1) \leq 3l-5 < 4l+1 < 4l+2 $\\ \hline
    \end{tabular}}
\caption{Consider Figure \ref{fig:73a2a+2k-5}.  $l\geq 2$. The equalities $2l+1=3l-5$,  $2l-8=-1$, and $2l-8=-2$ do not give rise to further considerations.}\label{Ta:fig:red73a2a+2k-5}
\end{center}
\end{table}

\bigbreak

\begin{figure}[!ht]
    \psfrag{3l}{\huge $3l$}
    \psfrag{5l-2}{\huge $5l-2$}
    \psfrag{4l-1}{\huge $4l-1$}
    \psfrag{2l-1}{\huge $2l-1$}
    \psfrag{l-2}{\huge $l-2$}
    \psfrag{3l-3}{\huge $3l-3$}
    \psfrag{6l-3}{\huge $6l-3$}
    \psfrag{l-1}{\huge $l-1$}
    \psfrag{l-5}{\huge $l-5$}
    \psfrag{2l-4}{\huge $2l-4$}
    \psfrag{6l+1}{\huge $\mathbf{6l+1}$}
    \centerline{\scalebox{.37}{\includegraphics{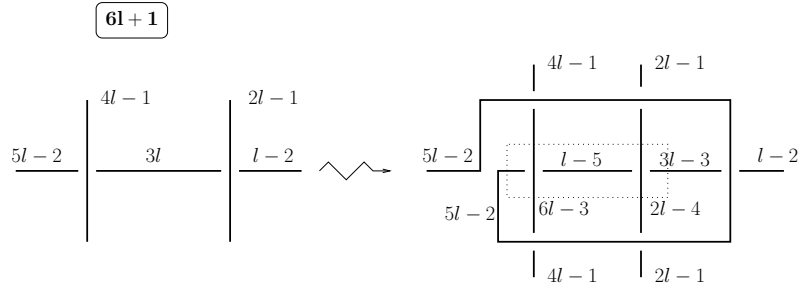}}}
    \caption{Removing color $k$ from an under-arc: the $b=2a+2$, $3a=-5$; $k=3l$ $a=4l-1$ mod $6l+1$ sub-instance. $a=4l-1(<6l-1<6l)$. The boxed region will be fixed for $l=3$ in Figure \ref{fig:red703a-5l=3}.}\label{fig:703a-5}
\end{figure}

\paragraph{The $b=2a+2$, $3a=-5$; $k=3l$ $a=4l-1$ mod $6l+1$ sub-instance}\label{parab2a+23a=-5k=3l} View Figure \ref{fig:703a-5}.

\bigbreak

\bigbreak

\begin{table}[!ht]
\begin{center}
    \begin{tabular}{| c |}\hline
$l\geq 2\qquad  l-5<2l-4<3l-3<k(=3l)\leq 6l-3  < 6l-1 < 6l $\\ \hline
    \end{tabular}
\caption{Consider Figure \ref{fig:703a-5}. $l\geq 2$. The equalities $6l-5=3l$,  $l-5=-1$, and $l-5=-2$ only give rise to further consideration with $l=3$. This corresponds to the the elimination of color $l-5=-2$, for $l=3$, in the boxed region in Figure \ref{fig:703a-5} and is dealt with in Figure \ref{fig:red703a-5l=3}.}\label{Ta:fig:red703a-5}
\end{center}
\end{table}

\bigbreak

\begin{figure}[!ht]
    \psfrag{13}{\huge $13$}
    \psfrag{15}{\huge $15$}
    \psfrag{-2}{\huge $-2$}
    \psfrag{2}{\huge $2$}
    \psfrag{6}{\huge $6$}
    \psfrag{8}{\huge $8$}
    \psfrag{10}{\huge $10$}
    \psfrag{19}{\huge $\mathbf{19}$}
    \centerline{\scalebox{.37}{\includegraphics{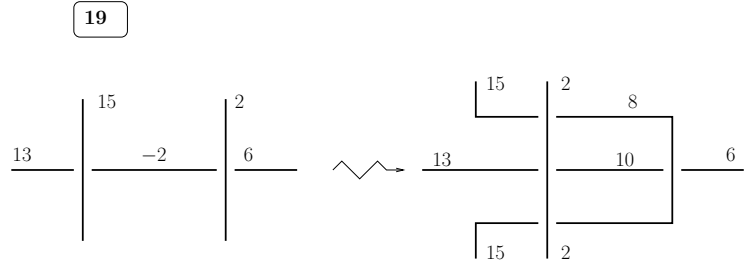}}}
    \caption{Removing color $k$ from an under-arc: the $b=2a+2$, $3a=-5$; $k=3l$ $a=4l-1$ mod $6l+1$ sub-instance. $a=4l-1(<6l-1<6l)$. Fixing the boxed region in Figure \ref{fig:703a-5} for $l=3$.}\label{fig:red703a-5l=3}
\end{figure}

\bigbreak

\begin{figure}[!ht]
    \psfrag{2l}{\huge $2l$}
    \psfrag{l-2}{\huge $l-2$}
    \psfrag{4l+2}{\huge $4l+2$}
    \psfrag{3l+2}{\huge $3l+2$}
    \psfrag{5l+2}{\huge $5l+2$}
    \psfrag{2l-1}{\huge $2l-1$}
    \psfrag{3l+1}{\huge $3l+1$}
    \psfrag{6l+2}{\huge $6l+2$}
    \psfrag{l-3}{\huge $l-3$}
    \psfrag{-3}{\huge $-3$}
    \psfrag{6l+5}{\huge $\mathbf{6l+5}$}
    \centerline{\scalebox{.35}{\includegraphics{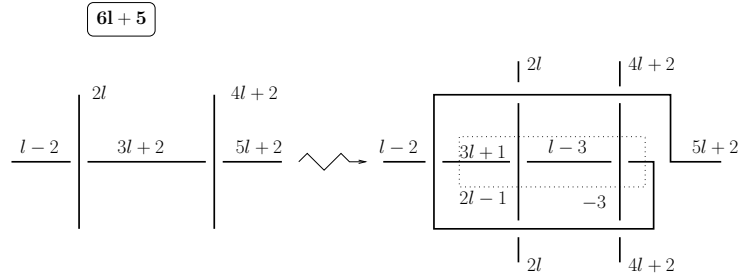}}}
    \caption{Removing color $k$ from an under-arc: the $b=2a+2$, $3a=-5$; $k=3l+2$, $a=2l$ mod $6l+5$ sub-instance. $a=2l(<k=3l+2<6l+3<6l+4)$. The boxed region will be fixed for $l=2$ in Figure \ref{fig:red713a-5l=3}.}\label{fig:red713a-5}
\end{figure}

\paragraph{The $b=2a+2$, $3a=-5$; $k=3l+2$, $a=2l$ mod $6l+5$ sub-instance}\label{parab2a+23a=-5k=3l+2} View Figure \ref{fig:red713a-5}.

\bigbreak

\begin{table}[!ht]
\begin{center}
    \begin{tabular}{| c |}\hline
$l\geq 1\qquad  l-3 < 2l-1 < 3l+1 < k(=3l+2) < 6l+3 < 6l+4$\\ \hline
    \end{tabular}
\caption{Consider Figure \ref{fig:red713a-5}. $l\geq 1$. The equalities $l-3=-1$, and $l-3=-2$ only give rise to further consideration with $l=2$. This corresponds to the boxed region in Figure \ref{fig:red713a-5} and is dealt with in Figure \ref{fig:red713a-5l=3}.}\label{Ta:fig:red713a-5}
\end{center}
\end{table}

\begin{figure}[!ht]
    \psfrag{-1}{\huge $-1$}
    \psfrag{3}{\huge $3$}
    \psfrag{7}{\huge $7$}
    \psfrag{9}{\huge $9$}
    \psfrag{11}{\huge $11$}
    \psfrag{12}{\huge $12$}
    \psfrag{14}{\huge $14$}
    \psfrag{17}{\huge $\mathbf{17}$}
    \centerline{\scalebox{.35}{\includegraphics{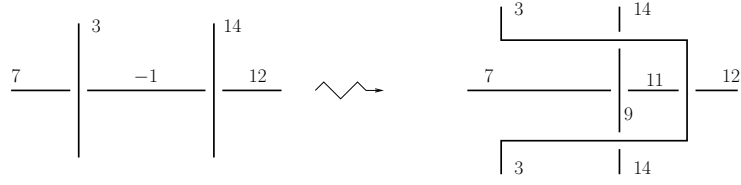}}}
    \caption{Removing color $k$ from an under-arc: the $b=2a+2$, $3a=-5$; $k=3l+2$, $a=2l$ mod $6l+5$ sub-instance. $a=2l(<k=3l+2<6l+3<6l+4)$. Fixing the boxed region for $l=2$ in Figure \ref{fig:red713a-5}.}\label{fig:red713a-5l=3}
\end{figure}

\subsubsection{The $b=2a-k$ instance}\label{subsubsectb2a-k} 
\begin{figure}[!ht]
    \psfrag{a}{\huge $a$}
    \psfrag{2a-k}{\huge $2a-k$}
    \psfrag{2a+2+k}{\huge $2a+2+k$}
    \psfrag{2a+1}{\huge $2a+1$}
    \psfrag{4a+k+2}{\huge $4a+k+2$}
    \psfrag{k}{\huge $k$}
    \psfrag{3a+1}{\huge $3a+1$}
    \centerline{\scalebox{.35}{\includegraphics{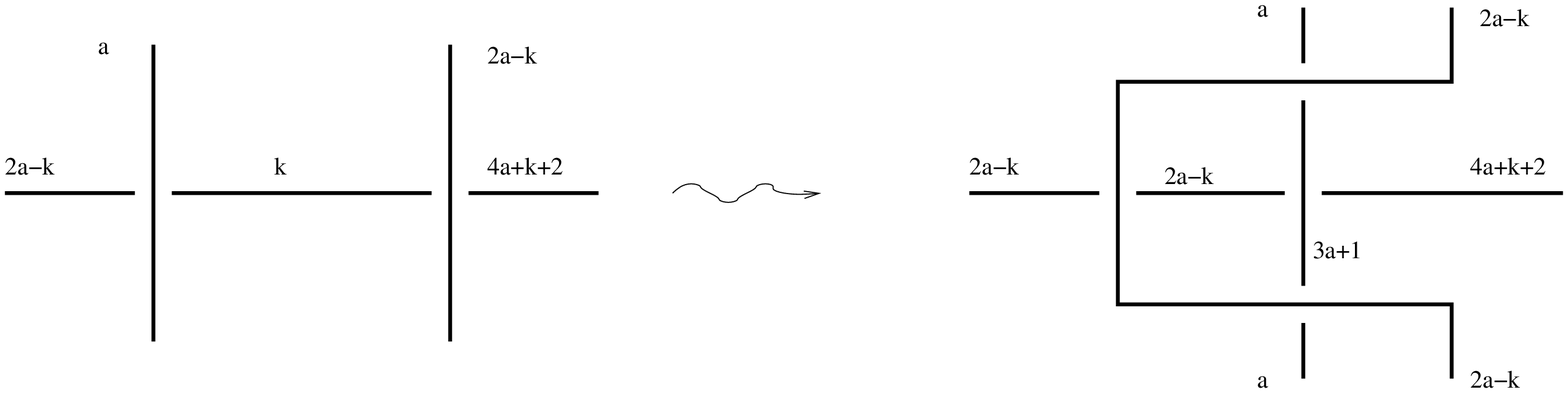}}}
    \caption{Removing color $k$ from an under-arc: the $b=2a-k$ instance}\label{fig:752a-k}
\end{figure}

View Figure \ref{fig:752a-k}.

\begin{table}[!ht]
\begin{center}
    \begin{tabular}{| c | c |  c | }\hline
3a+1=-1 & 3a+1=-2 & 3a+1=k(=3k+1)  \\ \hline
a=$_{6l+1}$4l or a=$_{6l+5}$2l+1  & a=-1 \, X &  a=k \, X    \\ \hline
    \end{tabular}
\caption{Equalities which should not occur in Figure \ref{fig:752a-k} (1st row) and their consequences (2nd row). $X$'s stand for  conflicts with assumptions thus not requiring further considerations. The other situations are dealt with in Figures \ref{fig:752a-kk=3l} and \ref{fig:752a-kk=3l+2}.}\label{Ta:fig:red752a-k}
\end{center}
\end{table}

\begin{figure}[!ht]
    \psfrag{3l}{\huge $3l$}
    \psfrag{4l}{\huge $4l$}
    \psfrag{5l}{\huge $5l$}
    \psfrag{l-1}{\huge $l-1$}
    \psfrag{3l-1}{\huge $3l-1$}
    \psfrag{4l-1}{\huge $4l-1$}
    \psfrag{5l-1}{\huge $5l-1$}
    \psfrag{6l+1}{\huge $\mathbf{6l+1}$}
    \centerline{\scalebox{.35}{\includegraphics{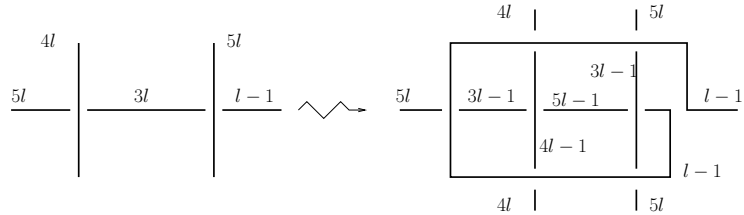}}}
    \caption{Removing color $k$ from an under-arc: the $b=2a-k$; $k=3l$, $a=4l$ mod $6l+1$ instance. $(k=3l<)a=4l(<6l-1<6l)$.}\label{fig:752a-kk=3l}
\end{figure}

\begin{table}[!ht]
\begin{center}
    \begin{tabular}{| c |  }\hline
$l\geq 2 \qquad 3l-1 < k(=3l) < 4l-1 < 5l-1 < 6l-1 < 6l $ \\ \hline
    \end{tabular}
\caption{Consider Figure \ref{fig:752a-kk=3l}.  $l\geq 2$. The equality $3l=4l-1$ does not give rise to further considerations. }\label{Ta:fig:red752a-kk=3l}
\end{center}
\end{table}

\begin{figure}[!ht]
    \psfrag{l}{\huge $l$}
    \psfrag{2l+1}{\huge $2l+1$}
    \psfrag{3l+2}{\huge $3l+2$}
    \psfrag{5l+3}{\huge $5l+3$}
    \psfrag{3l+1}{\huge $3l+1$}
    \psfrag{l-1}{\huge $l-1$}
    \psfrag{2l}{\huge $2l$}
    \psfrag{6l+5}{\huge $\mathbf{6l+5}$}
    \centerline{\scalebox{.35}{\includegraphics{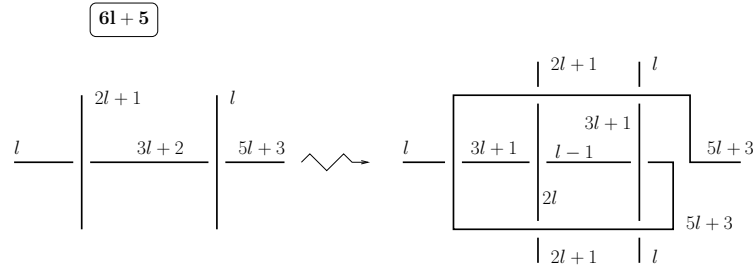}}}
    \caption{Removing color $k$ from an under-arc: the $b=2a-k$, $k=3l+2$ instance. $a=2l+1(<k=3l+2<6l+3<6l+4)$.}\label{fig:752a-kk=3l+2}
\end{figure}

\begin{table}[!ht]
\begin{center}
    \begin{tabular}{| c |  }\hline
$l\geq 1 \qquad    l-1 < 2l < 3l+1 < k(=3l+2) < 6l+3 < 6l+4   $ \\ \hline
    \end{tabular}
\caption{Consider Figure \ref{fig:752a-kk=3l+2}.  $l\geq 1$. There are no further considerations to be made here.}\label{Ta:fig:red752a-kk=3l+2}
\end{center}
\end{table}

\subsubsection{The $2a-2b+k=-1$ instance}\label{subsubsect2(b-a)=k+1}

\begin{figure}[!ht]
    \psfrag{0}{\huge $0$}
    \psfrag{a}{\huge $a$}
    \psfrag{2a-2l}{\huge $2a-2l$}
    \psfrag{2l}{\huge $2l$}
    \psfrag{a+3l+1}{\huge $a+3l+1$}
    \psfrag{a+2l+1}{\huge $a+2l+1$}
    \psfrag{2a+1}{\huge $2a+1$}
    \psfrag{4l+1}{\huge $\mathbf{4l+1}$}
    \centerline{\scalebox{.35}{\includegraphics{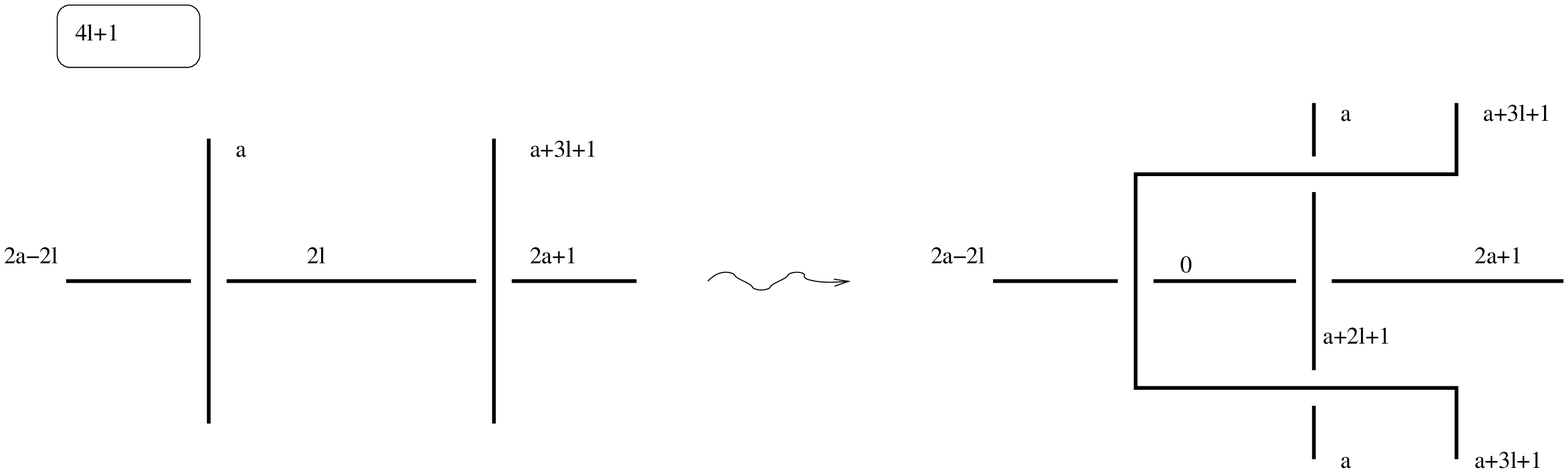}}}
    \caption{Removing color $k$ from an under-arc: the $2a-2b+k=-1$ instance with $k=2l$  and $b=a+3l+1$ mod $4l+1$.}\label{fig:redk2(b-a)=k+1k=2l}
\end{figure}

\paragraph{The $2a-2b+k=-1$ instance with $k=2l$  and $b=a+3l+1$ mod $4l+1$.}\label{para2(b-a)=k+1k=2l} View Figure \ref{fig:redk2(b-a)=k+1k=2l}.

\begin{table}[!ht]
\begin{center}
    \begin{tabular}{| c | c |  c | }\hline
a+2l+1=2l & a+2l+1=-1 & a+2l+1=-2   \\ \hline
a=-1 \, X & 2a+1=4l-1=-2 \, X &  a=2l-2    \\ \hline
    \end{tabular}
\caption{Equalities which should not occur in Figure \ref{fig:redk2(b-a)=k+1k=2l} (1st row) and their consequences (2nd row). $X$'s stand for  conflicts with assumptions thus not requiring further considerations.}\label{Ta:fig:redk2(b-a)=k+1k=2l}
\end{center}
\end{table}

\begin{figure}[!ht]
    \psfrag{2l-4}{\huge $2l-4$}
    \psfrag{2l-2}{\huge $2l-2$}
    \psfrag{2l}{\huge $2l$}
    \psfrag{l-2}{\huge $l-2$}
    \psfrag{-4}{\huge $-4$}
    \psfrag{2l-6}{\huge $2l-6$}
    \psfrag{2l-8}{\huge $2l-8$}
    \psfrag{3l-6}{\huge $3l-6$}
    \psfrag{4l-4}{\huge $4l-4$}
    \psfrag{4l+1}{\huge $\mathbf{4l+1}$}
    \centerline{\scalebox{.35}{\includegraphics{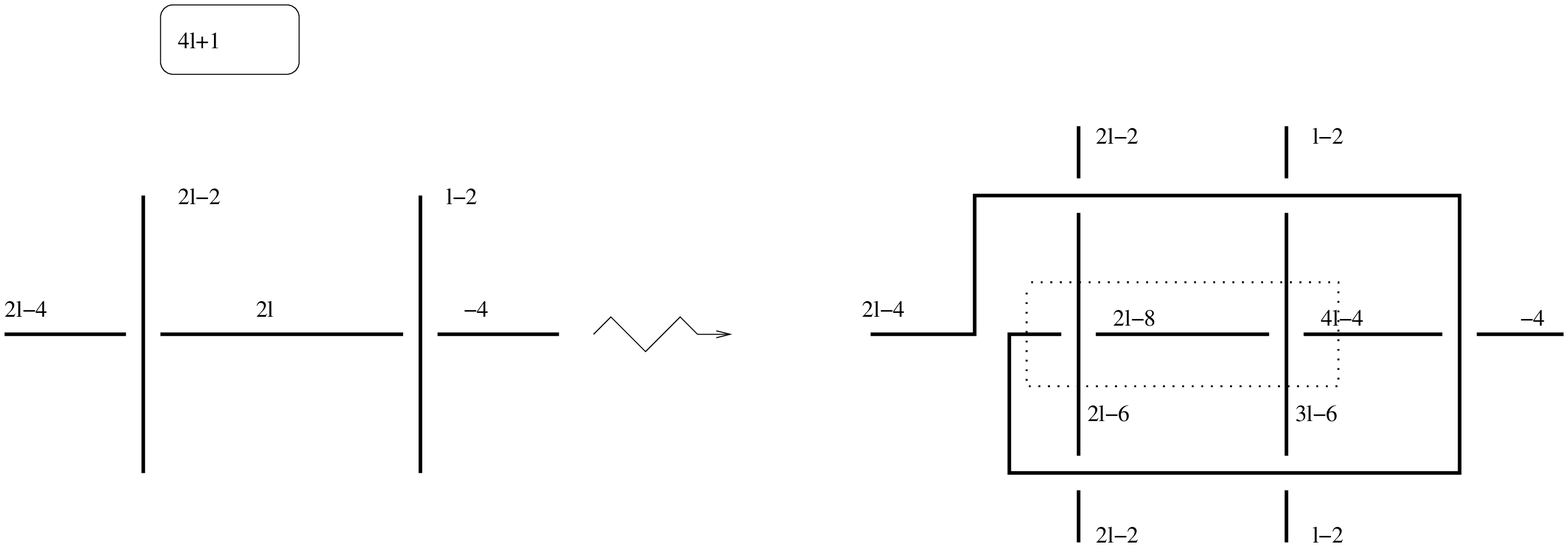}}}
    \caption{Removing color $k$ from an under-arc: the $2a-2b+k=-1$ instance with $k=2l$, $b=a+3l+1$ and $a=2l-2$, mod $4l+1$. The boxed region will be fixed for $l=3$ in Figure \ref{fig:752a-kk=3l+2}.}\label{fig:redk2(b-a)=k+1k=2la=2l-2}
\end{figure}

\begin{table}[!ht]
\begin{center}
    \scalebox{.74}{\begin{tabular}{| c |  }\hline
$l\geq 3 \qquad 2l-8 < 2l-6 < k(=2l)\leq 3l-6 < 4l-4 $\\ \hline
    \end{tabular}}
\caption{Consider Figure \ref{fig:redk2(b-a)=k+1k=2la=2l-2}. $l\geq 3$. The equalities $3l-6=2l$, $2l-8=-1$, and $2l-8=-2$ only give rise to further consideration with $l=3$. This corresponds to the boxed region in Figure \ref{fig:redk2(b-a)=k+1k=2la=2l-2} and is dealt with in Figure \ref{fig:redk2(b-a)=k+1k=2la=2l-2l=2}.}\label{Ta:fig:redk2(b-a)=k+1k=2la=2l-2}
\end{center}
\end{table}

\begin{figure}[!ht]
    \psfrag{1}{\huge $1$}
    \psfrag{2}{\huge $2$}
    \psfrag{0}{\huge $0$}
    \psfrag{3}{\huge $3$}
    \psfrag{4}{\huge $4$}
    \psfrag{5}{\huge $5$}
    \psfrag{6}{\huge $6$}
    \psfrag{8}{\huge $8$}
    \psfrag{9}{\huge $9$}
    \psfrag{-2}{\huge $-2$}
    \psfrag{10}{\huge $10$}
    \psfrag{13}{\huge $\mathbf{13}$}
    \centerline{\scalebox{.36}{\includegraphics{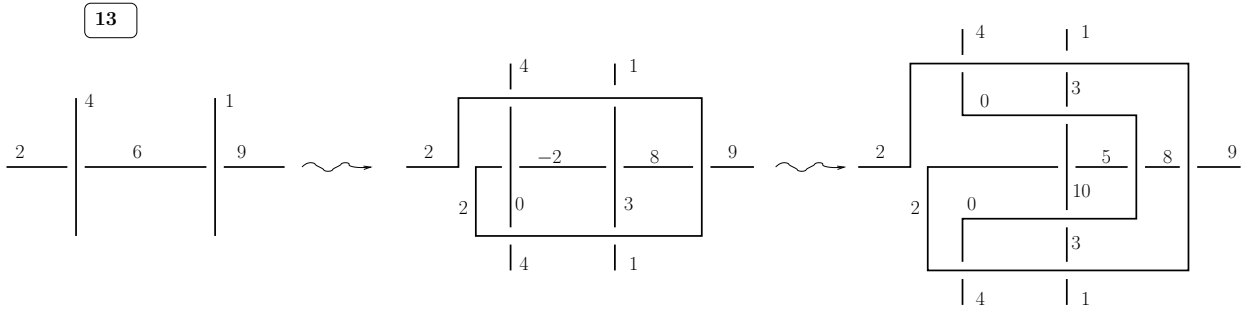}}}
    \caption{Removing color $k$ from an under-arc: the $2a-2b+k=-1$ instance with $k=2l$, $b=a+3l+1$ and $a=2l-2$, mod $4l+1$: fixing the boxed region for $l=3$ in Figure \ref{fig:redk2(b-a)=k+1k=2la=2l-2}.}\label{fig:redk2(b-a)=k+1k=2la=2l-2l=2}
\end{figure}

\begin{figure}[!ht]
    \psfrag{0}{\huge $0$}
    \psfrag{a}{\huge $a$}
    \psfrag{2a-2l}{\huge $2a-2l-1$}
    \psfrag{2l}{\huge $2l+1$}
    \psfrag{a+3l+1}{\huge $a+l+1$}
    \psfrag{a+2l+1}{\huge $a+2l+2$}
    \psfrag{2a+1}{\huge $2a+1$}
    \psfrag{4l+1}{\huge $\mathbf{4l+3}$}
    \centerline{\scalebox{.35}{\includegraphics{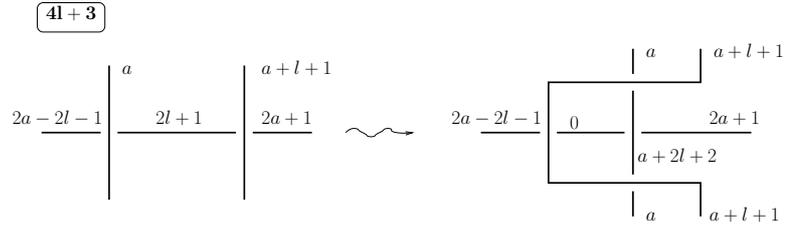}}}
    \caption{Removing color $k$ from an under-arc: the $2a-2b+k=-1$ instance with $k=2l+1$  and $b=a+l+1$ mod $4l+3$.}\label{fig:redk2(b-a)=k+1k=2l+1}
\end{figure}

\paragraph{The $2a-2b+k=-1$ instance with $k=2l+1$  and $b=a+l+1$ mod $4l+3$.}\label{para2(b-a)=k+1k=2l+1} View Figure \ref{fig:redk2(b-a)=k+1k=2l+1}.

\begin{table}[!ht]
\begin{center}
    \begin{tabular}{| c | c |  c | }\hline
a+2l+2=2l+1 & a+2l+2=-1 & a+2l+2=-2   \\ \hline
a=-1 \, X & 2a+1=4l+1=-2 \, X &  a=2l-1    \\ \hline
    \end{tabular}
\caption{Equalities which should not occur in Figure \ref{fig:redk2(b-a)=k+1k=2l+1} (1st row) and their consequences (2nd row). $X$'s stand for  conflicts with assumptions thus not requiring further considerations.}\label{Ta:fig:redk2(b-a)=k+1k=2l+1}
\end{center}
\end{table}

\begin{figure}[!ht]
    \psfrag{2l-1}{\huge $2l-1$}
    \psfrag{2l-3}{\huge $2l-3$}
    \psfrag{3l}{\huge $3l$}
    \psfrag{-4}{\huge $-4$}
    \psfrag{2l+1}{\huge $2l+1$}
    \psfrag{2l-4}{\huge $2l-4$}
    \psfrag{2l-2}{\huge $2l-2$}
    \psfrag{2l-6}{\huge $2l-6$}
    \psfrag{l-5}{\huge $l-5$}
    \psfrag{4l+3}{\huge $\mathbf{4l+3}$}
    \centerline{\scalebox{.35}{\includegraphics{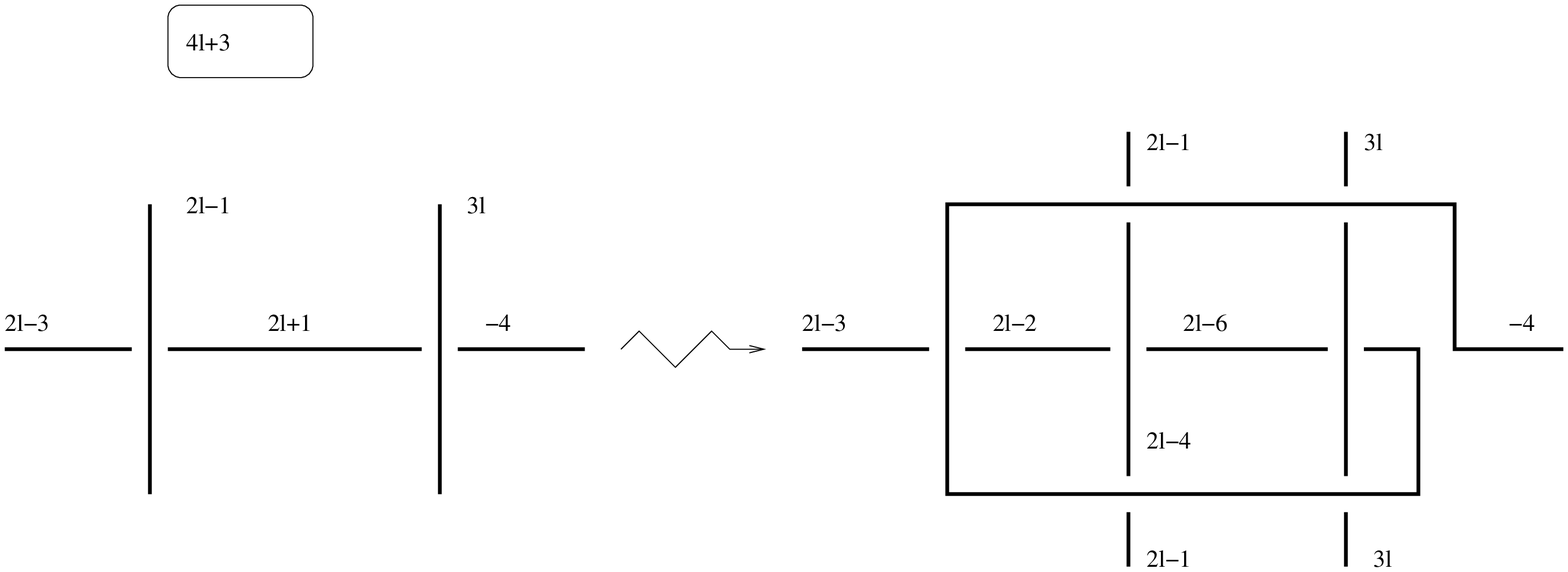}}}
    \caption{Removing color $k$ from an under-arc: the $2a-2b+k=-1$ instance with $k=2l+1$, $b=a+l+1$ and $a=2l-1$, mod $4l+3$. The boxed region will be fixed for $l=3$ in Figure \ref{fig:752a-kk=3l+2}.}\label{fig:redk2(b-a)=k+1k=2l+1a=2l-1}
\end{figure}

\begin{table}[h!]
\begin{center}
    \scalebox{.8}{\begin{tabular}{| c |  }\hline
$l\geq 2 \qquad l-5 \leq 2l-6 < 2l-4 < 2l-2 < k(=2l+1) < 4l+1 < 4l+3$\\ \hline
    \end{tabular}}
\caption{Consider Figure \ref{fig:redk2(b-a)=k+1k=2l+1a=2l-1}. $l\geq 2$. The equalities $l-5=2l-6$, $l-5=-1$, $l-5=-2$, $2l-6=-1$, and $2l-6=-2$ only give rise to further consideration for $l=2$ or for $l=4$. These situations are dealt with in Figures \ref{fig:752a-kk=3l+2l=2} and \ref{fig:752a-kk=3l+2l=4}.}\label{Ta:fig:redk2(b-a)=k+1k=2l+1a=2l-1}
\end{center}
\end{table}

\begin{figure}[!ht]
    \psfrag{1}{\huge $1$}
    \psfrag{2}{\huge $2$}
    \psfrag{0}{\huge $0$}
    \psfrag{3}{\huge $3$}
    \psfrag{4}{\huge $4$}
    \psfrag{5}{\huge $5$}
    \psfrag{6}{\huge $6$}
    \psfrag{8}{\huge $8$}
    \psfrag{9}{\huge $9$}
    \psfrag{-2}{\huge $-2$}
    \psfrag{7}{\huge $7$}
    \psfrag{11}{\huge $\mathbf{11}$}
    \centerline{\scalebox{.36}{\includegraphics{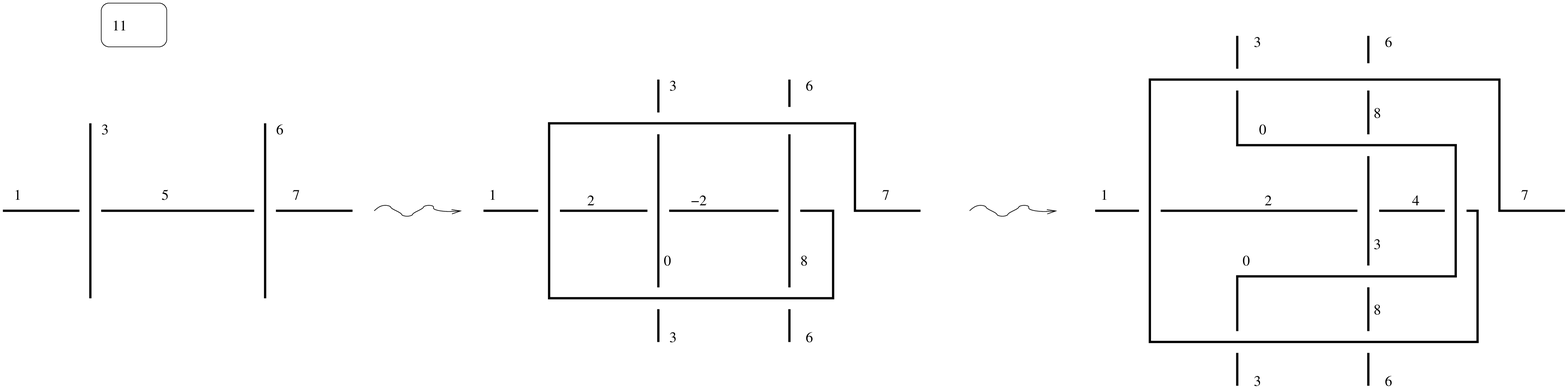}}}
    \caption{Removing color $k$ from an under-arc: the $2a-2b+k=-1$ instance with $k=2l+1$, $b=a+l+1$ and $a=2l-1$, mod $4l+3$: fixing the boxed region for $l=2$ in Figure \ref{fig:redk2(b-a)=k+1k=2l+1a=2l-1}.}\label{fig:752a-kk=3l+2l=2}
\end{figure}

\begin{figure}[!ht]
    \psfrag{1}{\huge $1$}
    \psfrag{2}{\huge $2$}
    \psfrag{0}{\huge $0$}
    \psfrag{3}{\huge $3$}
    \psfrag{4}{\huge $4$}
    \psfrag{5}{\huge $5$}
    \psfrag{6}{\huge $6$}
    \psfrag{12}{\huge $12$}
    \psfrag{9}{\huge $9$}
    \psfrag{-4}{\huge $-4$}
    \psfrag{7}{\huge $7$}
    \psfrag{19}{\huge $\mathbf{19}$}
    \centerline{\scalebox{.35}{\includegraphics{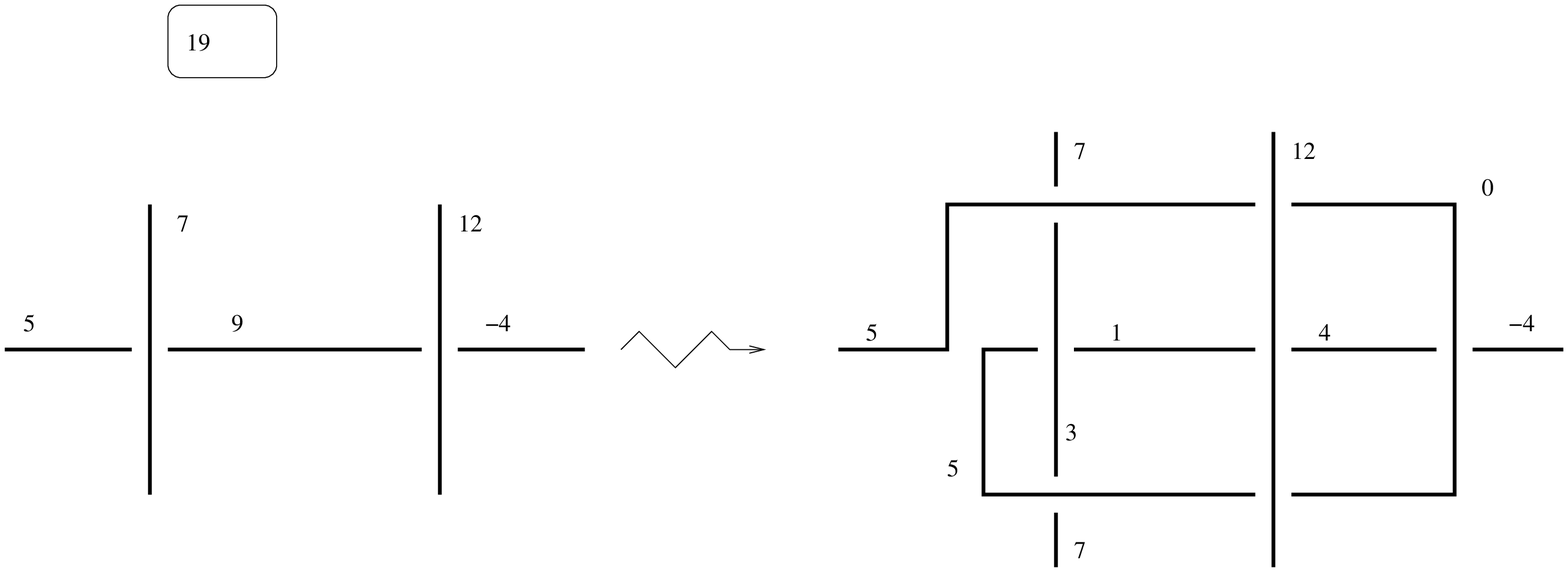}}}
    \caption{Removing color $k$ from an under-arc: the $2a-2b+k=-1$ instance with $k=2l+1$, $b=a+l+1$ and $a=2l-1$, mod $4l+3$: fixing the boxed region for $l=4$ in Figure \ref{fig:redk2(b-a)=k+1k=2l+1a=2l-1}.}\label{fig:752a-kk=3l+2l=4}
\end{figure}

\subsubsection{The $2a-2b+k=-2$ instance}\label{subsubsect2(b-a)=k+2}

\bigbreak

\begin{figure}[!ht]
    \psfrag{0}{\huge $1$}
    \psfrag{a}{\huge $a$}
    \psfrag{2a-2l}{\huge $2a-2l$}
    \psfrag{2l}{\huge $2l$}
    \psfrag{a+3l+1}{\huge $a+l+1$}
    \psfrag{a+2l+1}{\huge $a+2l+2$}
    \psfrag{2a+1}{\huge $2a+2$}
    \psfrag{4l+1}{\huge $\mathbf{4l+1}$}
    \centerline{\scalebox{.35}{\includegraphics{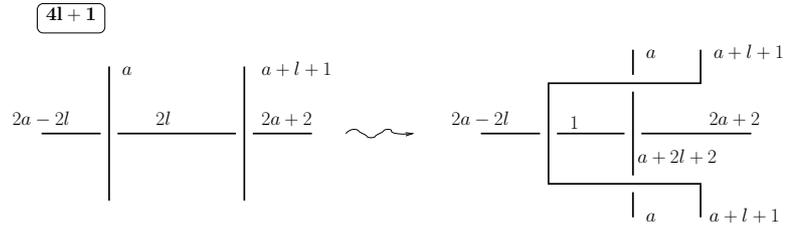}}}
    \caption{Removing color $k$ from an under-arc: the $2a-2b+k=-2$ instance with $k=2l$ and $b=a+l+1$ mod  $4l+1$.}\label{fig:redk2(b-a)=k+2k=2l}
\end{figure}

\paragraph{The $2a-2b+k=-2$ instance with $k=2l$ and $b=a+l+1$ mod  $4l+1$.}\label{para2(b-a)=k+12k=2l} View Figure \ref{fig:redk2(b-a)=k+2k=2l}.

\begin{table}[h!]
\begin{center}
    \begin{tabular}{| c | c |  c | }\hline
a+2l+2=2l & a+2l+2=-1 & a+2l+2=-2   \\ \hline
a=-2 \, X &  a=2l-2   & a=2l-3   \\ \hline
    \end{tabular}
\caption{Equalities which should not occur in Figure \ref{fig:redk2(b-a)=k+2k=2l} (1st row) and their consequences (2nd row). $X$'s stand for  conflicts with assumptions thus not requiring further considerations.}\label{Ta:fig:redk2(b-a)=k+2k=2l}
\end{center}
\end{table}

\bigbreak

\begin{figure}[!ht]
    \psfrag{2l-4}{\huge $2l-4$}
    \psfrag{2l-2}{\huge $2l-2$}
    \psfrag{2l}{\huge $2l$}
    \psfrag{l-2}{\huge $3l-1$}
    \psfrag{-4}{\huge $-3$}
    \psfrag{2l-6}{\huge $2l-6$}
    \psfrag{2l-8}{\huge $2l-8$}
    \psfrag{3l-6}{\huge $l-7$}
    \psfrag{4l-4}{\huge $4l-5$}
    \psfrag{4l+1}{\huge $\mathbf{4l+1}$}
    \centerline{\scalebox{.35}{\includegraphics{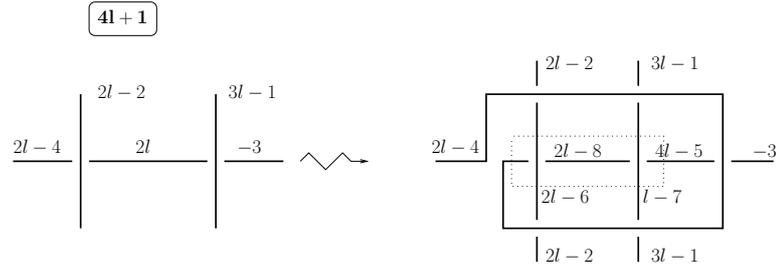}}}
    \caption{Removing color $k$ from an under-arc: the $2a-2b+k=-2$ instance with $k=2l$, $b=a+l+1$ and $a=2l-2$, mod  $4l+1$. The boxed region will be fixed for $l=3$ in Figure \ref{fig:redl=k2(b-a)=k+2k=2la=2l-2bistri}.}\label{fig:redk2(b-a)=k+2k=2la=2l-2tritri}
\end{figure}

\bigbreak

\begin{table}[h!]
\begin{center}
    \begin{tabular}{| c  |  }\hline
$l\geq 3 \qquad l-7 \leq 2l-8 < 2l-6 < k(=2l) < 4l-5 < 4l-1 < 4l$\\ \hline
    \end{tabular}
\caption{Consider Figure \ref{fig:redk2(b-a)=k+2k=2la=2l-2tritri}. $l\geq 3$. The equalities $l-7=-1$, $l-7=-2$, $2l-8=-1$, and $2l-8=-2$ only give rise to further consideration for $l=3$. This situation is dealt with in Figure \ref{fig:redl=k2(b-a)=k+2k=2la=2l-2bistri}.}\label{Ta:fig:redk2(b-a)=k+2k=2la=2l-2}
\end{center}
\end{table}

\bigbreak

\begin{figure}[!ht]
    \psfrag{1}{\huge $1$}
    \psfrag{2}{\huge $2$}
    \psfrag{0}{\huge $0$}
    \psfrag{3}{\huge $3$}
    \psfrag{4}{\huge $4$}
    \psfrag{5}{\huge $5$}
    \psfrag{6}{\huge $6$}
    \psfrag{8}{\huge $8$}
    \psfrag{9}{\huge $9$}
    \psfrag{-3}{\huge $-3$}
    \psfrag{-2}{\huge $-2$}
    \psfrag{7}{\huge $7$}
    \psfrag{13}{\huge $\mathbf{13}$}
    \centerline{\scalebox{.35}{\includegraphics{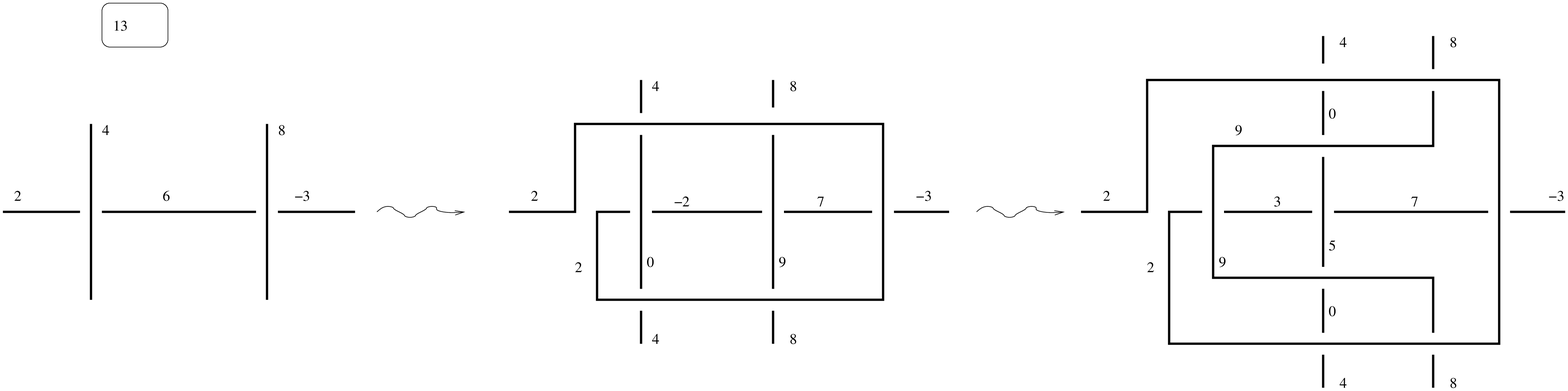}}}
    \caption{Removing color $k$ from an under-arc: the $2a-2b+k=-2$ instance with $k=2l$, $b=a+l+1$, and $a=2l-2$ mod  $4l+1$: fixing the boxed region for $l=3$ in Figure \ref{fig:redk2(b-a)=k+2k=2la=2l-2tritri}.}\label{fig:redl=k2(b-a)=k+2k=2la=2l-2bistri}
\end{figure}

\bigbreak

\bigbreak

\begin{figure}[!ht]
    \psfrag{2l-6}{\huge $2l-6$}
    \psfrag{2l-3}{\huge $2l-3$}
    \psfrag{3l-2}{\huge $3l-2$}
    \psfrag{2l}{\huge $2l$}
    \psfrag{l-7}{\huge $l-7$}
    \psfrag{2l-9}{\huge $2l-9$}
    \psfrag{4l-4}{\huge $4l-4$}
    \psfrag{4l+1}{\huge $\mathbf{4l+1}$}
    \centerline{\scalebox{.35}{\includegraphics{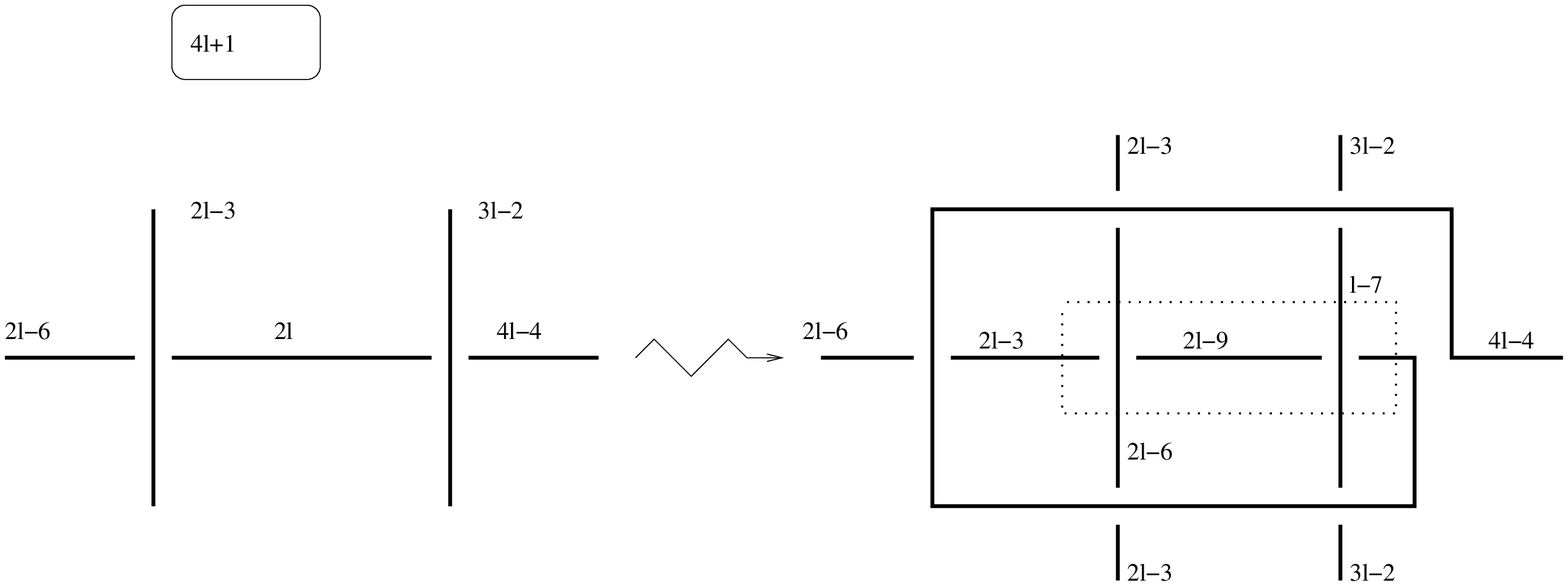}}}
    \caption{Removing color $k$ from an under-arc: the $2a-2b+k=-2$ instance with $k=2l$, $b=a+l+1$, and $a=2l-3$ mod  $4l+1$. The boxed region will be fixed for $l=4$ in Figure \ref{fig:752a-kk+2k=3l+2}.}\label{fig:redk2(b-a)=k+2k=2la=2l-2}
\end{figure}
\begin{table}[h!]
\begin{center}
    \begin{tabular}{ | c |  }\hline
$l\geq 3  \qquad l-7 \leq 2l-9 < k(=2l) < 4l-1 < 4l$\\ \hline
    \end{tabular}
\caption{Consider Figure \ref{fig:redk2(b-a)=k+2k=2la=2l-2}. $l\geq 3$. The equalities $l-7=-1$, $l-7=-2$, $2l-9=-1$, and $2l-9=-2$  give rise to further consideration for $l=4$. This situation is dealt with in Figure \ref{fig:752a-kk+2k=3l+2}.}\label{Ta:fig:redk2(b-a)=k+2k=2la=2l-2bis}
\end{center}
\end{table}

\bigbreak

\begin{figure}[!ht]
    \psfrag{1}{\huge $1$}
    \psfrag{2}{\huge $2$}
    \psfrag{0}{\huge $0$}
    \psfrag{10}{\huge $10$}
    \psfrag{4}{\huge $4$}
    \psfrag{5}{\huge $5$}
    \psfrag{14}{\huge $14$}
    \psfrag{8}{\huge $8$}
    \psfrag{9}{\huge $9$}
    \psfrag{-1}{\huge $-1$}
    \psfrag{12}{\huge $12$}
    \psfrag{7}{\huge $7$}
    \psfrag{17}{\huge $\mathbf{17}$}
    \centerline{\scalebox{.36}{\includegraphics{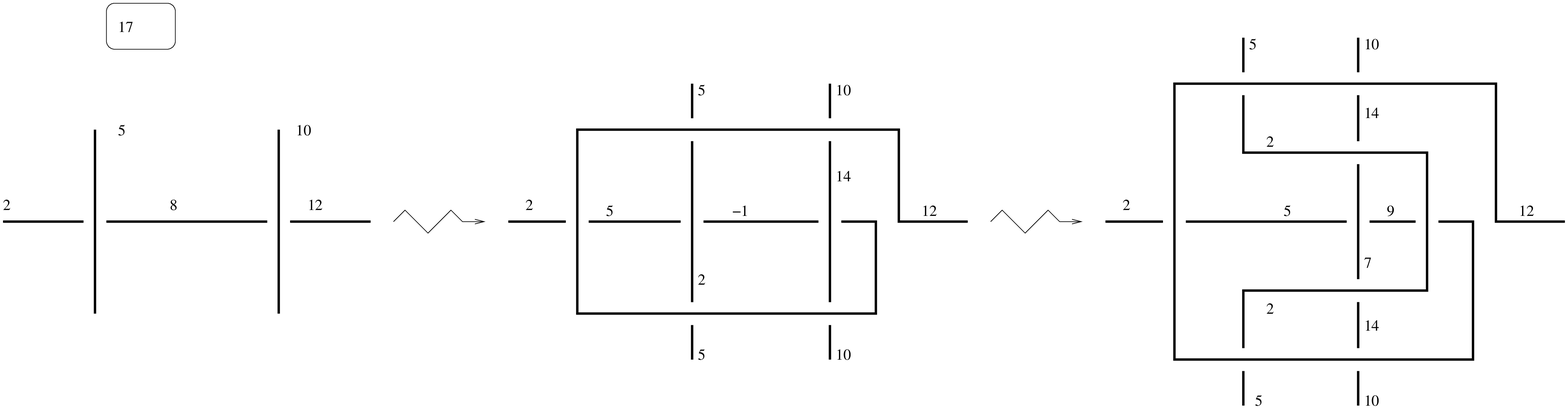}}}
    \caption{Removing color $k$ from an under-arc: the $2a-2b+k=-2$ instance with $k=2l$, $b=a+l+1$, and $a=2l-3$ mod  $4l+1$: fixing the boxed region for $l=4$ in Figure \ref{fig:redk2(b-a)=k+2k=2la=2l-2}.}\label{fig:752a-kk+2k=3l+2}
\end{figure}

\paragraph{The $2a-2b+k=-2$ instance with $k=2l+1$ and $b=a+3l+3$ mod  $4l+3$.}\label{para2(b-a)=k+2k=2l+1} View Figure \ref{fig:redk2(b-a)=k+2k=2l}.

\begin{figure}[!ht]
    \psfrag{1}{\huge $1$}
    \psfrag{a}{\huge $a$}
    \psfrag{2a+2l+2}{\huge $2a+2l+2$}
    \psfrag{2l+1}{\huge $2l+1$}
    \psfrag{a+3l+3}{\huge $a+3l+3$}
    \psfrag{a+2l+3}{\huge $a+2l+3$}
    \psfrag{2a+2}{\huge $2a+2$}
    \psfrag{4l+3}{\huge $\mathbf{4l+3}$}
    \centerline{\scalebox{.35}{\includegraphics{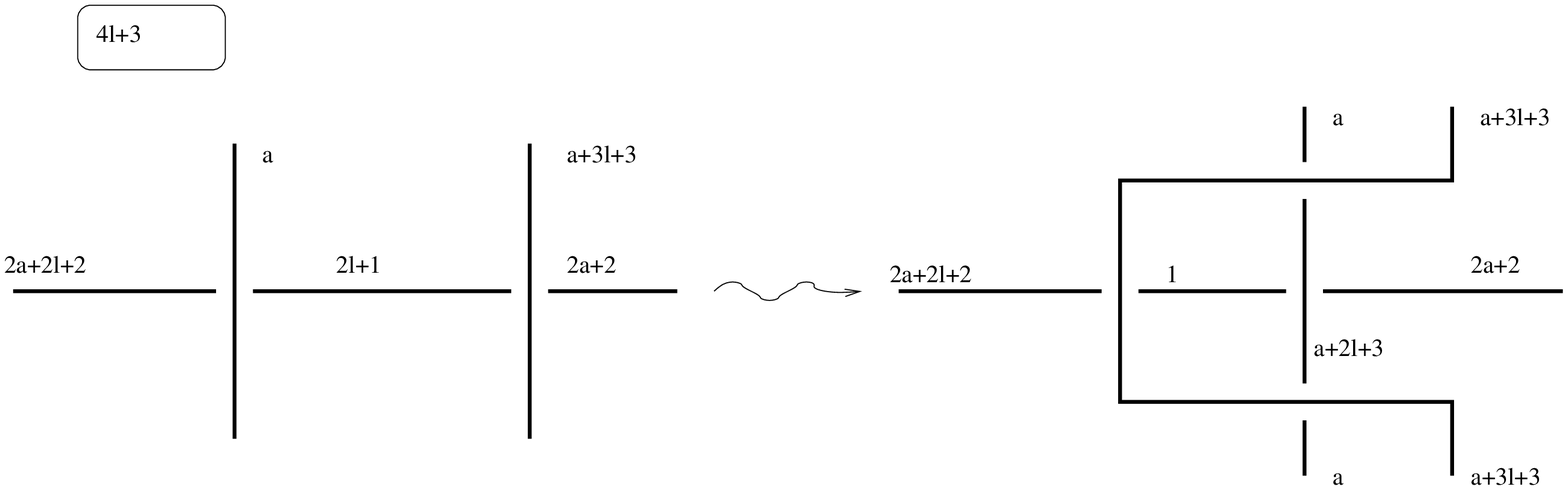}}}
    \caption{Removing color $k$ from an under-arc: the $2a-2b+k=-2$ instance with $k=2l+1$ and $b=a+3l+3$ mod  $4l+3$.}\label{fig:redk2(b-a)=k+2k=2l+}
\end{figure}

\begin{table}[h!]
\begin{center}
    \begin{tabular}{| c | c |  c | }\hline
a+2l+3=2l+1 & a+2l+3=-1 & a+2l+3=-2   \\ \hline
a=-2 \, X &  a=2l-1   & a=2l-2   \\ \hline
    \end{tabular}
\caption{Equalities which should not occur in Figure \ref{fig:redk2(b-a)=k+2k=2l+} (1st row) and their consequences (2nd row). $X$'s stand for  conflicts with assumptions thus not requiring further considerations.}\label{Ta:fig:redk2(b-a)=k+2k=2l+}
\end{center}
\end{table}

\begin{figure}[!ht]
    \psfrag{2l-3}{\huge $2l-3$}
    \psfrag{2l-1}{\huge $2l-1$}
    \psfrag{2l+1}{\huge $2l+1$}
    \psfrag{l-1}{\huge $l-1$}
    \psfrag{-3}{\huge $-3$}
    \psfrag{2l-5}{\huge $2l-5$}
    \psfrag{2l-7}{\huge $2l-7$}
    \psfrag{3l-5}{\huge $3l-5$}
    \psfrag{4l-3}{\huge $4l-3$}
    \psfrag{4l+3}{\huge $\mathbf{4l+3}$}
    \centerline{\scalebox{.35}{\includegraphics{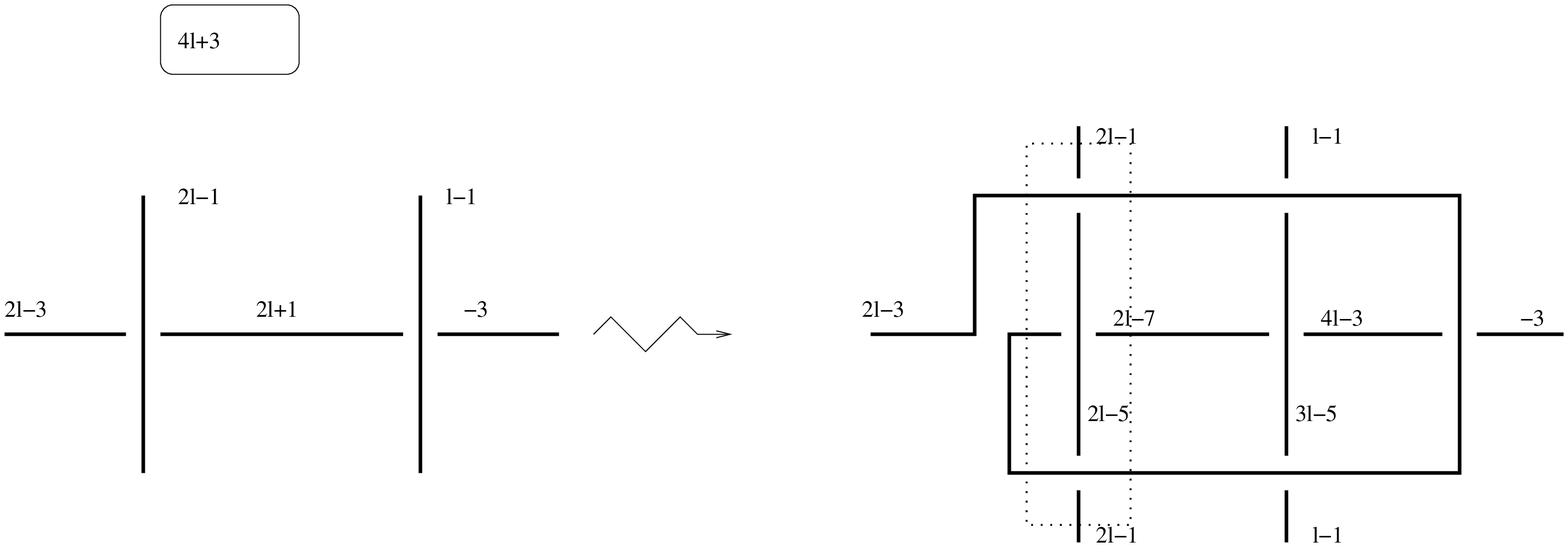}}}
    \caption{Removing color $k$ from an under-arc: the $2a-2b+k=-2$ instance with $k=2l+1$, $b=a+3l+3$, and $a=2l-1$, mod  $4l+3$ . The boxed region will be fixed for $l=2$ in Figure \ref{fig:red2(b-a)=k+2k=2l+1}.}\label{fig:redk2(b-a)=k+2k=2l+1a=2l-1}
\end{figure}

\bigbreak

\bigbreak

\bigbreak

\bigbreak

\bigbreak

\bigbreak

\begin{table}[h!]
\begin{center}
    \begin{tabular}{| c | }\hline
$l\geq 2 \qquad 2l-7 < 2l-5 < k(=2l+1) \leq 3l-5 < 4l-3 < 4l+1 < 4l+2$\\ \hline
    \end{tabular}
\caption{Consider Figure \ref{fig:redk2(b-a)=k+2k=2l+1a=2l-1}. $l\geq 2$. The equalities $2l+1=3l-5$, $2l-7=-1$, $2l-7=-2$, $2l-5=-1$, and $2l-5=-2$  give rise to further consideration for $l=2$. This situation is dealt with in Figure \ref{fig:red2(b-a)=k+2k=2l+1}.}\label{Ta:fig:redk2(b-a)=k+2k=2l+1a=2l-1}
\end{center}
\end{table}

\bigbreak

\bigbreak

\bigbreak

\begin{figure}[!ht]
    \psfrag{1}{\huge $1$}
    \psfrag{2}{\huge $2$}
    \psfrag{0}{\huge $0$}
    \psfrag{3}{\huge $3$}
    \psfrag{4}{\huge $4$}
    \psfrag{5}{\huge $5$}
    \psfrag{6}{\huge $6$}
    \psfrag{8}{\huge $8$}
    \psfrag{9}{\huge $9$}
    \psfrag{-3}{\huge $-3$}
    \psfrag{-2}{\huge $-2$}
    \psfrag{7}{\huge $7$}
    \psfrag{11}{\huge $\mathbf{11}$}
    \centerline{\scalebox{.35}{\includegraphics{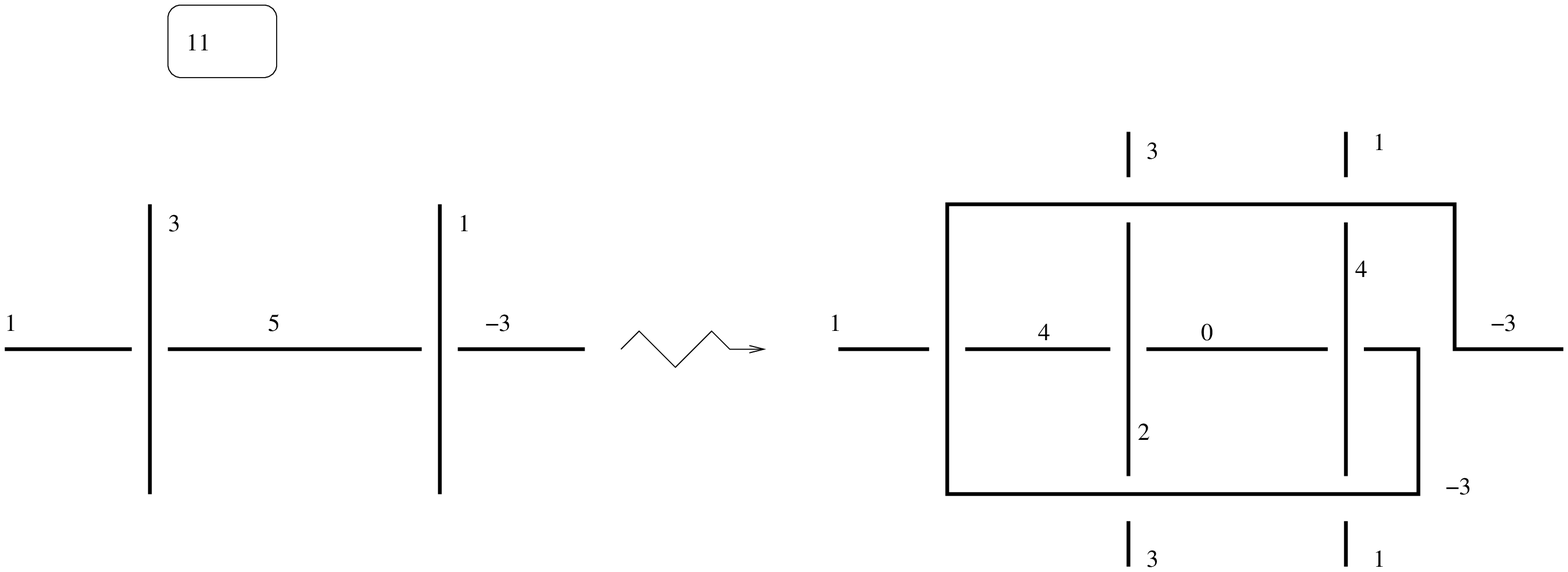}}}
    \caption{Removing color $k$ from an under-arc: the $2a-2b+k=-2$ instance with $k=2l+1$, $b=a+3l+3$, and $a=2l-1$, mod  $4l+3$ : fixing the boxed region for $l=2$ in Figure \ref{fig:redk2(b-a)=k+2k=2l+1a=2l-1}.}
\label{fig:red2(b-a)=k+2k=2l+1}
\end{figure}

\bigbreak

\bigbreak

\bigbreak

\bigbreak

\bigbreak

\begin{figure}[h]
    \psfrag{2l-5}{\huge $2l-5$}
    \psfrag{2l-2}{\huge $2l-2$}
    \psfrag{3l-5}{\huge $3l-5$}
    \psfrag{2l+1}{\huge $2l+1$}
    \psfrag{l-2}{\huge $l-2$}
    \psfrag{2l-8}{\huge $2l-8$}
    \psfrag{4l-2}{\huge $4l-2$}
    \psfrag{4l+3}{\huge $\mathbf{4l+3}$}
    \centerline{\scalebox{.35}{\includegraphics{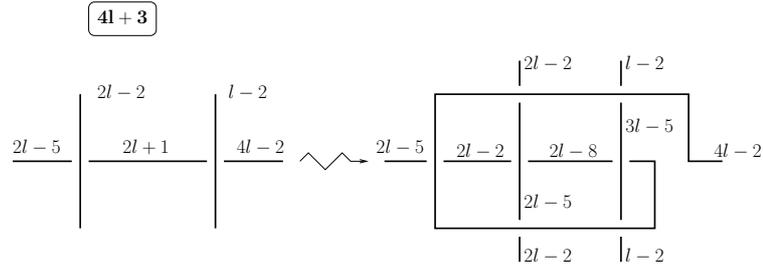}}}
    \caption{Removing color $k$ from an under-arc: the $2a-2b+k=-2$ instance with $k=2l+1$, $b=a+3l+3$, and $a=2l-2$, mod  $4l+3$.}\label{fig:redk2(b-a)=k+2k=2l+1a=2l-2}
\end{figure}

\bigbreak

\bigbreak

\bigbreak

\begin{table}[h!]
\begin{center}
    \begin{tabular}{ | c |  }\hline
 $l\geq 2 \qquad 2l-8 < k(=2l+1) \leq 3l-5 < 4l+1 < 4l+2 $\\ \hline
    \end{tabular}
\caption{Consider Figure \ref{fig:redk2(b-a)=k+2k=2l+1a=2l-2}. $l\geq 2$. The equalities $2l+1=3l-5$, $2l-8=-1$, and $2l-8=-2$  do not give rise to further considerations.}\label{Ta:fig:redk2(b-a)=k+2k=2l+1a=2l-2}
\end{center}
\end{table}

\bigbreak

\bigbreak

\bigbreak

\subsubsection{The $b=a$ instance}\label{subsubsectb=a} 
\begin{figure}[!ht]
    \psfrag{a}{\huge $a$}
    \psfrag{k}{\huge $k$}
    \psfrag{2a-k}{\huge $2a-k$}
    \psfrag{3a+1}{\huge $3a+1$}
    \psfrag{4a+2+k}{\huge $4a+k+2$}
    \centerline{\scalebox{.35}{\includegraphics{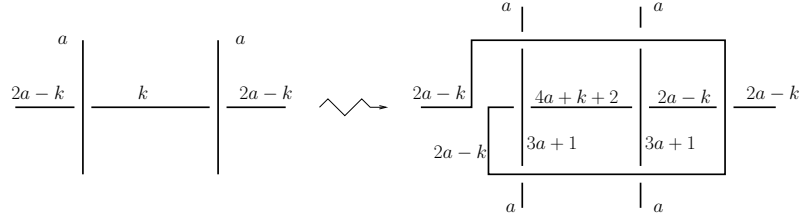}}}
    \caption{Removing color $k$ from an under-arc: recalling the setting for the $b=a$ instance.}\label{fig:red65}
\end{figure}

View Figure \ref{fig:red65}.

\bigbreak

\bigbreak

\bigbreak

\bigbreak

\bigbreak

\bigbreak

\bigbreak

\begin{table}[h!]
\begin{center}
   \scalebox{.6275}{ \begin{tabular}{ | c | c | c |  c | c | c |  }\hline
   3a+1=k &  3a+1=-1  &  3a+1=-2 & 4a+k+2=k  &  4a+k+2=-1  &  4a+k+2=-2\\ \hline
  a=k \, X &  a=$_{6l+1}$4l or  a=$_{6l+5}$2l+1  &  a=-1 \, X  &  a=k \, X &  a=$_{8l+1}$5l or a=$_{8l+3}$7l+2 or a=$_{8l+5}$l or a=$_{8l+7}$3l+2  &  a=$_{8l+1}$7l or a=$_{8l+3}$5l+1 or a=$_{8l+5}$3l+1 or a=$_{8l+7}$l\\ \hline
    \end{tabular}}
\caption{Equalities which should not occur in Figure \ref{fig:red65} (1st row) and their consequences (2nd row). $X$'s stand for  conflicts with assumptions thus not requiring further considerations. The other situations are dealt with in \ref{parab=a3a=-2k=3l}, \ref{parab=a3a=-2k=3l+2}, \ref{parab=a4a=k-2k=4l}, \ref{parab=a4a=k-2k=4l+1}, \ref{parab=a4a=k-2k=4l+2}, \ref{parab=a4a=k-2k=4l+3}, \ref{parab=a4a=k-3k=4l}, \ref{parab=a4a=k-3k=4l+1}, \ref{parab=a4a=k-3k=4l+2}, and in \ref{parab=a4a=k-3k=4l+3}.}\label{Ta:fig:red65}
\end{center}
\end{table}

\bigbreak
\bigbreak

\begin{figure}[!ht]
    \psfrag{5l}{\huge $5l$}
    \psfrag{4l}{\huge $4l$}
    \psfrag{3l}{\huge $3l$}
    \psfrag{6l}{\huge $6l$}
    \psfrag{3l+1}{\huge $3l+1$}
    \psfrag{l}{\huge $l$}
    \psfrag{2l+1}{\huge $2l+1$}
    \psfrag{4l+1}{\huge $4l+1$}
    \psfrag{5l+1}{\huge $5l+1$}
    \psfrag{3l+1}{\huge $3l+1$}
    \psfrag{6l+1}{\huge $\mathbf{6l+1}$}
    \centerline{\scalebox{.37}{\includegraphics{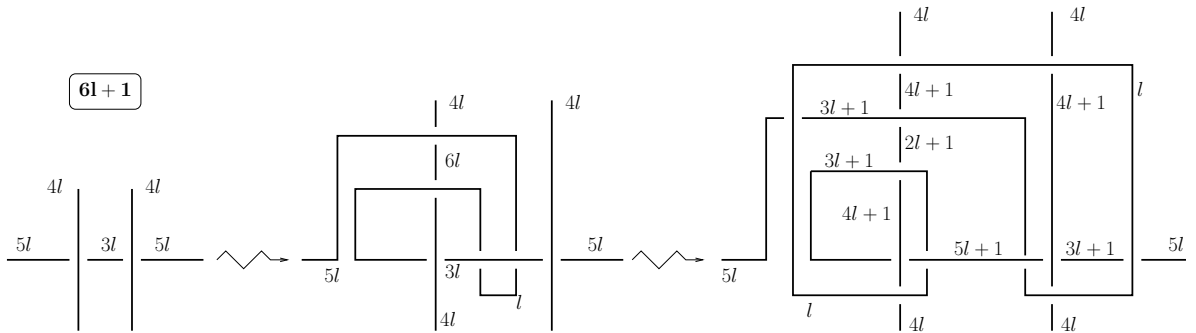}}}
    \caption{Removing color $k$ from an under-arc:  $b-a$ instance with $3a+1=-1$; and $k=3l$, $a=4l$ mod $6l+1$. $(k=3l<)a=4l(<6l-1<6l)$.}\label{fig:red116bis}
\end{figure}

\paragraph{The $b=a$ instance with $3a+1=-1$; and $k=3l$, $a=4l$ mod $6l+1$.}\label{parab=a3a=-2k=3l} View Figure \ref{fig:red116bis}.

\bigbreak

\bigbreak

\bigbreak

\bigbreak

\bigbreak

\begin{table}[h!]
\begin{center}
    \begin{tabular}{ | c |  }\hline
$l\geq 2 \qquad 2l+1 \leq 3l < 3l+1 < 4l+1 < 5l+1 \leq 6l-1 < 6l $\\ \hline
    \end{tabular}
\caption{Consider Figure \ref{fig:red116bis}. $l\geq 2$. The equalities $2l+1=3l$, $5l+1=6l-1$, and $5l+1=6l$  give rise to further considerations for $l=2$ alone. This situation is dealt with in Figure \ref{fig:red1161}.}\label{Ta:fig:red116}
\end{center}
\end{table}

\bigbreak

\bigbreak

\bigbreak

\begin{figure}[!ht]
    \psfrag{1}{\huge $1$}
    \psfrag{-1}{\huge $-1$}
    \psfrag{3}{\huge $3$}
    \psfrag{6}{\huge $6$}
    \psfrag{7}{\huge $7$}
    \psfrag{10}{\huge $10$}
    \psfrag{8}{\huge $8$}
    \psfrag{5}{\huge $5$}
    \psfrag{13}{\huge $\mathbf{13}$}
    \centerline{\scalebox{.32}{\includegraphics{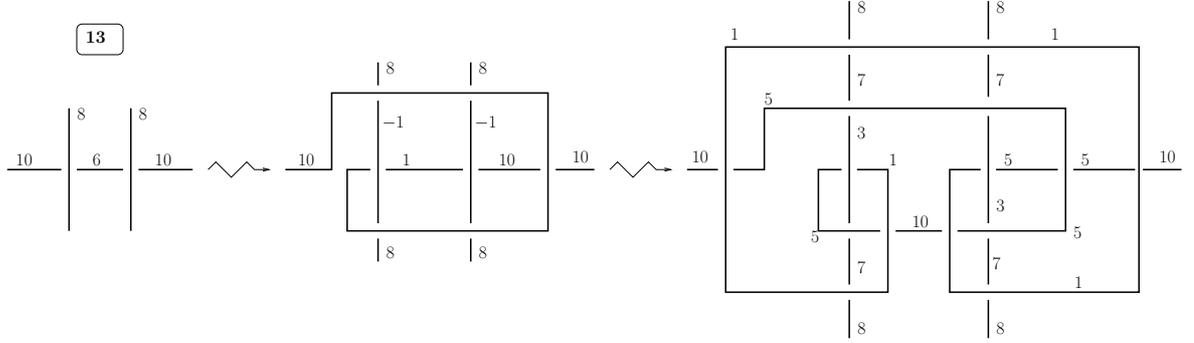}}}
    \caption{Removing color $k$ from an under-arc:   $b-a$ instance with $3a+1=-1$; and $k=3l$, $a=4l$ mod $6l+1$ and $l=2$.}\label{fig:red1161}
\end{figure}

\bigbreak

\bigbreak

\begin{figure}[!ht]
    \psfrag{-1}{\huge $-1$}
    \psfrag{5l+3}{\huge $5l+3$}
    \psfrag{4l+2}{\huge $4l+2$}
    \psfrag{3l+2}{\huge $3l+2$}
    \psfrag{2l}{\huge $2l$}
    \psfrag{3l+1}{\huge $3l+1$}
    \psfrag{l}{\huge $l$}
    \psfrag{2l+1}{\huge $2l+1$}
    \psfrag{4l+1}{\huge $4l+1$}
    \psfrag{5l+1}{\huge $5l+1$}
    \psfrag{3l+1}{\huge $3l+1$}
    \psfrag{6l+5}{\huge $\mathbf{6l+5}$}
    \centerline{\scalebox{.32}{\includegraphics{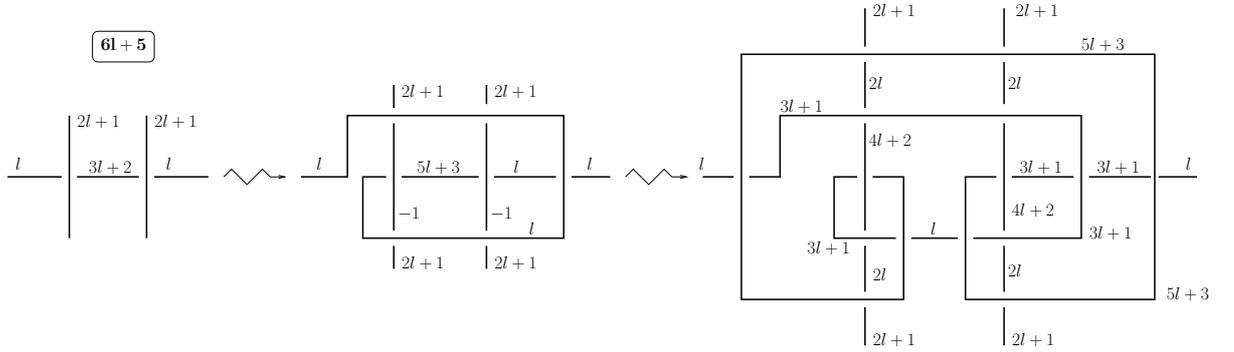}}}
    \caption{Removing color $k$ from an under-arc: the $b=a$ instance with $3a+1=-1$; and $k=3l+2$, $a=2l+1$ mod $6l+5$. $a=2l+1(<k=3l+2<6l+3<6l+4)$.}\label{fig:red111}
\end{figure}

\paragraph{The  $b-a$ instance with $3a+1=-1$; and $k=3l+2$, $a=2l+1$ mod $6l+5$ $a=2l+1(<k=3l<6l+3<6l+4)$.}\label{parab=a3a=-2k=3l+2} View Figure \ref{fig:red111}.

\bigbreak

\bigbreak

\bigbreak

\begin{table}[h!]
\begin{center}
    \begin{tabular}{ | c |   }\hline
$l\geq 1 \qquad 2l<3l+1<k(=3l+2)<4l+2<5l+3<6l+3<6l+4$ \\ \hline
    \end{tabular}
\caption{Consider Figure \ref{fig:red111}. $l\geq 2$. There are no further considerations, here.}\label{Ta:fig:red111}
\end{center}
\end{table}

\bigbreak

\bigbreak

\begin{figure}[!ht]
    \psfrag{5l}{\huge $5l$}
    \psfrag{4l}{\huge $4l$}
    \psfrag{2l}{\huge $2l$}
    \psfrag{6l}{\huge $6l$}
    \psfrag{7l}{\huge $7l$}
    \psfrag{0}{\huge $0$}
    \psfrag{2l+1}{\huge $2l+1$}
    \psfrag{7l+1}{\huge $7l+1$}
    \psfrag{5l+1}{\huge $5l+1$}
    \psfrag{6l+1}{\huge $6l+1$}
    \psfrag{8l+1}{\huge $\mathbf{8l+1}$}
    \centerline{\scalebox{.35}{\includegraphics{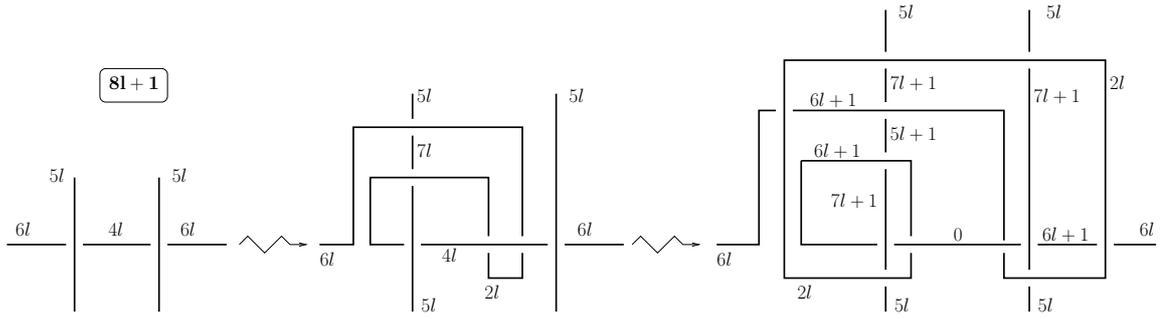}}}
    \caption{Removing color $k$ from an under-arc: the $b=a$ instance with $4a+k+2=-1$; and $k=4l$, $a=5l$ mod $8l+1$. $(k=4l<)a=5l(<8l-1<8l)$.}\label{fig:red126}
\end{figure}

\paragraph{The $b=a$ instance with $4a+k+2=-1$; and $k=4l$, $a=5l$ mod $8l+1$. so $(k=4l<)a=5l(<8l-1<8l)$.}\label{parab=a4a=k-2k=4l} View Figure \ref{fig:red126}.

\bigbreak

\bigbreak

\bigbreak

\bigbreak

\begin{table}[h!]
\begin{center}
    \begin{tabular}{ | c |}\hline
$l\geq 2 \qquad 2l<k(=4l)<5l+1<6l+1<7l+1\leq 8l-1<8l$\\ \hline
    \end{tabular}
\caption{Consider Figure \ref{fig:red126}. $l\geq 2$. The equalities $7l+1=8l-1$, and $7l+1=8l$ give rise to further considerations for $l=2$ alone. This situation is dealt with in Figures \ref{fig:red45a} and \ref{fig:red45b}.}\label{Ta:fig:red126}
\end{center}
\end{table}

\bigbreak

\bigbreak

\begin{figure}[!ht]
    \psfrag{10}{\huge $10$}
    \psfrag{12}{\huge $12$}
    \psfrag{8}{\huge $8$}
    \psfrag{4}{\huge $4$}
    \psfrag{14}{\huge $14$}
    \psfrag{17}{\huge $\mathbf{17}$}
    \centerline{\scalebox{.35}{\includegraphics{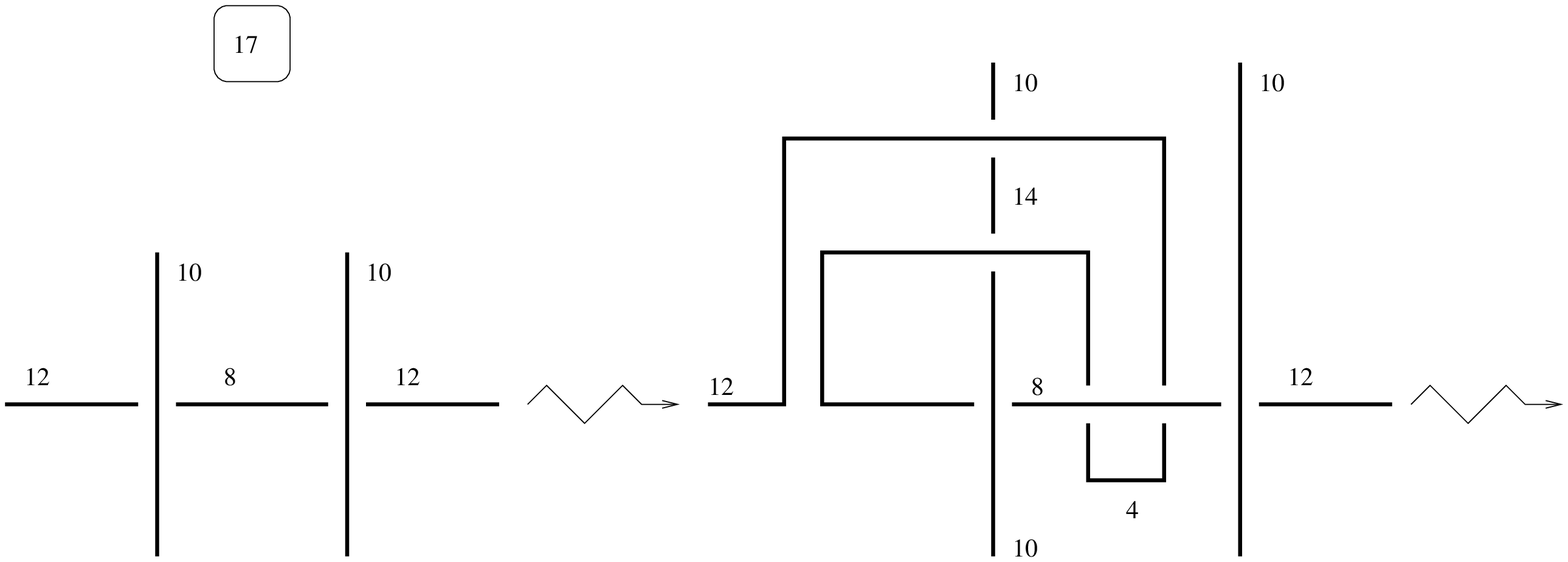}}}
    \caption{Removing color $k$ from an under-arc: the $b=a$ instance with $4a+k+2=-1$; and $k=4l$, $a=5l$ mod $8l+1$. Fixing the $l=2$ case (part $I$).}\label{fig:red45a}
\end{figure}

\bigbreak

\bigbreak

\begin{figure}[!ht]
    \psfrag{0}{\huge $0$}
    \psfrag{2}{\huge $2$}
    \psfrag{-2}{\huge $-2$}
    \psfrag{10}{\huge $10$}
    \psfrag{12}{\huge $12$}
    \psfrag{7}{\huge $7$}
    \psfrag{4}{\huge $4$}
    \psfrag{6}{\huge $6$}
    \psfrag{11}{\huge $11$}
    \psfrag{5}{\huge $5$}
    \psfrag{13}{\huge $13$}
    \psfrag{17}{\huge $\mathbf{17}$}
    \centerline{\scalebox{.37}{\includegraphics{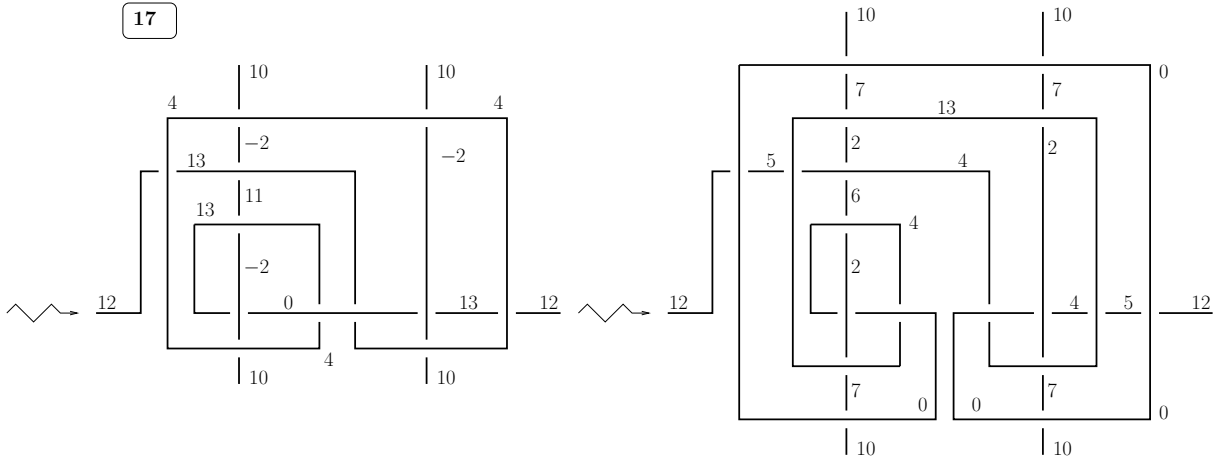}}}
    \caption{Removing color $k$ from an under-arc: the $b=a$ instance with $4a+k+2=-1$; and $k=4l$, $a=5l$ mod $8l+1$. Fixing the $l=2$ case (part $II$).}\label{fig:red45b}
\end{figure}

\bigbreak

\bigbreak

\begin{figure}[!ht]
    \psfrag{0}{\huge $0$}
    \psfrag{5l+1}{\huge $5l+1$}
    \psfrag{4l+1}{\huge $4l+1$}
    \psfrag{2l}{\huge $2l$}
    \psfrag{2l+1}{\huge $2l+1$}
    \psfrag{5l+2}{\huge $5l+2$}
    \psfrag{7l+3}{\huge $7l+3$}
    \psfrag{7l+2}{\huge $7l+2$}
    \psfrag{6l+2}{\huge $6l+2$}
    \psfrag{8l+3}{\huge $\mathbf{8l+3}$}
    \centerline{\scalebox{.35}{\includegraphics{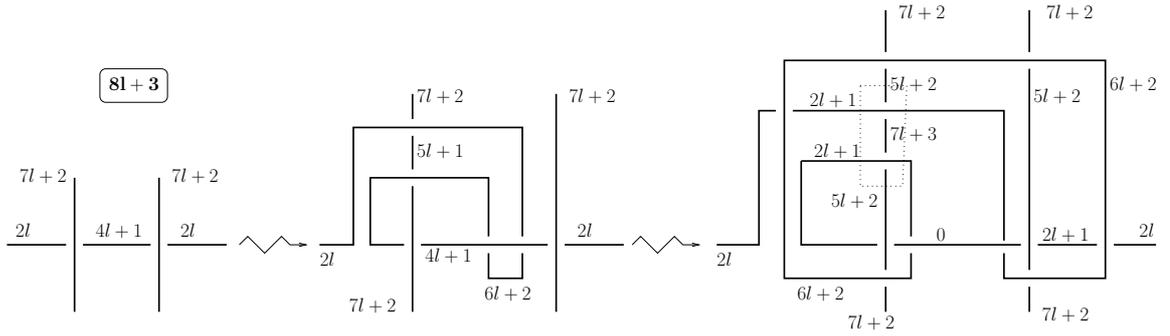}}}
    \caption{Removing color $k$ from an under-arc: the $b=a$ instance with $4a+k+2=-1$; and $k=4l+1$, $a=7l+2$ mod $8l+3$. $(k=4l+1<)a=7l+2(<8l+1)$. The boxed region is fixed for $l=2$ in Figure \ref{fig:red133a}.}\label{fig:red133}
\end{figure}

\paragraph{The $b=a$ instance with $4a+k+2=-1$; and $k=4l+1$, $a=7l+2$ mod $8l+3$. $(k=4l+1<)a=7l+2(<8l+1)$.}\label{parab=a4a=k-2k=4l+1} View Figure \ref{fig:red133}.

\bigbreak

\bigbreak

\bigbreak

\bigbreak

\begin{table}[h!]
\begin{center}
    \begin{tabular}{ | c | }\hline
$l\geq 1 \qquad 2l+1 < k(=4l+1) < 5l+2 < 6l+2 < 7l+3 \leq 8l+1 < 8l+2$\\ \hline
    \end{tabular}
\caption{Consider Figure \ref{fig:red133}. $l\geq 1$. The equalities $7l+3=8l+1$, and $7l+3=8l+2$ give rise to further considerations for $l=2$ alone. This situation is dealt with in Figure \ref{fig:red133a}.}\label{Ta:fig:red133}
\end{center}
\end{table}

\bigbreak

\bigbreak

\begin{figure}[!ht]
    \psfrag{0}{\huge $0$}
    \psfrag{12}{\huge $12$}
    \psfrag{7}{\huge $7$}
    \psfrag{5}{\huge $5$}
    \psfrag{-2}{\huge $-2$}
    \psfrag{19}{\huge $\mathbf{19}$}
    \centerline{\scalebox{.35}{\includegraphics{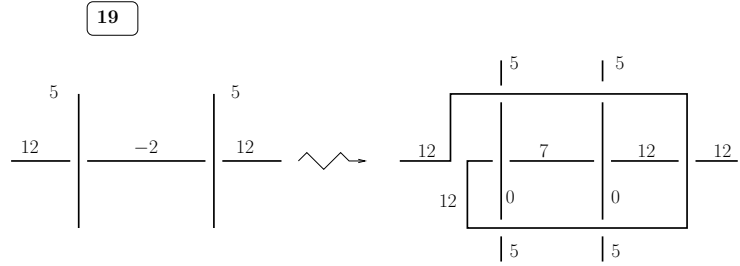}}}
    \caption{Removing color $k$ from an under-arc: the $b-a$ instance with $4a=k-2$ and $k=4l+1$ so $(k=4l+1<)a=7l+2(<8l+1)$. Fixing the boxed region for $l=2$ in Figure \ref{fig:red133}.}\label{fig:red133a}
\end{figure}

\bigbreak

\bigbreak

\begin{figure}[!ht]
    \psfrag{0}{\huge $0$}
    \psfrag{6l+3}{\huge $6l+3$}
    \psfrag{l}{\huge $l$}
    \psfrag{4l+2}{\huge $4l+2$}
    \psfrag{3l+1}{\huge $3l+1$}
    \psfrag{2l+1}{\huge $2l+1$}
    \psfrag{6l+4}{\huge $6l+4$}
    \psfrag{l+1}{\huge $l+1$}
    \psfrag{3l+2}{\huge $3l+2$}
    \psfrag{8l+5}{\huge $\mathbf{8l+5}$}
    \centerline{\scalebox{.36}{\includegraphics{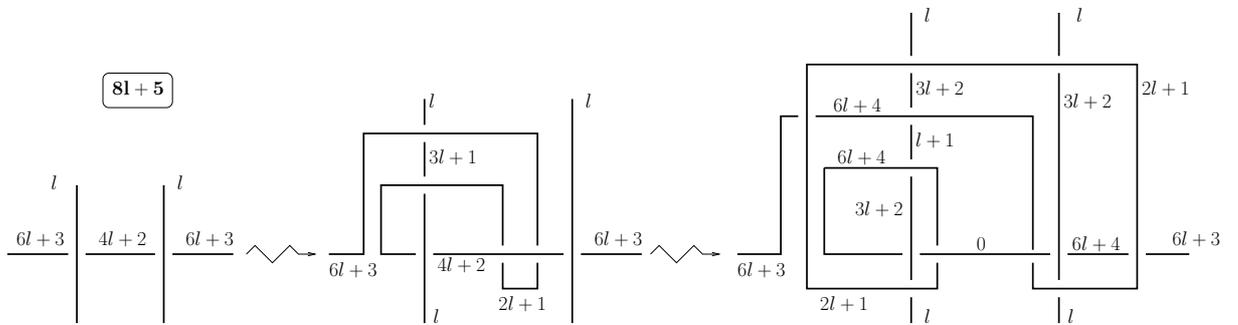}}}
    \caption{Removing color $k$ from an under-arc: the $b=a$ instance with $4a+k+2=-1$; and $k=4l+2$, $a=l$ mod $8l+5$. $a=l(<k=4l+2<8l+3<8l+4)$.}\label{fig:red134}
\end{figure}

\paragraph{The $b=a$ instance with $4a+k+2=-1$; and $k=4l+2$, $a=l$ mod $8l+5$. $a=l(<k=4l+2<8l+3<8l+4)$.}\label{parab=a4a=k-2k=4l+2} View Figure \ref{fig:red134}.

\bigbreak

\bigbreak

\begin{table}[h!]
\begin{center}
    \begin{tabular}{ | c | }\hline
$l\geq 1 \qquad l+1 < 2l+1 < 3l+2 < k(=4l+2) < 6l+4 < 8l+3 <  8l+4$\\ \hline
    \end{tabular}
\caption{Consider Figure \ref{fig:red134}. $l\geq 1$. There are no further considerations here.}\label{Ta:fig:red134}
\end{center}
\end{table}

\begin{figure}[!ht]
    \psfrag{0}{\huge $0$}
    \psfrag{4l+2}{\huge $4l+2$}
    \psfrag{4l+3}{\huge $4l+3$}
    \psfrag{2l+1}{\huge $2l+1$}
    \psfrag{6l+5}{\huge $6l+5$}
    \psfrag{2l+2}{\huge $2l+2$}
    \psfrag{l+1}{\huge $l+1$}
    \psfrag{l}{\huge $l$}
    \psfrag{3l+2}{\huge $3l+2$}
    \psfrag{3l+3}{\huge $3l+3$}
    \psfrag{8l+7}{\huge $\mathbf{8l+7}$}
    \centerline{\scalebox{.36}{\includegraphics{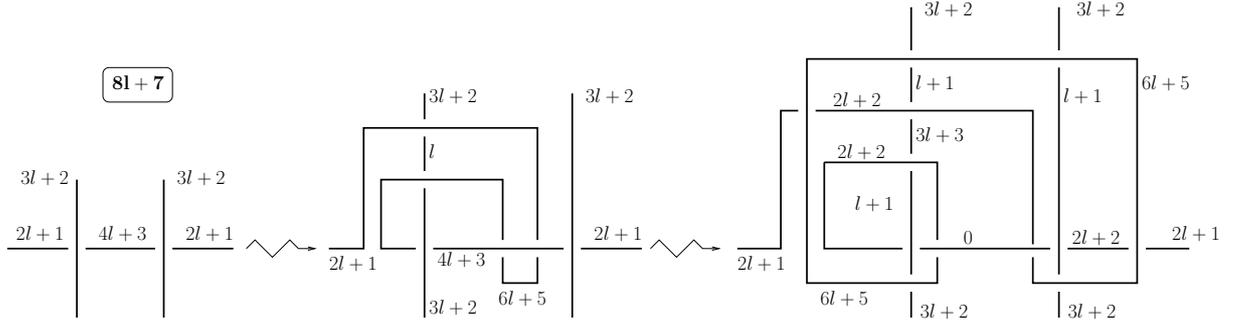}}}
    \caption{Removing color $k$ from an under-arc: the $b=a$ instance with $4a+k+2=-1$; and $k=4l+3$, $a=3l+2$ mod $8l+7$. $a=3l+2(<k=4l+3<8l+5<8l+6)$.}\label{fig:red138}
\end{figure}

\paragraph{The $b=a$ instance with $4a+k+2=-1$; and $k=4l+3$, $a=3l+2$ mod $8l+7$. $a=3l+2(<k=4l+3<8l+5<8l+6)$.}\label{parab=a4a=k-2k=4l+3} View Figure \ref{fig:red138}.

\bigbreak

\bigbreak

\begin{table}[h!]
\begin{center}
    \begin{tabular}{ | c | }\hline
$l\geq 2 \qquad l+1 < 2l+2 < 3l+3 < k(=4l+3) < 6l+5 < 8l+5 <  8l+6$\\ \hline
    \end{tabular}
\caption{Consider Figure \ref{fig:red138}. $l\geq 2$. There are no further considerations here.}\label{Ta:fig:red138}
\end{center}
\end{table}

\bigbreak

\bigbreak

\begin{figure}[!ht]
    \psfrag{5l-1}{\huge $5l-1$}
    \psfrag{4l}{\huge $4l$}
    \psfrag{2l-1}{\huge $2l-1$}
    \psfrag{6l}{\huge $6l$}
    \psfrag{7l}{\huge $7l$}
    \psfrag{0}{\huge $1$}
    \psfrag{2l+2}{\huge $2l+2$}
    \psfrag{5l+2}{\huge $5l+2$}
    \psfrag{7l+3}{\huge $7l+3$}
    \psfrag{5l+1}{\huge $5l+1$}
    \psfrag{6l+1}{\huge $6l+1$}
    \psfrag{8l+1}{\huge $\mathbf{8l+1}$}
    \centerline{\scalebox{.35}{\includegraphics{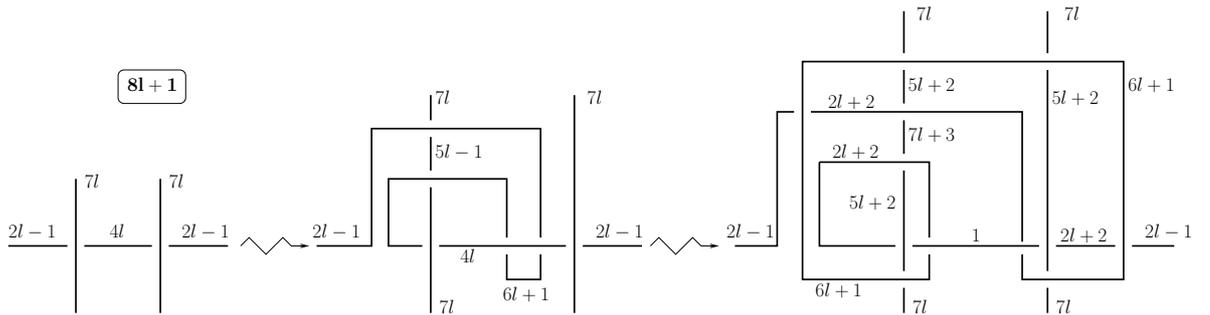}}}
    \caption{Removing color $k$ from an under-arc: the $b=a$ instance with $4a+k+2=-2$; and $k=4l$, $a=7l$ mod $8l+1$. $(k=4l<)a=7l(<8l)$.}\label{fig:red126a}
\end{figure}

\paragraph{The $b=a$ instance with $4a+k+2=-2$; and $k=4l$, $a=7l$ mod $8l+1$. $(k=4l<)a=7l(<8l)$.}\label{parab=a4a=k-3k=4l} View Figure \ref{fig:red126a}.

\bigbreak

\bigbreak

\bigbreak

\bigbreak

\begin{table}[h!]
\begin{center}
    \begin{tabular}{ | c |  }\hline
$l\geq 2 \qquad 2l+2\leq k(=4l) < 5l+2 \leq 6l+1 < 7l+3 \leq 8l-1 < 8l$\\ \hline
    \end{tabular}
\caption{Consider Figure \ref{fig:red126a}. $l\geq 2$. The equalities $2l+2= 4l$, $6l+1=8l-1$, $7l+3=8l-1$, and $7l+3=8l$ do not give rise to further considerations.}\label{Ta:fig:red126bis}
\end{center}
\end{table}

\bigbreak

\bigbreak

\begin{figure}[!ht]
    \psfrag{1}{\huge $1$}
    \psfrag{5l+1}{\huge $5l+1$}
    \psfrag{4l+1}{\huge $4l+1$}
    \psfrag{2l}{\huge $2l$}
    \psfrag{2l+1}{\huge $2l+1$}
    \psfrag{5l+4}{\huge $5l+4$}
    \psfrag{7l+3}{\huge $7l+3$}
    \psfrag{7l+4}{\huge $7l+4$}
    \psfrag{6l+4}{\huge $6l+4$}
    \psfrag{7l+1}{\huge $7l+1$}
    \psfrag{6l+1}{\huge $6l+1$}
    \psfrag{8l+3}{\huge $\mathbf{8l+3}$}
    \centerline{\scalebox{.35}{\includegraphics{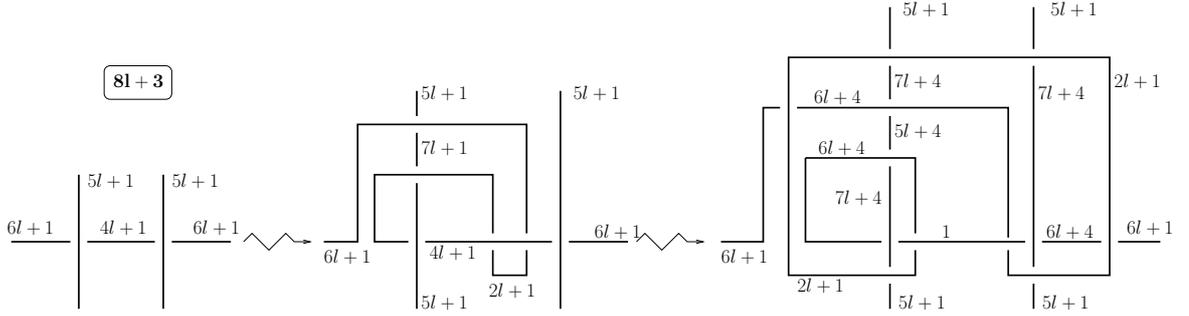}}}
    \caption{Removing color $k$ from an under-arc: the $b-a$ instance with $4a=k-2$ and $k=4l+1$ so $(k=4l+1<)a=5l+2(<8l+1<8l+2)$ (continued in Figure \ref{fig:red133beta}.}\label{fig:red133alpha}
\end{figure}

\paragraph{The $b=a$ instance with $4a+k+2=-2$; and $k=4l+1$, $a=5l+1$ mod $8l+3$. $(k=4l+1<)a=5l+1(<8l+1<8l+2)$.}\label{parab=a4a=k-3k=4l+1} View Figure \ref{fig:red133alpha}.

\bigbreak

\bigbreak

\bigbreak

\bigbreak

\begin{figure}[!ht]
    \psfrag{1}{\huge $1$}
    \psfrag{5l+1}{\huge $5l+1$}
    \psfrag{3l+1}{\huge $3l+1$}
    \psfrag{l+1}{\huge $l+1$}
    \psfrag{2l}{\huge $2l$}
    \psfrag{2l+1}{\huge $2l+1$}
    \psfrag{2l+4}{\huge $2l+4$}
    \psfrag{7l+3}{\huge $7l+3$}
    \psfrag{3l+4}{\huge $3l+4$}
    \psfrag{6l+4}{\huge $6l+4$}
    \psfrag{7l+1}{\huge $7l+1$}
    \psfrag{6l+1}{\huge $6l+1$}
    \psfrag{8l+3}{\huge $\mathbf{8l+3}$}
    \centerline{\scalebox{.35}{\includegraphics{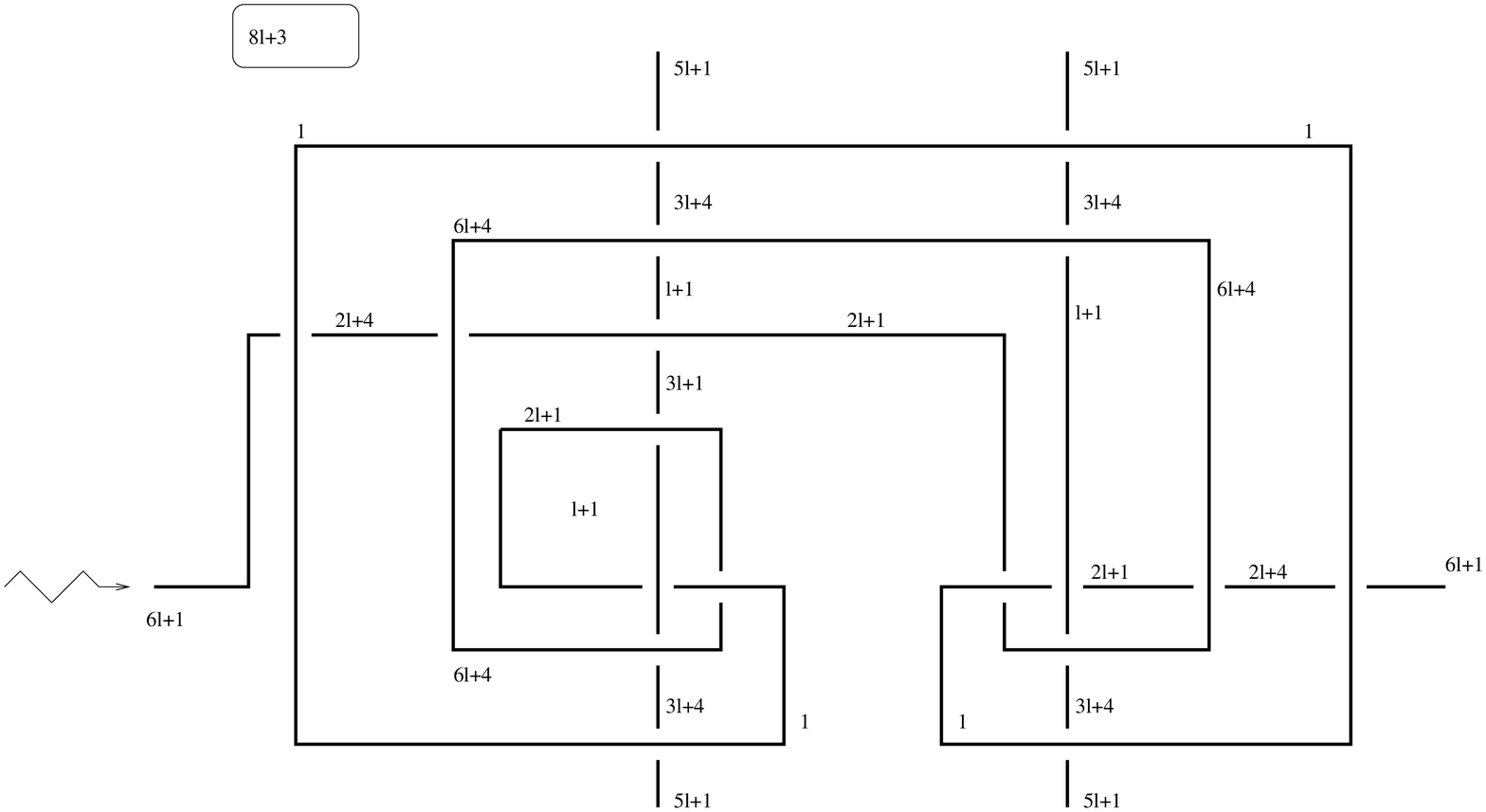}}}
    \caption{Removing color $k$ from an under-arc: the $b-a$ instance with $4a=k-2$ and $k=4l+1$ so $(k=4l+1<)a=5l+2(<8l+1<8l+2)$ (continued form Figure \ref{fig:red133alpha}.}\label{fig:red133beta}
\end{figure}

\bigbreak

\begin{table}[h!]
\begin{center}
    \begin{tabular}{ | c | }\hline
$l\geq 1 \qquad l+1 < 2l+1 < 2l+4 < 3l+4  \leq k(=4l+1) < 5l+1 < 6l+4 \leq 8l+1 < 8l+2 $  \\ \hline
    \end{tabular}
\caption{Consider Figure \ref{fig:red133beta}. $l\geq 2$. The equalities $2l+4= 4l+1$, $3l+4=4l+1$, $6l+4=8l+1$, and $6l+4=8l+2$ do not give rise to further considerations.}\label{Ta:fig:red133bis}
\end{center}
\end{table}

\bigbreak

\begin{figure}[!ht]
    \psfrag{3}{\huge $3$}
    \psfrag{5}{\huge $5$}
    \psfrag{6}{\huge $6$}
    \psfrag{7}{\huge $7$}
    \psfrag{8}{\huge $8$}
    \psfrag{11}{\huge $\mathbf{11}$}
    \centerline{\scalebox{.5}{\includegraphics{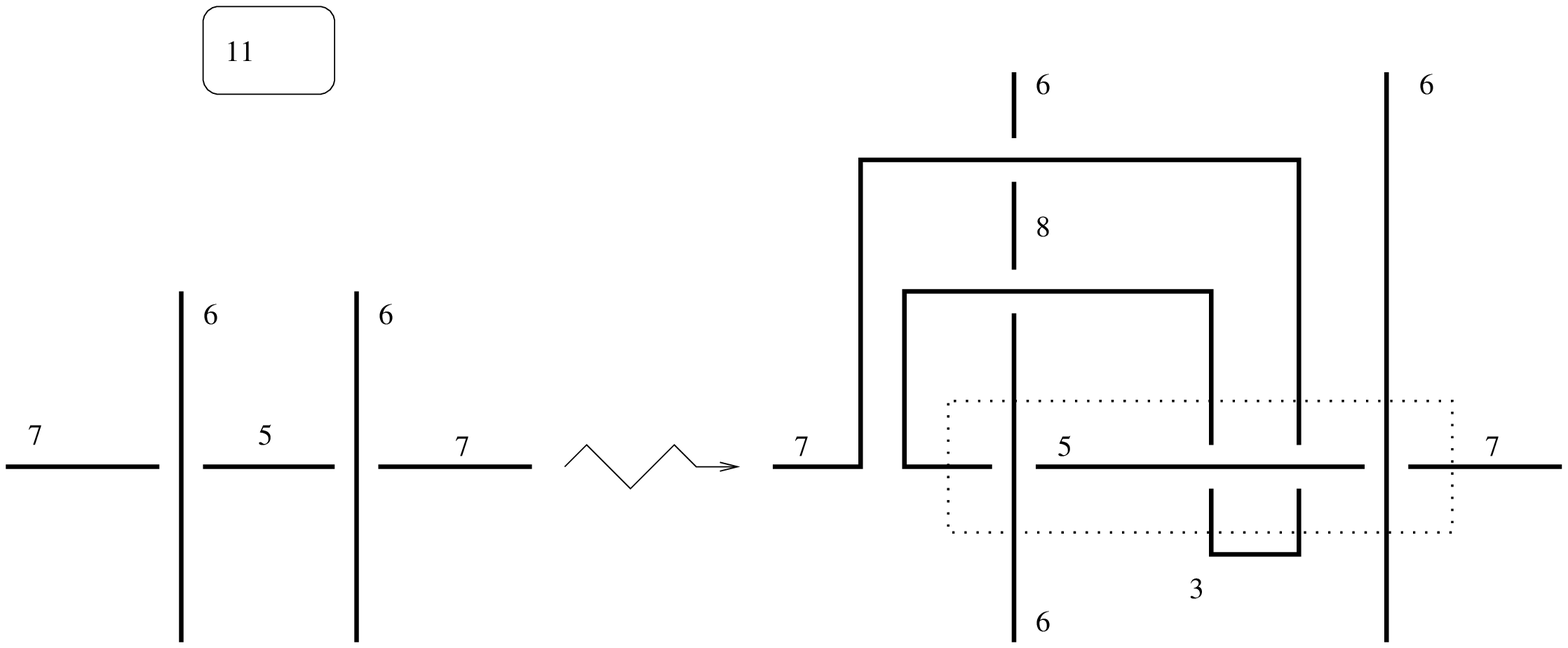}}}
    \caption{Fixing the $l=1$ situation in Figure \ref{fig:red133beta}.}\label{fig:red133betaalpha}
\end{figure}

\bigbreak

\begin{figure}[!ht]
    \psfrag{3}{\huge $3$}
    \psfrag{0}{\huge $0$}
    \psfrag{6}{\huge $6$}
    \psfrag{7}{\huge $7$}
    \psfrag{1}{\huge $1$}
    \psfrag{-1}{\huge $-1$}
    \psfrag{11}{\huge $\mathbf{11}$}
    \centerline{\scalebox{.5}{\includegraphics{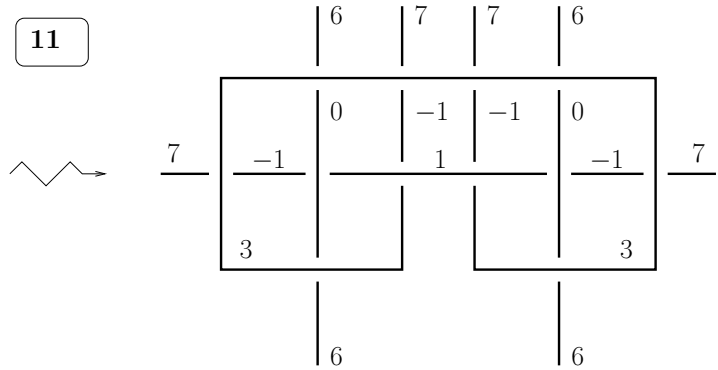}}}
    \caption{Fixing the $l=1$ situation in Figure \ref{fig:red133beta} (cont'd).}\label{fig:red133betabeta}
\end{figure}

\bigbreak

\begin{figure}[!ht]
    \psfrag{4}{\huge $4$}
    \psfrag{2}{\huge $2$}
    \psfrag{3}{\huge $3$}
    \psfrag{0}{\huge $0$}
    \psfrag{6}{\huge $6$}
    \psfrag{7}{\huge $7$}
    \psfrag{1}{\huge $1$}
    \psfrag{-1}{\huge $-1$}
    \psfrag{11}{\huge $\mathbf{11}$}
    \centerline{\scalebox{.5}{\includegraphics{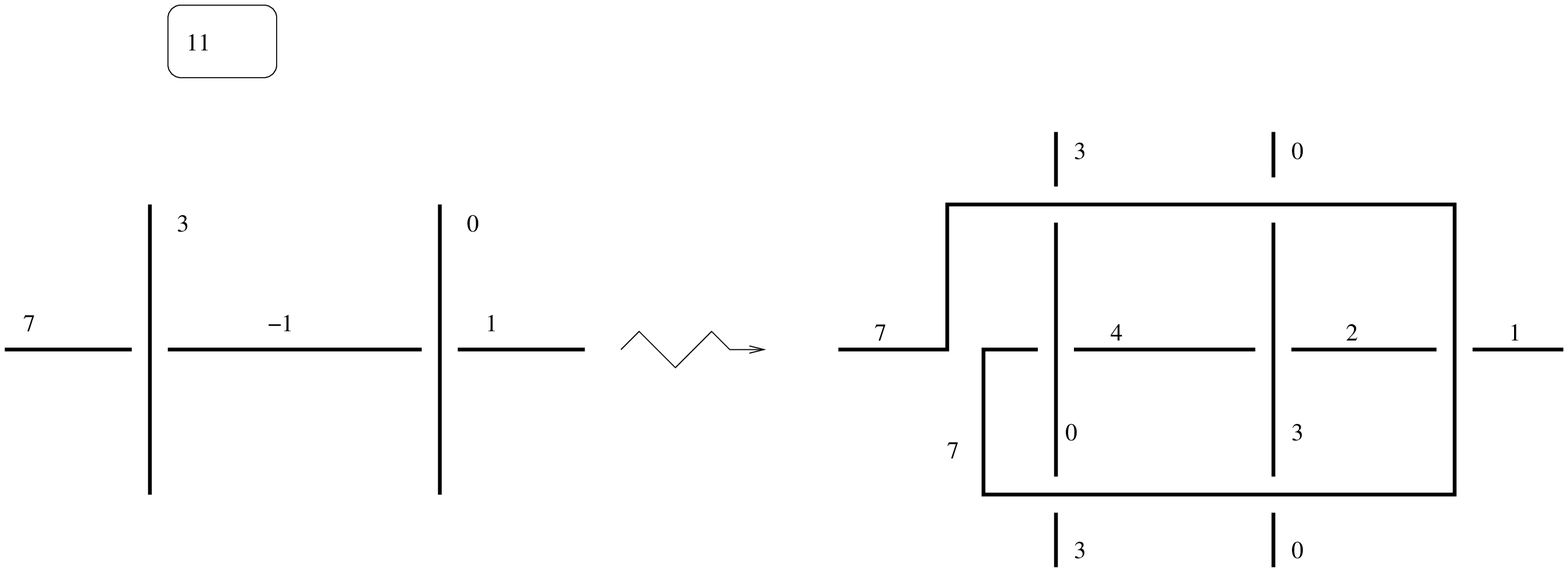}}}
    \caption{Fixing the $l=1$ situation in Figure \ref{fig:red133beta}.}\label{fig:red133betabetabeta}
\end{figure}

\begin{figure}[!ht]
    \psfrag{4}{\huge $4$}
    \psfrag{0}{\huge $0$}
    \psfrag{2}{\huge $2$}
    \psfrag{1}{\huge $1$}
    \psfrag{3}{\huge $3$}
    \psfrag{7}{\huge $7$}
    \psfrag{-1}{\huge $-1$}
    \psfrag{11}{\huge $\mathbf{11}$}
    \centerline{\scalebox{.5}{\includegraphics{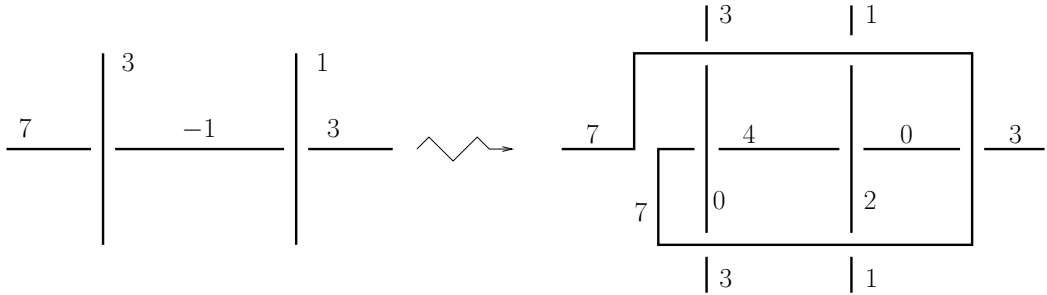}}}
    \caption{Fixing the $l=1$ situation in Figure \ref{fig:red133beta}.}\label{fig:red133betabetabetabis}
\end{figure}

\begin{figure}[!ht]
    \psfrag{1}{\huge $1$}
    \psfrag{2l+3}{\huge $2l+3$}
    \psfrag{l}{\huge $l$}
    \psfrag{4l+2}{\huge $4l+2$}
    \psfrag{3l+1}{\huge $3l+1$}
    \psfrag{2l}{\huge $2l$}
    \psfrag{6l+4}{\huge $6l+4$}
    \psfrag{l-1}{\huge $l-1$}
    \psfrag{3l+2}{\huge $3l+2$}
    \psfrag{l+2}{\huge $l+2$}
    \psfrag{3l+4}{\huge $3l+4$}
    \psfrag{8l+5}{\huge $\mathbf{8l+5}$}
    \centerline{\scalebox{.36}{\includegraphics{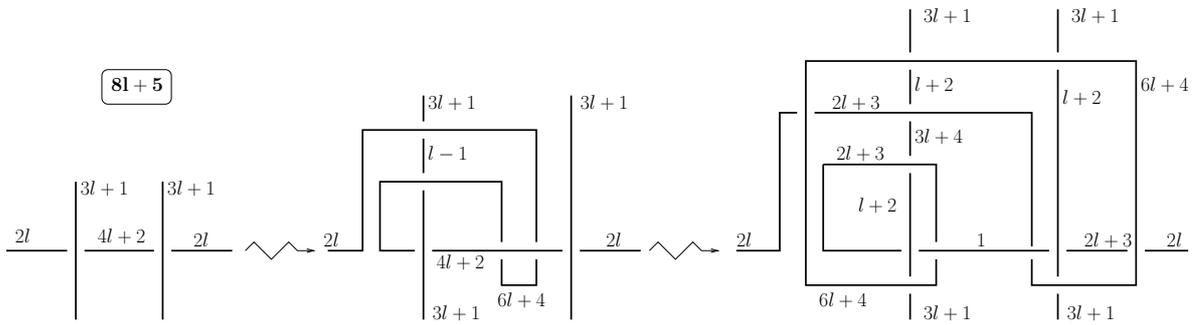}}}
    \caption{Removing color $k$ from an under-arc: the $b=a$ instance with $4a+k+2=-2$; and $k=4l+2$, $a=3l+1$ mod $8l+5$. $a=3l+1(<k=4l+2<8l+3<8l+4)$.}\label{fig:red134alpha}
\end{figure}

\paragraph{The $b=a$ instance with $4a+k+2=-2$; and $k=4l+2$, $a=3l+1$ mod $8l+5$. $a=3l+1(<k=4l+2<8l+3<8l+4)$.}\label{parab=a4a=k-3k=4l+2} View Figure \ref{fig:red134alpha}.

\begin{table}[h!]
\begin{center}
    \begin{tabular}{ | c | }\hline
$l\geq 1 \qquad l+2 < 2l+3 < 3l+4 \leq 4l+2 < 8l+3 < 8l+4$ \\ \hline
    \end{tabular}
\caption{Consider Figure \ref{fig:red134alpha}. $l\geq 1$. The equality $3l+4= 4l+2$ does not give rise to further considerations.}\label{Ta:fig:red134alpha}
\end{center}
\end{table}

\begin{figure}[!ht]
    \psfrag{1}{\huge $1$}
    \psfrag{3l+4}{\huge $3l+4$}
    \psfrag{4l+3}{\huge $4l+3$}
    \psfrag{2l+1}{\huge $2l+1$}
    \psfrag{6l+4}{\huge $6l+4$}
    \psfrag{2l+2}{\huge $2l+2$}
    \psfrag{6l+7}{\huge $6l+7$}
    \psfrag{l}{\huge $l$}
    \psfrag{3l+1}{\huge $3l+1$}
    \psfrag{l+3}{\huge $l+3$}
    \psfrag{8l+7}{\huge $\mathbf{8l+7}$}
    \centerline{\scalebox{.36}{\includegraphics{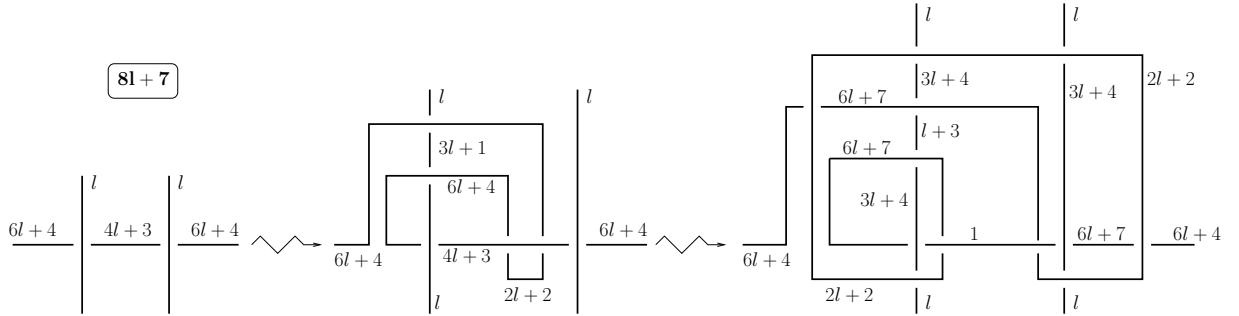}}}
    \caption{Removing color $k$ from an under-arc: the $b=a$ instance with $4a+k+2=-2$; and $k=4l+2$, $a=l$ mod $8l+7$. $a=3l+2(<k=4l+3<8l+5<8l+6)$.}\label{fig:red138alpha}
\end{figure}

\paragraph{The $b=a$ instance with $4a+k+2=-2$; and $k=4l+2$, $a=l$ mod $8l+7$. $a=l(<k=4l+3<8l+5<8l+6)$.}\label{parab=a4a=k-3k=4l+3} View Figure \ref{fig:red138alpha}.

\begin{table}[h!]
\begin{center}
    \begin{tabular}{ | c |  }\hline
$l\geq 2 \qquad l+3 \leq 2l+2 < 3l+4 \leq 4l+3 < 6l+7 \leq 8l+5 < 8l+6$  \\ \hline
    \end{tabular}
\caption{Consider Figure \ref{fig:red138alpha}. $l\geq 2$. The equalities $3l+4= 4l+3$, $6l+7= 8l+5$, and $6l+7= 8l+6$ do not give rise to further considerations.}\label{Ta:fig:red138alpha}
\end{center}
\end{table}

\end{document}